\documentclass[12pt]{article}

\usepackage{amsmath,amssymb,mathtools,amsthm}
\mathtoolsset{showonlyrefs} % blends out all labels of equations that are not referenced
\usepackage{authblk}
\usepackage{yhmath}
\usepackage{dsfont}
\usepackage{faktor}
\usepackage{graphicx}
\usepackage{subcaption}
\usepackage{enumerate}
\usepackage[dvipsnames]{xcolor}
\usepackage{color}
\usepackage{tikz}
\usepackage{tikz-cd}
\usepackage{overpic}
\usepackage[utf8]{inputenc}
\usepackage[
    inner=2.3cm,
    outer=2.3cm,
    top=2.7cm,
    marginparwidth=2.5cm,
    marginparsep=0.1cm
    ]
    {geometry}
\usepackage[pdffitwindow=false,
            plainpages=false,
            pdfpagelabels=true,
            pdfpagemode=UseOutlines,
            pdfpagelayout=SinglePage,
            bookmarks=false,
            colorlinks=true,
            hyperfootnotes=false,
            linkcolor=blue,
            urlcolor=blue!30!black,
            citecolor=green!50!black]{hyperref}
    
\renewcommand{\H}{\mathcal{H}}
\newcommand{\R}{\mathds{R}}
\newcommand{\N}{\mathds{N}}
\newcommand{\C}{\mathds{C}}
\newcommand{\Z}{\mathds{Z}}

\renewcommand{\L}{\mathcal{L}}

\newcommand{\norm}[1]{\left\lVert#1\right\rVert}
\newcommand{\abs}[1]{\left|#1\right|}
\newcommand{\1}{\mathds{1}}

\newcommand{\G}{\mathbf{G}}
\newcommand{\f}{\mathbf{f}}

\renewcommand{\P}{\mathds{P}}
\newcommand{\tor}{\mathds{T}}
\newcommand{\Per}{\mathcal{P}}
\newcommand{\ran}{\textup{ran}}
\newcommand{\Id}{\textup{Id}}

\newtheorem{theorem}{Theorem}[section]
\newtheorem{proposition}[theorem]{Proposition}
\newtheorem{lemma}[theorem]{Lemma}

\newtheorem{definition}[theorem]{Definition}

\theoremstyle{definition}

\theoremstyle{definition}
\newtheorem{remark}[theorem]{Remark}

\newenvironment{nalign}{
    \begin{equation}
    \begin{aligned}
}{
    \end{aligned}
    \end{equation}
    \ignorespacesafterend
}

\title{Extracting coherent sets in aperiodically driven flows from generators of Mather semigroups}

\author[1]{Robin Chemnitz}
\author[1,2]{Maximilian Engel}
\author[3]{P\'eter Koltai}

\affil[1]{Freie Universität Berlin}
\affil[2]{University of Amsterdam}
\affil[3]{Universität Bayreuth}

\date{}

\begin{document}

\maketitle

\begin{abstract}
    Coherent sets are time-dependent regions in the physical space of nonautonomous flows that exhibit little mixing with their neighborhoods, robustly under small random perturbations of the flow. They thus characterize the global long-term transport behavior of the system. We propose a framework to extract such time-dependent families of coherent sets for nonautonomous systems with an ergodic driving dynamics and (small) Brownian noise in physical space. Our construction involves the assembly and analysis of an operator on functions over the augmented space of the associated skew product that, for each fixed state of the driving, propagates distributions on the corresponding physical-space fibre according to the dynamics. This time-dependent operator has the structure of a semigroup (it is called the Mather semigroup), and we show that a spectral analysis of its generator allows for a trajectory-free computation of coherent families, simultaneously for all states of the driving. Additionally, for quasi-periodically driven torus flows, we propose a tailored Fourier discretization scheme for this generator and demonstrate our method by means of three examples of two-dimensional flows.
\end{abstract}

\section{Introduction}

Understanding the transport and mixing behavior of complicated nonautonomous flows is still an outstanding problem on the interface of dynamical systems theory and its many applications. We will consider in this work the problem of finding time-evolving families of sets in the state space of the dynamics that mix little with their surroundings. This way we provide a macroscopic view on prevalent transport-related characteristics of aperiodically driven flows.
In the classical theory of deterministic dynamical systems, invariant manifolds organize the state space in terms of transport: 
they form impenetrable barriers to purely advective transport. For nonautonomous systems, this led to the concept of lobe dynamics~\cite{KW90}, of general transport barriers such as Lagrangian Coherent Structures~\cite{haller2000lagrangian, shadden2005definition, rypina2007lagrangian}, and related objects~\cite{haller2013coherent, ma2014differential, haller2016defining, balasuriya2018generalized}.

\subsection{Background}

\paragraph{Small noise, metastability, and coherence.}
In a physically motivated modeling step, small random perturbations can be added to the (otherwise deterministic) state evolution, while the questions remain: What are the characteristic state-space objects governing mixing and the typical timescales associated to them?
In a natural approach to this question, subsets of the state space can be sought that mix particularly little with their surroundings. For autonomous systems, such sets are called \emph{almost-invariant}~\cite{dellnitz1999approximation} or \emph{metastable}~\cite{davies1982metastableI,davies1982metastableII,schutte2013metastability,bovier2016metastability}. Their presence implies the existence of persistent patterns, such as slowly decaying concentration fields of quantities passively transported by the flow~\cite{liu2004strange}. The corresponding objects in nonautonomous systems are called \emph{coherent sets}, and this term loosely refers to time-varying sets that have little dynamical interaction with their neighborhood. The concept has initially been developed in a finite-time setting for systems subject to noise~\cite{FSM10,froyland2013analytic,FPG14}; it has later been generalized separately to deterministic dynamics~\cite{froyland2015isoperimetry} and to the periodically-forced infinite-time setting~\cite{froyland2017estimating}. We note that the deterministic and noisy notions can be linked in the vanishing-diffusion limit~\cite{froyland2015isoperimetry,karrasch2020geometric,schilling2021heat}.

\paragraph{Escape rates.}
If finite-time considerations are replaced by infinite-time ones, the notions of metastability and coherence need to be adapted as well. \emph{Escape rates} become a natural measure of dynamical persistence of a set, and the notion has been extensively studied for autonomous \emph{open systems}~\cite{pianigiani1979expanding,chernov1998conditionally,liverani2003lasota,demers2005escape, DemersTodd2017, DemersWrightYoung2012,froyland2010escape,bunimovich2011place,bose2014ulam} and for time-homogeneous (Markovian) systems alike~\cite{collet1997pianigiani,froyland2013estimating} 
The nonautonomous generalization of escape rates from a family of time-dependent sets was considered in~\cite{balasuriya2014nonautonomous} for deterministic flows in terms of time-varying flux. It was used to extract coherent sets in~\cite{froyland2017estimating,froyland2020computation} for specific driving dynamics, and has appeared in the random dynamical systems literature for open systems~\cite{atnip2023equilibrium} 
to measure persistence.

\paragraph{Practical characterization and computation: Transfer operators.}
In computational approaches to metastability and coherence a functional view is prevailing. Sets are relaxed to (signed) distributions and these in turn are evolved by \emph{transfer operators} associated with the dynamics. 
Metastable sets of autonomous systems can be connected to the dominant eigenspectrum of transfer operators~\cite{davies1982metastableI,davies1982metastableII,dellnitz1999approximation,deuflhard2000identification,bovier2002metastability,deuflhard2005robust,huisinga2006metastability,schutte2013metastability}, while finite-time coherent sets can be extracted from a singular value decomposition of associated Markov operators~\cite{FSM10,froyland2013analytic,FPG14,denner2016computing}.
Intuitively, the eigen-(singular) values here are commensurate with the metastability (coherence) of the associated sets: The closer these values are to~1, the less dynamical leakage do the corresponding metastable (coherent) sets exhibit.

The connection to escape rates arises if we consider instead of one-step quantities (such as eigen- or singular values) of transfer operators an iterated quantity: \emph{decay rates}. Naturally, for autonomous systems the decay of a function is governed by its decomposition into eigenfunctions of the transfer operator (assuming it has one). However, nonautonomous systems require a generalization of this concept, which is predominantly formalized in the framework of transfer operator \emph{cocycles} and multiplicative ergodic theory, where the spectra and eigenfunctions are replaced by Lyapunov spectra (exponents) and associated Lyapunov filtrations (in special cases, Oseledets splittings)~\cite{froyland2010coherent,FrLlQu10,gonzalez2018multiplicative}.
While our considerations below are intimately related to Lyapunov spectra and Oseledets splittings, we defer rigorous elaborations on this connection to future studies. The purpose of the present study is to highlight other aspects.
A description of linear stability in nonautonomous dynamical systems that more broadly mimics the autonomous distinction into stable and unstable directions is given by the Sacker--Sell (or dichotomy) spectrum \cite{SackerSell} which typically contains the Lyapunov spectrum and allows for constructions of dynamically relevant subsets such as invariant manifolds \cite{Potzsche2012, PotzscheRuss}. We will make use of this spectrum within this work.
 
\paragraph{Augmented space.}
The same way as time-augmentation transforms a nonautonomous system into an equivalent autonomous one, we ask whether metastability as a persistence concept for the skew-product associated to a nonautonomous system reveals coherent sets of the latter. This research thread was initiated in~\cite{froyland2017estimating} for periodic driving and in~\cite{froyland2020computation} for finite-time systems, and here we consider the natural next step, where a general ergodic flow is driving a stochastic differential equation (SDE) on physical space. We set up the transfer operator of this skew-product system, a so-called Mather (evolution) semigroup~\cite[Section~6.2]{chicone1999evolution}, and use its spectrum to extract coherent sets. More precisely, we consider the generator of this semigroup, since this has advantages over the numerical approximation of the entire semigroup---as we elaborate below.
A phenomenologically ``dual'' construction in~\cite{GiDa20} considers the evolution semigroup of Koopman operators and extracts so-called ``coherent patterns'', with the crucial difference to our construction that they apply diffusion both in the driving and in physical space. The spectral properties of evolution semigroups (also called Howland semigroup~\cite{howland1974stationary} if the parameter space is the real line representing time) and Mather semigroups in relation to the dynamical properties of the driving system (the ``base'') are comprehensively elaborated in~\cite{chicone1999evolution}.

\subsection{Contributions}

\paragraph{The setting.}

To describe our approach and contributions more precisely, we consider a nonautonomous SDE on the physical space $M\subset \R^d$ (or $M=\tor^d$) with additive noise (and reflecting boundary conditions on $\partial M$, unless $M=\tor^d$) that is driven by some base dynamics on a parameter space $\Theta$,
\begin{equation}
    \label{eq:SDE_intro}
    \begin{aligned}
    d\theta_t &= \Psi(\theta_t) dt, \\
    dx_t &= v(\theta_t, x_t)dt + \varepsilon dw_t,
    \end{aligned}
\end{equation}
i.e., $x_t \in M$ and~$\theta_t \in \Theta$. Here, $\Theta$ is a compact smooth manifold equipped with a continuous vector field~$\Psi$. 
We assume that the unique solution to the ODE on $\Theta$ governed by $\Psi$ is an invertible {ergodic} flow $\phi^t:\Theta \to \Theta$ and that $v:\Theta \times M \to \R^d$ is a smooth, divergence-free vector field.

Ensembles of states in $M$, which will be generalized to be viewed as functions $f\in L^2(M;\C)$, are evolved by~\eqref{eq:SDE_intro} as described by the transfer operator cocycle~$\Per_{\theta}^t$, acting on~$L^2(M;\C)$ equipped with the usual norm~$\|\cdot\|_2$. This means, if $\theta_0 = \theta$ and the random physical-space initial condition is $f$-distributed (we write $x_0\sim f$), then~$x_t \sim \Per_{\theta}^t f$. We note that $\Per_{\theta}^t$ is the solution operator of the Fokker--Planck (or advection-diffusion) equation associated with the SDE for~$x_t$.
The transfer operators satisfy the cocycle identity $\Per^t_{\phi^s \theta} \circ \Per_\theta^s = \Per_\theta^t$. 

In general, nonautonomous transfer operators depend on an initial time $s$ and a time duration~$t$. For a cocycle, the initial time is replaced by the parameter $\theta \in \Theta$ describing the state of the driving system at initial time. Thus, $\Theta$ can be viewed as a compactification of the set~$\R_{\geq 0}$ of initial times~$s$. This compactification allows for the study of asymptotics as $t\to \infty$ without (numerically) emulating this limit. Such asymptotics, like the exponential decay rate of $\norm{\Per_\theta^t f}_2$ for a distribution $f\in L^2(M, \C)$, are commonly studied using multiplicative ergodic theory. In this work, we take a different approach and study the long-term behavior of the transfer operator cocycle by constructing an augmented generator that captures the nonautonomous dynamics in a single time-independent operator. An analogous approach 
has previously been applied to study periodically driven flows in~\cite{froyland2017estimating}. The extension of this approach to aperiodically driven flows introduces significant complexity to the analysis and requires a different theoretical framework which we will elaborate on in the following.

\paragraph{Our results.}

Our central object of study is the linear operator $\mathbf{M}^t$ acting on a function space $\mathfrak F$, which is either $L^2(\Theta , L^2(M, \C))$ or $C(\Theta, L^2(M, \C))$, describing the evolution of some $\f \in\mathfrak F$ fibre-wise according to~$\Per$: 
\[
\big[ \mathbf{M}^t \f \big](\theta) = \Per_{\phi^{-t}\theta}^t \f(\phi^{-t}\theta).
\]
In Section~\ref{sec:cont_time_cocycles}, we introduce a set of assumptions (I)--(IV) for the transfer operators $\Per_{\theta}^t$, guaranteeing that $(\mathbf{M}^t)_{t\geq 0}$ is a strongly continuous semigroup, which is called the \emph{Mather semigroup}. The study of the Mather semigroup and its generator, in particular their spectrum and spectral subspaces, provides insight into the long-term behavior of the transfer operator cocycle~$\Per$. We establish a spectral mapping relation in Theorem~\ref{thm:spectral_mapping}.
In Theorem~\ref{thm:PF_cocycle}, we ensure that the operators $\Per_\theta^t$ of the Fokker--Planck equation related to~\eqref{eq:SDE_intro}, indeed, form a well-defined linear cocycle satisfying assumptions (I)--(IV). 

The purpose of the current work is to study in theory and by several numerical examples how spectral objects of $\mathbf{M}^t$ and its generator $\G$ give rise to families of coherent sets~$A_{\theta}^t \subset M$, $\theta\in\Theta$, $t\ge 0$. To measure coherence of such a family, we will consider two measures. The first is the escape rate, 
which is an asymptotic quantity related to the tail behavior of the survival probabilities $\P(x_s \in A_\theta^s, \: \forall s\in [0, t])$ for~$t\to\infty$. The second is the cumulative survival probability, 
which is an integral quantity of the survival probabilities, and, in particular, includes its finite-time characteristics. 

The spectral objects of $\mathbf{M}^t$ which we extract coherent sets from are eigenfunctions, approximate eigenfunctions and spectral subspaces. We provide rigorous bounds on the measure of coherence of the resulting family $A_\theta^t$ for eigenfuctions (Theorem~\ref{thm:CohFamFromFun}) and approximate eigenfunctions (Proposition~\ref{prop:CohFamfromApprox}) and additionally derive heuristic methods that can be applied in numerical computations.

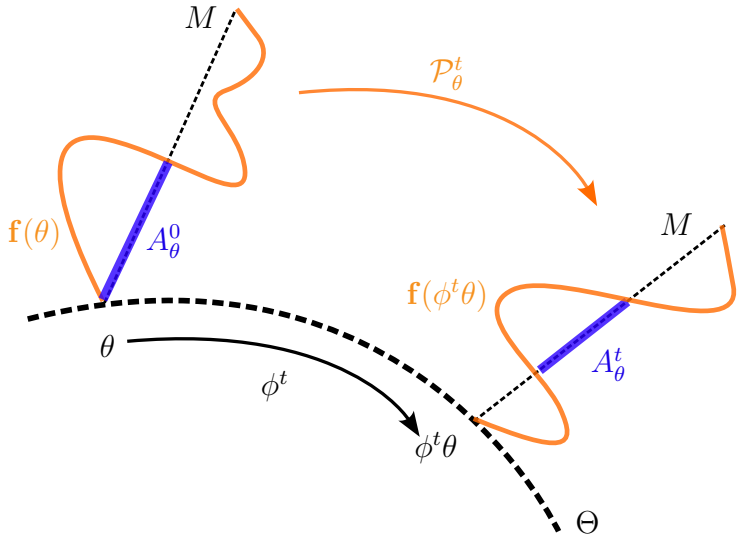
\begin{figure}[htb]
  \centering
\begin{overpic}[width=0.6\linewidth]{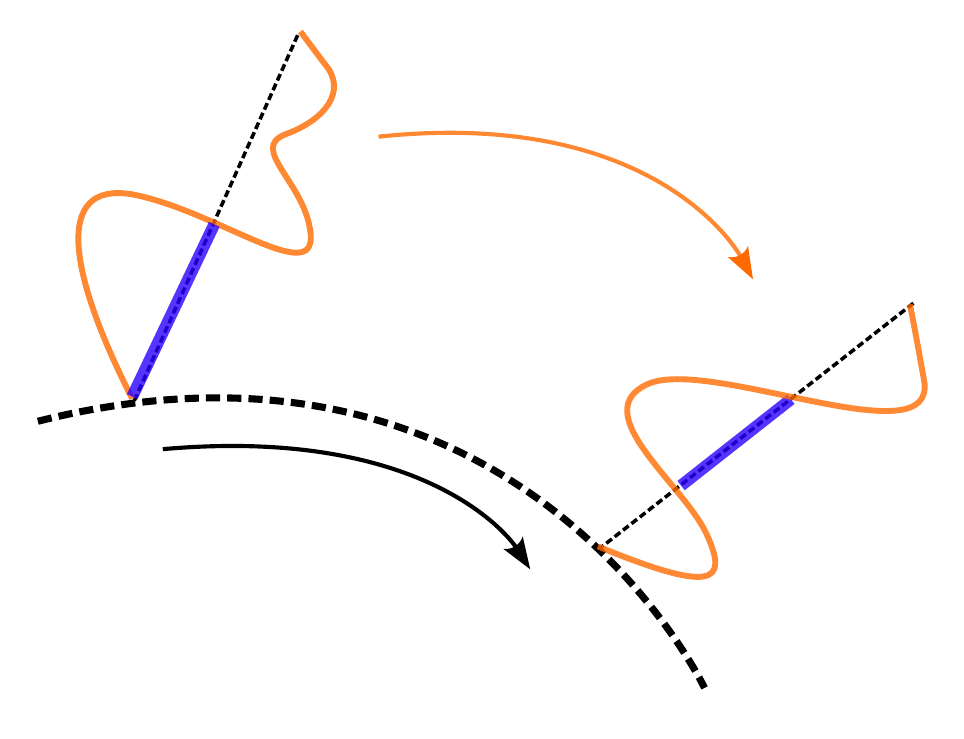}
\colorlet{ColFun}{BurntOrange}
\colorlet{ColSet}{blue!80!violet}
\put(12,28){ $\theta$ }
\put(0,43){ \color{ColFun} $\f(\theta)$ }
\put(52,35){ \color{ColFun} $\f(\phi^t\theta)$ }
\put(53,15){ $\phi^t \theta$ }
\put(33,23){ $\phi^t$ }
\put(55,64){ \color{ColFun} $\Per_{\theta}^t$ }
\put(23,71){ $M$ }
\put(85,44){ $M$ }
\put(18,42){ \color{ColSet} $A_{\theta}^0$ }
\put(76,26){ \color{ColSet} $A_{\theta}^t$ }
\put(74,5){ $\Theta$ }
    \end{overpic}
        \caption{
        Sketch of the temporal evolution of the coherent set $A_\theta^t=\{\f(\phi^t\theta) \geq 0\} \subset M$, where $\f$ is an eigenfunction of the Mather semigroup with real eigenvalue. The distribution $\f(\theta)$ is mapped to (a multiple of) $\f(\phi^t\theta)$ under the cocycle $\Per_\theta^t$ over the driving~$\phi^t: \Theta \to \Theta$.
        }
    \label{fig:ErgCoh_schematic}
\end{figure}

To illustrate the central idea of how we construct coherent sets, assume that $\mathbf f$ is an eigenfunction of the Mather semigroup to a real eigenvalue $\lambda \in \R$, i.e. $\Per_\theta^t \mathbf f(\theta) = e^{\lambda t} \mathbf f(\phi^t\theta)$ for all $t\geq 0$ and $\theta \in \Theta$. Define the family of sets $A_{\theta}^t := \{ \f(\phi^t\theta) \geq 0\} \subset M$ which is the non-negative support of the fibres of $\mathbf f$. Since the transfer operator $\Per_\theta^t$ maps $\mathbf f(\theta)$ to (a multiple of) $\mathbf f(\phi^t \theta)$, particles that were initialized in the set $A_\theta^0$ have a high probability to lie inside $A_\theta^t$ after time~$t$. A schematic representation of this setting is found in Figure~\ref{fig:ErgCoh_schematic}.

To evaluate the applicability of our methods, we numerically compute coherent sets and measure their coherence for three example dynamics: translated gyres, oscillating gyres and alternating shears.
We restrict our considerations to the case where both the parameter space $\Theta$ and the physical space $M$ are given by tori~$\tor^d$. In order to extract coherent sets from respective eigenfuncions, we discretize the Mather semigroup and its generator $\G$ by Galerkin projection onto selected Fourier modes. Based on the infinite-dimensional matrix representation of the augmented generator $\G$ in Proposition~\ref{prop:matrix_Gamma}, we introduce an informed selection process to boost accuracy even in our numerically challenging multidimensional setting, where functions over $\Theta \times M$ need to be discretized.
For simplicity, the driving dynamics on $\Theta$ is assumed to be a quasi-periodic torus rotation. The extension of our methods to general domains and driving dynamics is discussed in the outlook. For all three examples, we compute coherent sets whose survival probabilities decay close to the theoretical bounds. An implementation of our method in Matlab, including all examples presented below, is available under
\begin{center}
    \url{https://github.com/RobinChemnitz/MatherCoherent}
\end{center}

\paragraph{Outline.}
The rest of the paper is structured as follows. Section~\ref{sec:preliminaries} formally introduces the objects which we will work with in the later sections; namely linear cocycles, their Sacker--Sell spectrum, and the Mather semigroup. In Section~\ref{sec:transfer_operator}, we study the SDE \eqref{eq:SDE_intro} through its transfer operators~$\Per_\theta^t$. The central result of this section is that the operators $\Per_\theta^t$ form a linear cocycle  that enables the construction of the associated Mather semigroup and its generator. In Section~\ref{sec:coherent_sets}, we introduce coherent sets and derive methods to extract coherent sets from spectral objects of the Mather semigroup and its generator. In Section \ref{sec:numerics}, we represent the augmented generator $\G$ in Fourier coordinates and construct a discretized generator using a Galerkin approximation onto carefully selected Fourier modes. We extract coherent sets from eigenfunctions of the discretized generator and test their coherence experimentally for three examples 
in Section~\ref{sec:examples}. We conclude with a discussion of our results as well as an outlook.

\section{Preliminaries}\label{sec:preliminaries}
\subsection{Continuous-time cocycles}\label{sec:cont_time_cocycles}
In this section, we introduce the notion of a linear cocycle, its Sacker--Sell spectrum, and the Mather semigroup in a general setting. 

Let $\Theta$ be a compact metric space equipped with its Borel $\sigma$-algebra and let $\H$ be a Hilbert space. We denote the space of bounded linear operators on $\H$ by $\L(\H)$. Let $\phi^t:\Theta \to \Theta$, where $t \in \R$, be a continuous and invertible flow on~$\Theta$. That is, $(\phi^t)_{t\in \R}$ is a one-parameter group on $\Theta$ such that the map $(\theta, t) \mapsto \phi^t \theta$ is continuous. Let $\mu$ be an ergodic measure of $\phi$, that is $\mu(A) \in \{0,1\}$ for any set $A\subset \Theta$ that is $\phi^t$-invariant for all $t\geq 0$. Since $\Theta$ is compact at least one such measure exists, which can be shown using the Banach--Alaoglu theorem. W.l.o.g.~we assume $\textup{supp}(\mu) = \Theta$, i.e.~any open set $U\subset \Theta$ has $\mu(U)>0$. There are two cases that can occur. Either $\Theta \cong S^1$ consists of a single periodic orbit, or the set of aperiodic points with respect to $\phi$ is dense in~$\Theta$. Throughout, we assume that the set of aperiodic points is dense in~$\Theta$.

A linear cocycle in $\H$ with driving $\phi$ is a map
\begin{align*}
\begin{split}
    \Phi: \Theta\times \R_{\geq 0} &\to \L(\H) \\
    (\theta, t) &\mapsto \Phi_\theta^t,
\end{split}
\end{align*}
with the following defining cocycle properties:
\begin{enumerate}[(i)]
    \item $\Phi_\theta^0=\Id, \quad \forall \theta\in \Theta$;
    \item $\Phi_\theta^{s+t} = \Phi_{\phi^s\theta}^t \circ \Phi_\theta^s, \quad \forall \theta \in \Theta, \: \forall s,t \in \R_{\geq 0}$.
\end{enumerate}
We introduce additional assumptions which we impose on (and verify for specific) continuous-time cocycles $\Phi$ throughout this work.
\begin{enumerate}[(I)]
    \item \textit{compact.} The operators $\Phi_\theta^t\in \L(\H)$ are compact $\forall \theta \in \Theta, \: t> 0$;
    \item \textit{norm-continuous in $\theta$.} For each $t\geq 0$, the map $\theta \mapsto \Phi_\theta^t \in \L(\H)$ is norm-continuous;
    \item \textit{strongly-continuous.} For each $x\in \H$, the map $(\theta, t) \mapsto \Phi_\theta^t x \in \H$ is continuous.
    \item \textit{exponentially bounded.} There are constants $K>0$ and $L>0$ such that $\norm{\Phi_\theta^t} \leq K e^{Lt}$, $\forall t\geq 0$.
\end{enumerate}
Assumptions (I) and (II) are necessary to obtain a nice representation of the Sacker--Sell spectrum, cf.\ Proposition~\ref{prop:Sacker_Sell_segments}. Assumptions (III) and (IV) are needed to define the Mather semigroup, cf.\ Section~\ref{sec:Mather semigroup}.

\subsection{The Sacker--Sell spectrum}
This section introduces the Sacker--Sell spectrum, also called dichotomy spectrum, of a continuous-time linear cocycle $\Phi$ satisfying conditions (I)--(IV). For the most part, we follow \cite[Chapter 6.1]{chicone1999evolution}. 
\begin{definition}\label{def:exp_dichotomy}
    We say that a cocycle $\Phi$ has an \emph{exponential dichotomy at} $\lambda\in \R$ if there are constants $\beta>0$ and $C>0$ and a strongly continuous, projection-valued function $\Pi:\Theta\to \L(\H)$ with complemented subspaces
    \begin{equation*}
        S(\theta) := \ran(\Pi(\theta)), \qquad U(\theta) := \ker(\Pi(\theta)),
    \end{equation*}
    that together satisfy the following properties for every $\theta \in \Theta$ and $t\geq 0$:
    \begin{enumerate}[(i)]
        \item $\Pi(\phi^t \theta) \Phi_\theta^t=\Phi_\theta^t \Pi(\theta)$;
        \item The restriction $\Phi_\theta^t |_U : U(\theta) \to U(\phi^t \theta)$ is invertible;
        \item $\norm{\Phi_\theta^t|_S}\leq C e^{(\lambda-\beta) t}, \qquad \big\lVert\left(\Phi_\theta^t|_U \right)^{-1}\big\rVert^{-1} \geq C^{-1} e^{(\lambda+\beta) t},$ \\[7pt]
        where $\Phi_\theta^t|_S$ is defined analogously to $\Phi_\theta^t |_U$.
    \end{enumerate}
    
\end{definition}
The letters $S$ and $U$ stand for the stable and unstable bundles, respectively. 
\begin{remark} 
    One can verify that the first property is equivalent to $\Phi_\theta^t S(\theta) \subset S(\phi^t \theta)$ and $\Phi_\theta^t U(\theta) \subset U(\phi^t \theta)$. Since $\Phi_\theta^t |_U$ is assumed to be invertible and thereby surjective, we can even conclude $\Phi_\theta^t U(\theta) = U(\phi^t \theta)$. The third property implies $\norm{\Phi_\theta^t u} \geq C^{-1} e^{(\lambda + \beta)t}$ for all $u\in U(\theta)$ and $t\geq 0$. If $U(\theta)$ was infinite-dimensional, this would imply that $\Phi_\theta^t$ is not compact. Hence, for any cocycle $\Phi$ satisfying the compactness assumption (I), that admits an exponential dichotomy, the unstable bundle $U(\theta)$ is finite-dimensional.
\end{remark}

\begin{definition}
    The Sacker--Sell spectrum $\Sigma(\Phi)$ is defined as
    \begin{equation*}
        \Sigma(\Phi) := \{\lambda \in \R \:|\: \Phi \textup{ does not have an exponential dichotomy at }\lambda \}.
    \end{equation*}
\end{definition}

By definition, $\Sigma(\Phi)\subset \R$ is a closed set. For cocycles $\Phi$ satisfying assumptions (I)--(IV), the Sacker--Sell spectrum is of the following form. For a proof, see \cite[Theorem 8.12]{chicone1999evolution}.
\begin{proposition}\label{prop:Sacker_Sell_segments}
    The Sacker--Sell spectrum $\Sigma(\Phi)$ consists of a finite or countably infinite number of closed segments
    \begin{equation*}
        \Sigma(\Phi) = \bigcup_{k=1}^N [r_k^-, r_k^+],
    \end{equation*}
    where $N \in \N_0 \cup \{\infty\}$ and $-\infty \leq \hdots <r_2^- \leq r_2^+ < r_1^- \leq r_1^+$.
    
    Assume $\dim(\H)=\infty$. If $N<\infty$, then $r^-_N=-\infty$. If $N=\infty$, we find  
    \begin{equation*}
        \lim_{k \to \infty} r_k^\pm = -\infty.
    \end{equation*}
\end{proposition}
Note that $r_k^- = r_k^+$ is allowed, such that the segment $[r_k^-, r_k^+]$ consists of a single point. In general, the cases $\Sigma(\Phi) = \emptyset$ and $\Sigma(\Phi) = (-\infty, r_1^+]$ can occur. 
In Section \ref{sec:autonomous_case} we further specify the structure of the Sacker--Sell spectrum in the special case where $\Phi$ is autonomous, i.e.~the operators $\Phi_\theta^t$ do not depend on~$\theta$.

\begin{remark}
    In the literature, e.g.~\cite{chicone1999evolution}, exponential dichotomies are commonly only defined at $\lambda = 0$. A cocycle $\Phi$ has an exponential dichotomy at $\lambda$ if and only if the cocycle defined by $\Psi_\theta^t := e^{-\lambda t} \Phi_\theta^t$ has an exponential dichotomy at $0$. This type of rescaling argument is used on multiple occasions throughout this work.
\end{remark}

\subsection{The Mather semigroup}
\label{sec:Mather semigroup}

In this section we consider the evolution of elements under a cocycle $\Phi$ for all $\theta \in \Theta$ simultaneously. Let $\mathfrak{F}$ be either the Bochner space $L^2(\Theta, \H)$ or $C(\Theta, \H)$, where the latter is the space of continuous functions equipped with the supremum norm~$\norm{\cdot}_\infty$. Functions in $\mathfrak{F}$ are denoted by bold letters. For each $t \geq 0$, define the linear operator~$\mathbf{M}^t : \mathfrak{F} \to\mathfrak{F}$ by
\begin{equation}\label{eq:def_Mather}
    \big[ \mathbf{M}^t \f \big](\theta) = \Phi_{\phi^{-t}\theta}^t \f(\phi^{-t}\theta).
\end{equation}
The operators $(\mathbf{M}^t)_{t\ge 0}$ satisfy the semigroup property, i.e.~$\mathbf{M}^t \circ \mathbf{M}^s = \mathbf{M}^{s+t}$, for all $s, t \geq 0$. The semigroup $(\mathbf{M}^t)_{t\geq 0}$ is called the \emph{Mather semigroup}. The action of the Mather semigroup can be interpreted as evolving $\f$ fibre-wise through the cocycle $\Phi$. Hence, the $\theta$-fibre of $\f$ is mapped to the $\phi^t \theta$-fibre of $\mathbf{M}^t \f$. Different fibres of $\f$ do not interact.

In this section, we state our results for the Mather semigroup constructed over either $L^2(\Theta, \H)$ or $C(\Theta, \H)$. In Section \ref{sec:numerics}, we discretize the operators $\mathbf{M}^t$ using a Galerkin projection onto Fourier modes, which requires us to consider the Mather semigroup over the Hilbert space $L^2(\Theta, \H)$. In Section \ref{sec:coherent_sets}, we need higher regularity than $L^2$, for which considering $C(\Theta, \H)$ is necessary. We note that the chosen function space only plays a role in the theoretical study, since a computational implementation of our methods is unable to distinguish continuous and measurable functions. 
For the study of the Mather semigroup over different function spaces in a more general setting, see~\cite{barreira2019hyperbolicity}. 

\begin{figure}
    \centering
    \begin{tikzpicture}
        \definecolor{myblue}{rgb}{0, 0.4470, 0.7410} % the default Matlab color "#0072BD" used in the eigenvalue plots
        \tikzstyle{bands}=[myblue!70];
        \tikzstyle{lines}=[myblue!70, line width=1.5pt];

        \draw[thick, {Latex}-{Latex}] (0, -3.2) -- (0, 3.2) node[anchor=west]{Im};
        \draw[thick, {Latex}-{Latex}] (-5, 0) -- (1.2, 0)  node[anchor=north]{Re};
        \fill[bands] (-4.5, -3) rectangle (-4, 3);
        \fill[bands] (-1.7, -3) rectangle (-0.5, 3);
        \fill[bands] (0.4, -3) rectangle (0.7, 3);
         
        \draw[lines] (-4.5, 3) -- (-4.5, -3);
        \draw[lines] (-4, 3) -- (-4, -3);
        \draw[lines] (-3, -3) -- (-3, 3); 
        \draw[lines] (-1.7, 3) -- (-1.7, -3);
        \draw[lines] (-0.5, 3) -- (-0.5, -3);
        \draw[lines] (0.4, 3) -- (0.4, -3);
        \draw[lines] (0.7, 3) -- (0.7, -3);

        \draw[thick, {Latex}-{Latex}] (6, -3.2) -- (6, 3.2) node[anchor=west]{Im};
        \draw[thick, {Latex}-{Latex}] (2.8, 0) -- (9.2, 0)  node[anchor=north]{Re}; 
        \fill [bands,even odd rule] (6,0) circle[radius=2.7607cm] circle[radius=2.5999cm] circle[radius=2.1716cm] circle[radius=1.7082cm]  circle[radius=1.0784cm] circle[radius=0.9758cm];
        
        \draw[lines] (6,0) circle[radius=2.7607cm] circle[radius=2.5999cm] circle[radius=2.1716cm]  circle[radius=1.7082cm] circle[radius=1.0784cm] circle[radius=0.9758cm];
        \draw[lines] (6,0) circle[radius=1.3171cm];
        \draw[dashed] (6,0) circle[radius=2.4cm];
    \end{tikzpicture}
    \caption{A schematic representation of the spectrum of the generator $\G$ (left), and the spectrum of the Mather operator $\mathbf{M}^1$ (right). The unit circle $S^1$ is marked as a dashed circle. In this example, the Sacker--Sell spectrum consists of four intervals, one of which is a single point.} 
    \label{fig:schematic_spectrum}
\end{figure}
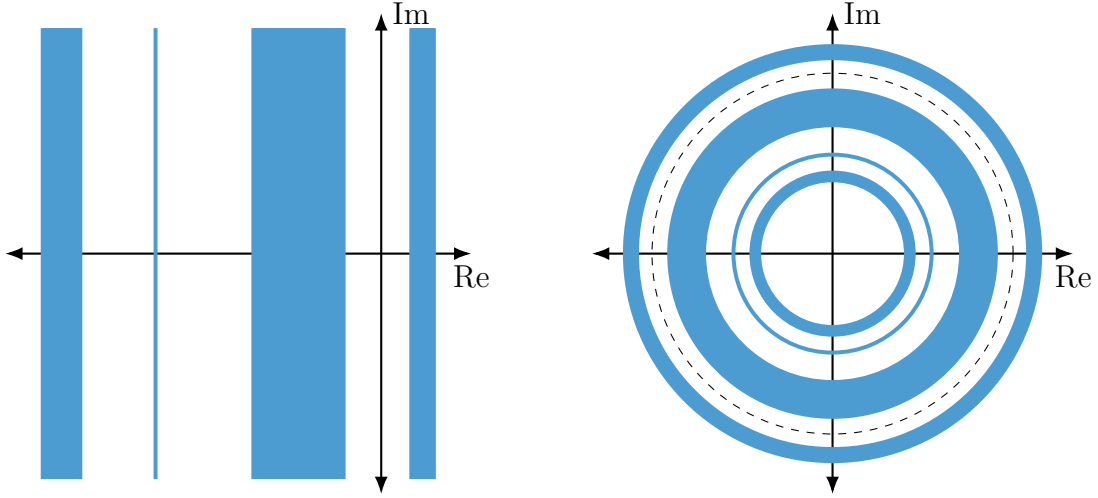

By \cite[Theorem~6.20 and Lemma~6.33]{chicone1999evolution}, $(\mathbf{M}^t)_{t\in \R_{\geq 0}}$ is a strongly continuous semigroup, given that $\Phi$ is strongly continuous and exponentially bounded, which is satisfied with assumptioins (III) and~(IV). Hence, the Mather semigroup has a generator defined pointwise by
\begin{equation*}
    \G = \lim_{t\to 0} \frac{\mathbf{M}^t - \Id}{t}.
\end{equation*}
The (potentially unbounded) operator $\G$ is closed and has a dense domain $\mathcal D(\G)\subset \mathfrak{F}$, see e.g.~\cite[Corollary 2.5]{pazy2012semigroups}. We call $\G$ the \emph{augmented generator}. The following theorem characterizes the spectrum of $\mathbf{M}^t$ and its generator $\G$ in terms of the Sacker--Sell spectrum $\Sigma(\Phi)$. A schematic representation of the spectra $\sigma(\G)$ and $\sigma(\mathbf{M}^t)$ is given in Figure \ref{fig:schematic_spectrum}. We note that even when the spectrum of $\mathbf{M}^t$ and $\G$ is non-empty, their point spectrum may be empty, i.e.~they might not possess eigenfunctions.

\begin{theorem}\label{thm:spectral_mapping}
    The spectra of $\mathbf{M}^t$, where $t>0$, and $\G$ are given by
    \begin{nalign}
    \sigma(\mathbf{M}^t) \setminus \{0\} &= \{e^{\lambda t + \eta i} \hspace{1.24mm}\: | \: \lambda \in \Sigma(\Phi), \:\eta\in [0, 2\pi) \}, \\
        \sigma(\G) &= \{\lambda + \eta i \: | \: \lambda \in \Sigma(\Phi) , \: \eta \in \R \}.
    \end{nalign}
    In other words, $\sigma(\mathbf{M}^t)$ consists of annuli around the center, and $\sigma(\G)$ consists of bands parallel to the imaginary axis.
\end{theorem}
\begin{proof}
    The spectral mapping theorem \cite[Theorem 6.30, Theorem 6.37]{chicone1999evolution} states that for all $t>0$ we find
    \begin{equation}\label{eq:spectral_mapping}
        \sigma(\mathbf{M}^t) \setminus \{0\} = e^{t \sigma(\G)}.
    \end{equation}
    Additionally, it states that $\sigma(\mathbf{M}^t)$ is invariant under rotations around the center and that $\sigma(\G)$ is invariant under translations by imaginary values. The cocycle $\Phi$ has an exponential dichotomy at $0$ if and only if $\sigma(\mathbf{M}^t) \cap S^1 = \emptyset$. For $\mathfrak F= L^2(\Theta, \H)$, this follows from \cite[Assertion 1.4]{latushkin1991weighted}. For $\mathfrak F = C(\Theta, \H)$, the statement is shown in \cite[Theorem 6.41]{chicone1999evolution}. By rescaling, the cocycle $\Phi$ has an exponential dichotomy at $\lambda \in \R$ if and only if $\sigma(\mathbf{M}^t) \cap e^{\lambda t}S^1 = \emptyset$. Since the Sacker--Sell spectrum $\Sigma(\Phi)$ consists exactly of those $\lambda \in \R$ for which $\Phi$ does not have an exponential dichotomy, we conclude
    \begin{equation*}
        \sigma(\mathbf{M}^t) \setminus \{0\} = \{e^{\lambda t + \eta i} \hspace{1.24mm}\: | \: \lambda \in \Sigma(\Phi), \:\eta \in [0, 2\pi) \}.
    \end{equation*}
    The characteriziation of the spectrum of $\G$ follows directly from \eqref{eq:spectral_mapping} together with the fact that $\sigma(\G)$ is invariant under translations by imaginary values.
\end{proof}

\subsection{Spectral projections}\label{sec:spectral_projections}
Theorem \ref{thm:spectral_mapping} provides a  direct link between the spectrum of the Mather semigroup $\mathbf{M}^t$ and the Sacker--Sell spectrum $\Sigma(\Phi)$. In this section, we show  that this connection extends to the level of the corresponding subspaces, in the sense that the spectral projections of $\mathbf{M}^t$ can be described in terms of the projection-valued functions that make up the exponential dichotomies of $\Phi$, cf.~\eqref{eq:Riesz_range}. In the following, we give a brief introduction to spectral projections.

In a finite-dimensional vector space $V$, each eigenvalue $\lambda$ of a linear operator $T:V\to V$ has a corresponding generalized eigenspace $E_\lambda$. These generalized eigenspaces decompose the vector space $V$, and the spectrum of the restriction $T|_{E_\lambda}$ is the set $\{\lambda\}$. For an infinite-dimensional Banach space $\mathds{B}$ and a closed linear operator $T:\mathds{B} \to \mathds{B}$, a similar decomposition is obtained using Riesz-projections \cite{riesz2012functional}, also called \emph{spectral projections}. For a more recent introduction, we refer to \cite[Chapter 6]{van2022functional}. Let $\Lambda\subset \sigma(T)$ be an isolated part of the spectrum of $T$. We call $\Lambda$ a \emph{spectral set}.  Consider a closed contour $\Gamma \subset \rho(T)$ around $\Lambda$ with positive orientation. The Riesz-projection $\mathfrak{P}_\Lambda:\mathds{B} \to \mathds{B}$ is defined by
\begin{equation}
    \mathfrak{P}_\Lambda = -\frac{1}{2\pi i} \oint_{\Gamma} (T- z )^{-1} dz.
\end{equation}
The linear operator $\mathfrak{P}_\Lambda:\mathds{B} \to \mathds{B}$ is a bounded projection that splits the Banach space into its range $\ran(\mathfrak{P}_\Lambda) \subset \mathds{B}$ and its kernel $\ker(\mathfrak{P}_\Lambda) \subset \mathds{B}$. The subspaces $\ran(\mathfrak{P}_\Lambda)$ and $\ker(\mathfrak{P}_\Lambda)$ are invariant under $T$ and the spectrum of their restrictions is given by $\sigma(T|_{\ran(\mathfrak{P}_\Lambda)}) = \Lambda$ and $\sigma(T|_{\ker(\mathfrak{P}_\Lambda)}) = \sigma(T) \setminus \Lambda$. Hence, the subspace $\ran(\mathfrak{P}_\Lambda)$ can be considered as an infinite-dimensional version of the generalized eigenspace of the spectral set $\Lambda \subset \sigma(T)$. We call $\ran(\mathfrak{P}_\Lambda)\subset \mathds{B}$ the \emph{spectral subspace} of $\Lambda$.

We use Riesz-projections to describe a connection between the spectrum of the Mather semigroup $\mathbf{M}^t$ and the stable/unstable bundles of the cocycle $\Phi$. Let $[r^-, r^+]$ be an isolated segment of the Sacker--Sell spectrum $\Sigma(\Phi)$, and let $\gamma^- < r^- \leq r^+ < \gamma^+$ be such that $[\gamma^-, \gamma^+] \cap \Sigma(\Phi) = [r^-, r^+]$. By Theorem \ref{thm:spectral_mapping}, the annulus $\Lambda^t:=\{e^{\lambda t + \eta i} \:|\: \lambda \in [r^-, r^+], \: \eta\in [0, 2\pi)\}$ is an isolated component of the spectrum of $\mathbf{M}^t$, and the rings of radius $e^{\gamma^- t}$ and $e^{\gamma^+ t}$ separate $\Lambda^t$ from the rest of the spectrum of $\mathbf{M}^t$. Let $\mathds{D}\subset\C$ be the unit disk. The Riesz-projection onto the spectral sets $\sigma(\mathbf{M}^t) \cap e^{\gamma^- t} \mathds{D}$, respectively $\sigma(\mathbf{M}^t) \cap e^{\gamma^+ t} \mathds{D}$, are given by
\begin{equation*}
    \mathfrak{P}_{e^{\gamma^- t}} := -\frac{1}{2\pi i} \oint_{e^{\gamma^- t} S^1} (\mathbf{M}^t - z )^{-1} dz , \quad \mathfrak{P}_{e^{\gamma^+ t}} := -\frac{1}{2\pi i} \oint_{e^{\gamma^+ t} S^1} (\mathbf{M}^t - z )^{-1} dz.
\end{equation*}
In particular, the projection onto the spectral set $\Lambda^t$ is given by
\begin{equation*}
    \mathfrak{P}_{\Lambda^t} = \mathfrak{P}_{e^{\gamma^+ t}} - \mathfrak{P}_{e^{\gamma^- t}}.
\end{equation*}
Since $\gamma^-$ and $\gamma^+$ lie outside of $\Sigma(\Phi)$, the cocycle $\Phi$ has exponential dichotomies at $\gamma^-$ and~$\gamma^+$. Therefore, there are strongly continuous projections $\Pi_-(\theta)$, $\Pi_+(\theta)$ that split the space $\H$ into the respective stable and unstable bundles $S_\pm(\theta) := \ran(\Pi_\pm (\theta))$ and $U_\pm(\theta) := \ker(\Pi_\pm (\theta))$. Note that $S_-(\theta) \subset S_+(\theta)$ and $U_+(\theta) \subset U_-(\theta)$. In particular, for each $\theta \in \Theta$, we obtain a splitting $\H = S(\theta) \oplus E(\theta) \oplus U(\theta)$, where $S(\theta) := S_-(\theta)$, $E(\theta):= S_+(\theta) \cap U_-(\theta)$, and~$U(\theta):= U_+(\theta)$. The operator $\Pi(\theta):=\Pi_+(\theta)-\Pi_-(\theta)$ is a projection onto $E(\theta)$ along $S(\theta)\oplus U(\theta)$. Vectors $s\in S(\theta)$ decay uniformly faster than the rate $\gamma^-$ under $\Phi_\theta^t$, vectors $u \in U(\theta)$ grow uniformly faster than rate the~$\gamma^+$, and vectors $x\in E(\theta)$ decay/grow uniformly at rate between $\gamma^-$ and~$\gamma^+$.

The connection between the spectral projection $\mathfrak{P}_{\Lambda^t}$ and the dichotomy projections $\Pi_\pm (\theta)$ is of the following form: The projections $\mathfrak{P}_{e^{\gamma^\pm t}}$ are given by
\begin{equation*}
    [\mathfrak{P}_{e^{\gamma^\pm t}} \f] (\theta) = \Pi_\pm(\theta) \f(\theta).
\end{equation*}
For $\mathfrak F = L^2(\Theta, \H)$, this is shown in \cite[Assertion 1.4]{latushkin1991weighted}. For $\mathfrak F = C(\Theta, \H)$, the statement is found in \cite[Theorem 6.38]{chicone1999evolution}.
We conclude that the projection $\mathfrak{P}_{\Lambda^t}$ is given by
\begin{equation*}
    [\mathfrak{P}_{\Lambda^t} \f](\theta) = \Pi(\theta)\f(\theta).
\end{equation*}
Therefore, the spectral subspace of $\Lambda^t$ is given by
\begin{nalign}\label{eq:Riesz_range}
    \ran(\mathfrak{P}_{\Lambda^t}) &= \{\f \in \mathfrak F \: | \: \f(\theta) \in E(\theta), \: \forall \theta \in \Theta\}.
\end{nalign}
In the case $\mathfrak F = L^2(\Theta, \H)$, functions $\f\in \mathfrak F$ are only defined up to $\mu$-null sets. In that case, \eqref{eq:Riesz_range} is to be understood as an $\mu$-a.e.~statement.
We conclude that the spectral subspace of an annulus $\Lambda^t$ is made up of functions $\f \in \mathfrak F$ whose fibres decay/grow uniformly at rate between $r_-$ and $r_+$.

\subsection{The autonomous case}\label{sec:autonomous_case}
In this section we study the special case where the operators $\Phi_\theta^t$ do not depend on $\theta \in \Theta$, but only on~$t\geq 0$. We will encounter such cocycles in examples of Section~\ref{sec:examples}. Since $\Phi_\theta^t$ does not depend on $\theta$, we write $\Phi^t\in\L(\H)$ for the time-$t$-operator. By the cocycle property, $(\Phi^t)_{t\geq 0}$ is a one-parameter semigroup, i.e.~$\Phi^0=\Id$, and for all $s,t\geq 0$ we find $\Phi^{s+t}=\Phi^t \circ \Phi^s$. Since $\Phi$ is a cocycle satifying assumptions (I)--(IV), $(\Phi^t)_{t\geq 0}$ is a strongly continuous, exponentially bounded semigroup that is compact for $t>0$. For the general study of one-parameter semigroups, we refer to standard textbooks~\cite{engel2000one,pazy2012semigroups}.

By \cite[Proposition V.1.15]{engel2000one}, the cocycle $\Phi$ has an exponential dichotomy at $0$ if and only if $\sigma(\Phi^t) \cap S^1 = \emptyset$ for one, and therefore all, $t>0$. Consequently, $\Phi$ has an exponential dichotomy at $\lambda \in \R$ if and only if $\sigma(\Phi^t) \cap e^{\lambda t }S^1 = \emptyset$ for one, and therefore all,~$t>0$. Given that $\Phi$ has an exponential dichotomy at $\lambda \in \R$, the corresponding projection-valued function $\Pi:\Theta \to \L(\H)$ is constant. For $t>0$, the operator $\Phi^t$ is compact. By the Riesz--Schauder theorem, cf.~\cite[Theorem 7.11]{van2022functional}, the spectrum $\sigma(\Phi^t)$ consists of discrete points $(z_k)_{1\leq k \leq N}$ for $N\in \N \cup \{\infty\}$ that may only accumulate at $0$. Additionally, each nonzero element of the spectrum is an eigenvalue, i.e.~lies in the point spectrum. We conclude that the Sacker--Sell spectrum of $\Phi$ is given by
\begin{equation}\label{eq:autonomous Sacker--Sell}
    \Sigma(\Phi) = \big\{ \log(\abs{z_k}) \mid 1\leq k \leq N \big\}. 
\end{equation}
Hence, $\Sigma(\Phi)$ consists only of discrete points. Let $\hdots<\lambda_2 < \lambda_1$ be the discrete values such that $\Sigma(\Phi) = \{\lambda_k \mid 1\leq k \leq N\}$. For each $\lambda_k$ there is a vector $x\in \H$ and a complex value $z_k\in \C$ with $\abs{z_k}=\lambda_k$ such that for each $t\geq 0$ we find~$\Phi^t x = e^{z_kt} x$. 

Theorem \ref{thm:spectral_mapping} implies that the spectrum of $\mathbf{M}^t$ consists of circles around the origin at radii $e^{\lambda_k t}$ and the spectrum of $\G$ consists of lines parallel to the imaginary axis at real parts $\lambda_k$. Furthermore, the existence of eigenfunctions is guaranteed. Let $x\in \H$ such that $\Phi^t x = e^{z_kt} x$ for all $t \geq 0$. Then, the function $\f\in \mathfrak F$ that admits the constant value $x$, i.e.~$\f(\theta)=x$ for all $\theta\in \Theta$, is an eigenfunction of the Mather operators $\mathbf{M}^t$ and of the generator $\G$ with eigenvalue $e^{z_k t}$ and~$z_k$, respectively. 

\section{The transfer-operator cocycle}\label{sec:transfer_operator}
We consider the evolution of particles on a compact domain advected by a nonautonomous vector field and subject to small diffusion. The goal of our analysis is to characterize and numerically compute \emph{coherent sets}, which are time-dependent regions in space that only mix slowly with their surrounding. This problem has been studied for finite-time horizons in \cite{denner2016computing,froyland2020computation} and for periodically driven vector fields in~\cite{froyland2017estimating}. In this work, we apply a similar approach to find coherent sets in vector fields that are driven by an ergodic base dynamics.

Let $\Theta$ be a compact smooth manifold equipped with a continuous vector field $\Psi$. Let $M\subset \R^d$ be an open, bounded domain with piecewise $C^4$ boundary (such that $\overline{M}\subset \R^d$ is compact) or 
a torus $M = \tor^d\cong[0,1]^d_{/\sim}$. We call $\Theta$ the \emph{parameter space} and $M$ the \emph{physical space}. We consider the nonautonomous SDE 
\begin{nalign}\label{eq:SDE}
    d\theta_t &= \Psi(\theta_t) dt,\\
    dx_t &= v(\theta_t, x_t)dt + \varepsilon\, dw_t,
\end{nalign}
where~$\varepsilon>0$ and $w_t$ is a $d$-dimensional standard Wiener process.
We impose reflecting boundary conditions, unless $M=\tor^d$.  We assume that $\Psi$ is a continuous vector field on $\Theta$ that generates an invertible flow $\phi^t:\Theta \to \Theta$, for which the aperiodic points are dense, and fix an ergodic measure $\mu$ on $\Theta$. The vector field $v:\Theta \times \overline{M} \to \R^d$ is assumed to be divergence-free and smooth in both variables. For all $\theta \in \Theta$, the vector field $v(\theta, \cdot)$ should have no flow across the boundary of $M$, i.e.~$v(\theta, x) \cdot n(x) = 0$ for all $x\in \partial M$, where $n(x)$ is the outer normal unit vector in $x$. 

We study how distributions of particles evolve in time under the SDE (\ref{eq:SDE}) on the domain $M$ with reflecting boundary conditions. For a formal introduction to the Fokker--Planck equation and its unique solvability, we refer to \cite[Chapter 2]{stahn2022augmented}. In the following, we consider particle distributions represented by functions $f\in L^2(M, \R)$. The temporal evolution of an initial distribution $f_0\in L^2(M, \R)$ under the SDE (\ref{eq:SDE}) with starting parameter $\theta_0 \in \Theta$ is governed by the Fokker--Planck equation
\begin{equation}\label{eq:FP}
    \partial_t f(t,x) = \frac{1}{2} \varepsilon^2 \Delta_x f(t, x) - \text{div}_x\big(f(t, \cdot) v(\phi^t \theta_0, \cdot)\big)(x), \quad f(0, \cdot) = f_0, \quad \frac{\partial f(t, \cdot)}{\partial n} \Bigr|_{\partial M} = 0.
\end{equation}
The third term is the reflecting boundary condition, which requires that the derivative of $f(t, \cdot)$ in the direction of the outer normal unit vector on $\partial M$ vanishes for all $t>0$. We use the convention $f(t) := f(t, \cdot) \in L^2(M, \R)$. The Fokker--Planck equation can be written as a nonautonomous abstract Cauchy problem (NACP) on the space $C(\R_{\geq 0}, L^2(M, \R))$ 
\begin{nalign}\label{eq:NACP}
    \partial_t f(t) &= G(\phi^t\theta_0) f(t),\\
    f(0) &= f_0.
\end{nalign}
The generator $G(\theta)$ is the unbounded operator defined on a dense domain $\mathcal{D}(G(\theta))\subset L^2(M, \R)$, and is given by
\begin{nalign}\label{eq:FP_generator}
   \big[G(\theta) f\big] (x) &= \frac{1}{2} \varepsilon^2\Delta_x f(x) - \textup{div}_x\left(v(\theta, \cdot)f\right)(x) \\
    &= \frac{1}{2} \varepsilon^2 \Delta_x f(x) - \nabla_x f(x) \cdot v(\theta, x).
\end{nalign}
The simplification of the divergence term is possible since $v$ was assumed to be divergence-free. The implementation of reflecting boundary conditions is discussed below.
\begin{definition} 
A classical solution to the NACP \eqref{eq:NACP} is a function $f\in C(\R_{\geq 0}, L^2(M, \R)) \cap C^1(\R_{>0}, L^2(M, \R))$ with $f(0) = f_0$ that, for all $t>0$ satisfies $\partial_t f(t) = G(\phi^t\theta_0) f(t)$ and $f(t) \in \mathcal{D}(G(\phi^t\theta_0))$.
\end{definition}

The reflecting boundary condition can be encoded in the domain of the operator. We choose, independently of $\theta$, 
\begin{equation}\label{eq:generator_domain}
    \mathcal{D}(G) = \mathcal{D}(G(\theta)) := \Big\{f\in H^2(M) \,\Big\vert\,\, \frac{\partial f}{\partial n} \Bigr|_{\partial M} = 0\Big\},
\end{equation}
where $H^2(M) = W^{2,2}(M)$ the usual Sobolev space. In the case $M=\tor^d$ there is no boundary and the domain of the generator is simply $H^2(M)$. Note that $\mathcal{D}(G)\subset L^2(M, \R)$ is dense. With this choice of the domain of $G$, any classical solution to the NACP \eqref{eq:NACP} naturally satisfies the boundary condition, and hence, is a solution to the Fokker--Planck equation \eqref{eq:FP} and vice versa. Therefore, we call both \eqref{eq:FP} and \eqref{eq:NACP} the Fokker--Planck equation and use whichever formulation is suitable in the situation at hand. The following theorem asserts the existence of a unique classical solution to the Fokker--Planck equation.
\begin{theorem}\label{thm:classical_solution}
    For any $f_0 \in L^2(M, \R)$ and $\theta_0 \in \Theta$, the NACP \eqref{eq:NACP} with domain $\mathcal{D}(G)=\{f\in H^2(M) \:| \: \frac{\partial f}{\partial n} |_{\partial M} = 0\}$ has a unique classical solution.
\end{theorem}
\begin{proof}
     See \cite[Theorem 3.12]{stahn2022augmented}. 
     For further details on parabolic PDEs generating NACPs we refer to the classical textbooks~\cite{lunardi2012analytic,tanabe2017functional}.
\end{proof}

For a starting parameter $\theta \in \Theta$ and a time $t>0$, we define $\Per_{\theta}^{t}: L^2(M, \R) \to L^2(M, \R)$ to be the solution operator, also called Perron--Frobenius operator or Kolmogorov forward operator, for time $t$ of the Fokker--Planck equation, i.e.~$\Per_\theta^t f_0 = f(t)$.  By convention, we set $\Per^0_{\theta} = \Id$. It is well-known, that the operators $\Per_\theta^t$ are Markov operators, i.e.~
\begin{equation*}
    f_0\geq 0 \quad \Rightarrow \quad \Per_\theta^t f_0 \geq 0, \quad \text{and} \quad \int_M [\Per_\theta^t f_0](x) \: dx = \int_M f_0(x) \:dx.
\end{equation*}
Since $v$ is divergence-free, the constant function $\1\in L^2(M, \R)$ is an invariant solution to the Fokker--Planck equation, i.e.~$\Per_\theta^t \1 = \1$ for all $\theta\in \Theta$ and $t\geq 0$. Due to diffusion, $\Per_\theta^t f$ converges to a constant function as $t\to \infty$ for any $f\in L^2(M, \R)$. Hence, functions that integrate to 0 converge to the 0-function under~$\Per$. Proposition~\ref{prop:spectral_gap} below gives a quantitative formulation of this fact.

We can extend $\Per_\theta^t$ by linearity to act on complex-valued functions. Define the Hilbert space
\begin{equation}
    \H := L^2(M, \C).
\end{equation}
By construction, the operators $\Per_\theta^t$ are real, i.e.~they map real-valued functions to real-valued functions. We also define the space of functions whose integral vanishes
\begin{equation}
    \H_0:= \{f \in L^2(M, \C) \mid \textstyle{\int_M} f(x) dx = 0\}.
\end{equation}
Since the operators $\Per_\theta^t$ leave the integral of a function invariant, the subspace $\H_0$ is forward invariant, i.e.~$\Per_\theta^t \H_0 \subset \H_0$. From now on, we restrict our attention to the subspace $\H_0$ and study the solution operators $\Per_\theta^t$ as operators on $\H_0$.
\begin{proposition}\label{prop:spectral_gap}
    There is a constant $\varrho<0$ such that the transfer operators $\Per_\theta^t : \H_0 \to \H_0$ satisfy 
    \begin{equation*}
        \norm{\Per_\theta^t} \leq e^{\varrho t},
    \end{equation*}
    for all $\theta \in \Theta$ and $t\geq 0$.
\end{proposition}
\begin{proof}
    See Appendix \ref{ap:proofs}.
\end{proof}

We address the regularity of the operators $\Per_\theta^t$ in both $\theta$ and $t$. The Perron--Frobenius transfer operator of a deterministic dynamic is strongly continuous with respect to the vector field and to time, but not norm continuous. In the case of additive noise, \cite{koltai2019frechet} showed that $\Per_\theta^t$ is not only norm continuous in $\theta$, but even Fr\'echet-differentiable. The following lemma characterizes the regularity of $\Per_\theta^t$ in both $\theta$ and $t$ simultaneously.
\begin{lemma}\label{lem:uniform_regularity}
    Fix $t\geq 0$ and $\theta_1,$ $\theta_2\in \Theta$. Let $\delta > 0$ such that
    \begin{equation*}
        \norm{v(\phi^s \theta_2, x) - v(\phi^s \theta_1, x)}_e \leq \delta,
    \end{equation*}
    in the Euclidean norm, for all $s\in [0,t]$ and $x \in M$. Then, the following bound in operator norm holds:
    \begin{equation*}
        \norm{\Per_{\theta_1}^t - \Per_{\theta_2}^t} \leq \frac{\delta}{\varepsilon} \sqrt{t}.
    \end{equation*}
\end{lemma}
\begin{proof}
    See Appendix \ref{ap:proofs}.
\end{proof}

We verify that the solution operators $\Per$ form a cocycle over $\phi$ satisfying assumptions (I)--(IV).
By invariance of the subspace $\C \mathds{1}$, analogous assertions hold when $\Per$ is considered as a cocycle over $\H$ instead of~$\H_0$. 

\begin{theorem}\label{thm:PF_cocycle}
    The Perron--Frobenius operators $\Per_\theta^t:\H_0 \to \H_0$ of the Fokker--Planck equation~\eqref{eq:FP} form a well-defined linear cocycle over the driving $\phi$ that satisfies assumptions (I)--(IV) introduced in Section~\ref{sec:cont_time_cocycles}. 
\end{theorem}
\begin{proof}
    To apply the results of \cite{stahn2022augmented}, we restrict ourselves to the time interval $[0,T]$ for now and fix a starting parameter $\theta \in \Theta$. By \cite[Theorem 3.12]{stahn2022augmented}, the solution operator to the Fokker--Planck equation \eqref{eq:FP} is given by a two-parameter evolution family $P_{s,t} : \H \to \H$, where $P_{s,t}$ is the solution operator from time $s\geq 0$ to time $t\geq s$. By definition, we find $\Per_\theta^t = P_{0, t}|_{\H_0}$, as well as $\Per_{\phi^s \theta}^t = P_{s, s+t}|_{\H_0}$. The evolution property of $P$ yields
    \begin{equation*}
        \Per_\theta^{s+t} = P_{0, s+t}|_{\H_0} = P_{s, s+t}|_{\H_0} \circ P_{0, s}|_{\H_0} = \Per_{\phi^s \theta}^t \circ \Per_\theta^s.
    \end{equation*}
    This holds for all $\theta \in \Theta$ and $0\leq s \leq t \leq T$. Since $T>0$ was arbitrary, the operators $\Per$ form a well-defined cocycle over $\phi$. 
    
    We show that the cocycle $\Per$ satisfies assumptions (I)--(IV) defined in Section \ref{sec:cont_time_cocycles}.
    \begin{enumerate}[(I)]
        \item \textit{compact.} Compactness of $\Per_\theta^t$ for $t>0$ follows from \cite[Theorem 3.12]{stahn2022augmented}.
        \item \textit{norm-continuous in $\theta$.} Fix $t\geq 0$ and let $\theta_n \to \theta \in \Theta$. The vector field $v$ is smooth in both variables and the flow $\phi$ is continuous. Hence, as $\theta_n \to \theta$ the smallest number $\delta_n>0$ that satisfies
        \begin{equation*}
            \norm{v(\phi^s \theta_n, x) - v(\phi^s \theta, x)}_e \leq \delta_n,
        \end{equation*}
        for all $s\in [0,t]$ and $x \in M$,
        approaches $0$. The statement of Lemma \ref{lem:uniform_regularity} implies that $\norm{\Per_{\theta_n}^t - \Per_\theta^t}\to 0$. This shows that $\theta \mapsto \Per_\theta^t$ is norm-continuous.
        \item \textit{strongly continuous.} Fix $f\in \H_0$. Let $\theta_n \to \theta \in \Theta$ and $t_n \to t \geq 0$. We compute
        \begin{equation} \label{eq:strong_continuity_sum}
            \norm{\Per_{\theta_n}^{t_n} f - \Per_\theta^tf} \leq \norm{\Per_{\theta_n}^{t_n} f - \Per_\theta^{t_n}f}  + \norm{\Per_{\theta}^{t_n} f - \Per_\theta^tf}. 
        \end{equation}
        The sequence $t_n$ is bounded from above by some value $T>0$. As $\theta_n \to \theta$, the smallest number $\delta_n>0$ that satisfies
        \begin{equation*} 
            \norm{v(\phi^s \theta_n, x) - v(\phi^s \theta, x)}_e \leq \delta_n,
        \end{equation*}
        for all $s\in [0,T]$ and $x \in M$,
        approaches $0$. By Lemma \ref{lem:uniform_regularity}, we find $\norm{\Per_{\theta_n}^{t_n} - \Per_\theta^{t_n}} \to 0$ as $n\to \infty$. In particular, the first summand in \eqref{eq:strong_continuity_sum} tends to $0$ as $n\to \infty$.

        Since $f(s) = \Per_\theta^s f$ is the solution to the NACP \eqref{eq:NACP} with starting parameter $\theta$, and $f(s)$ is continuous in $s$, the second summand in \eqref{eq:strong_continuity_sum} vanishes as well as $n\to \infty$. Hence, $(\theta, t) \mapsto \Per_\theta^t f$ is continuous.
        \item \textit{exponentially bounded.} By Proposition \ref{prop:spectral_gap}, we find $\norm{\Per_\theta^t}\leq e^{\varrho t}$, for some $\varrho<0$. Hence, the cocycle $\Per$ is exponentially bounded with $K=1$ and $L = \varrho$.
        
    \end{enumerate}
This finishes the proof.
\end{proof}

Since the Perron--Frobenius cocycle $\Per$ satisfies assumptions (I)--(IV), the theory of the previous sections is applicable. Let $\mathfrak F$ be either $L^2(\Theta, \H_0)$ or $C(\Theta, \H_0)$. The Mather semigroup $\mathbf{M}^t:\mathfrak F \to \mathfrak F$ defined by \eqref{eq:def_Mather} is a strongly continuous semigroup with closed generator $\G$. Since $\Per$ is the strongly continuous cocycle that is the solution to the NACP \eqref{eq:NACP}, we apply \cite[Example 6.22, Proposition 6.23]{chicone1999evolution}\footnote{In \cite{chicone1999evolution}, the statement is formulated for the Mather semigroup on $C(\Theta, \H_0)$. However, the proof is directly transferable to $L^2(\Theta, \H_0)$.} to compute the augmented generator~$\G$:
\begin{nalign}\label{eq:Mather_generator}
    [\G \f] (\theta) = G(\theta)\f(\theta) - [\mathbf{d} \f](\theta),
\end{nalign}
where $\mathbf{d}$ is the generator of the translation semigroup on $\mathfrak F$ in the direction of $\Psi$. The operator is defined on its maximal domain and given by
\begin{equation*}
    [\mathbf{d}\f](\theta) = \partial_t \f(\phi^t \theta) |_{t=0} = \partial_\theta \f(\theta) \cdot \Psi(\theta).
\end{equation*}
We note that for $\mathfrak F = L^2(\Theta, \H_0)$ the derivative is in the weak sense. The generator $\G$ is defined on
\begin{equation}\label{eq:Mather_generator_domain}
    \mathcal{D}(\mathcal{\G}) := \{\f \in \mathfrak F \mid \f \in \mathcal{D}(\mathbf{d}),\: \f:\Theta \to \mathcal{D}(G), \: G(\cdot)\f(\cdot) - \mathbf{d}\f \in \mathfrak F\}.
\end{equation}
The key property of this domain which will be used later is the fact that for any function $\f\in \mathcal{D}(\G)$ the fibres $\f(\theta)$ lie in $\mathcal{D}(G)$. An explicit form of the augmented generator is given by
\begin{equation}\label{eq:Mather_generator_explicit}
    [\G \f](\theta, x) = \frac{1}{2} \varepsilon^2 \Delta_x \f(\theta, x) - \nabla_x \f(\theta, x) \cdot v(\theta, x) - \partial_\theta \f (\theta,x) \cdot \Psi(\theta).
\end{equation}
By Theorem \ref{thm:spectral_mapping}, the spectrum of $\G$ is given by 
\begin{equation*}
    \sigma(\G) = \{\lambda + \eta i \: | \: \lambda \in \Sigma(\Per) , \: \eta \in \R \}.
\end{equation*}
Since $\norm{\Per_\theta^t}\leq e^{\varrho t}$ for all $\theta\in \Theta$ and $t\geq 0$, we find $\Sigma(\Per) \subset (-\infty, \varrho ]$. 
In general, the existence of eigenvalues of $\G$ cannot be guaranteed. The existence of so-called approximate eigenfunctions is discussed in Section \ref{sec:coherent_from_approx}. 

\begin{remark}
    The space $L^2(\Theta, \H)$ is isomorphic to $L^2(\Theta\times M, \C)$. The Mather operators $\mathbf{M}^t$, interpreted as operators on the space $L^2(\Theta\times M, \C)$ are, in fact, transfer operators of the dynamics on $\Theta \times M$ defined in \eqref{eq:SDE}. Hence, the Mather semigroup on $L^2(\Theta, \H_0)$ can be interpreted as the semigroup of transfer operators of the autonomous (skew product) SDE \eqref{eq:SDE} on $\Theta \times M$ in which the noise only acts in the dimensions of $M$.
\end{remark}

\section{Coherent sets}\label{sec:coherent_sets}
Given a nonautonomous particle flow like \eqref{eq:SDE}, there are regions in space that experience more mixing than other regions. Due to the nonautonomous nature of the system, these regions are in general time-dependent. We call a time-evolving region that experiences little mixing a \emph{coherent set}. We characterize coherent sets by the property that particles that start inside a coherent set stay inside the coherent set for a long time. Given a starting parameter $\theta \in \Theta$, consider a family of sets $A_\theta^\bullet:=(A_\theta^t)_{t\geq 0}$, with $A_\theta^t \subset M$.
We say that $A_\theta^\bullet$ is coherent if the survival probability
\begin{equation}
    \P_{x_0 \sim A_\theta^0}(x_s \in A_\theta^s, \: \forall s\in [0, t])
\end{equation}
decays slowly, or alternatively we call $A_\theta^\bullet$ a coherent set. The initial condition $x_0 \sim A_\theta^0$ asserts that $x_0$ is uniformly distributed on the set~$A_\theta^0$. This is not a rigorous definition of a coherent set until we have made more precise what slow decay means. To this end,
we introduce two quantities that quantify the coherence of~$A_\theta^\bullet$. The first notion is the \emph{escape rate}\footnote{In the literature, the escape rate is sometimes defined using the~$\limsup$. Hence, our definition corresponds to the ``pessimistic'' escape rate.} 
\begin{equation}
    \label{eq:def_escape_rate}
    E(A_\theta^\bullet) := \liminf_{t\to \infty} \frac{1}{t} \log\left( \P_{x_0 \sim A_\theta^0}(x_s \in A_\theta^s, \: \forall s\in [0, t]) \right),
\end{equation}
which considers the asymptotic decay rate of the survival probability. By definition, the escape rate is non-positive. An escape rate close to $0$ corresponds to a highly coherent set. The escape rate $E(A_\theta^\bullet)$ has been studied for periodically driven systems in~\cite{froyland2017estimating}. Note that the escape rate does not guarantee coherence on any finite time horizon. It might happen that the survival probability quickly approaches 0, before transitioning to an exponential decay at rate $E(A_\theta^\bullet)$. For this reason, we introduce another quantity of coherence, the \emph{cumulative survival probability}
\begin{equation}
    \label{eq:def_cumul_survival}
    C(A_\theta^\bullet) := \int_0^\infty  \P_{x_0 \sim A_\theta^0}(x_s \in A_\theta^s, \: \forall s\in [0, t])\: dt.
\end{equation}

The goal of this section is to elaborate how coherent sets $A_\theta^\bullet$ with an escape rate $E(A_\theta^\bullet)$ close to $0$ or a high cumulative survival probability~$C(A_\theta^\bullet)$ can be extracted from spectral objects associated with $\mathbf{M}^t$ or the cocycle~$\Per$. Ultimately, we want to compute coherent sets $A_\theta^\bullet$ for every starting parameter $\theta\in \Theta$ simultaneously, since $\mathbf{M}^t$ allows to simultaneously evolve functions $f_{\theta}$ by all $\Per_{\theta}^t$, for all $\theta\in\Theta$, respectively. Hence, we are looking for a family of sets $A_\bullet^\bullet := (A_\theta^t)_{\theta \in \Theta, t\geq 0}$ such that $A_\theta^\bullet := (A_\theta^t)_{t\geq 0}$ is coherent for each $\theta \in \Theta$. In some cases it is possible to choose the sets $A_\theta^t$ such that $\smash{ A_\theta^t=A_{\phi^t \theta}^0 }$ for each $\theta \in \Theta$ and $t\geq 0$. Then, the family $A_\bullet^\bullet$ is fully characterized by the sets $A_\theta^0$ for $\theta \in \Theta$. If this is the case, we omit the superscript $0$ and simply write $A_\bullet = (A_\theta)_{\theta \in \Theta}$. Such a representation is particularly nice in practice, since once the sets $A_\theta$ are determined for all $\theta \in \Theta$, one obtains a coherent set for any starting parameter~$\theta$.

In \cite{froyland2017estimating} it was shown that for fixed $\theta \in \Theta$, coherent sets can be extracted from functions $f\in \H_0$ that decay slowly under $\Per$. This construction requires the point-wise evaluation of $\Per_\theta^t f$ in $x$. Since we work with $L^2$-functions, point-wise evaluations are a priori not well-defined. However, we know that $\Per_\theta^t f \in \mathcal{D}(G) \subset H^2(M)$ such that there is a unique continuous representative of $\Per_\theta^t f$. We state an even stronger regularity result. To apply a suited Sobolev embedding, the dimension of the physical space $M$ should not exceed~3.
\begin{proposition}\label{prop:continuous_regularity}
    Assume that the dimension $d$ of $M$ is at most 3. Let $f_0\in \mathcal{D}(G)$ and let $\theta \in \Theta$. Then $f:\R_{\geq 0}  \to \H_0$, defined by $f(t) = [\Per_\theta^t f_0]$, is represented by a continuous function $f \in C(\R_{\geq 0}\times M, \C)$.
\end{proposition}
\begin{proof}
    See \cite[Theorem A.4]{froyland2020computation}. The assumption $d\leq 3$ is needed such that the Sobolev embedding $W^{2,2}(M) \subset C^\alpha(M)$ holds for $\alpha = 2 - \frac{d}{2}$, where $C^\alpha(M)$ is the usual Hölder space. 
\end{proof}

The following theorem is almost identical to \cite[Theorem 19]{froyland2017estimating} with the difference that the latter result requires slightly stronger assumptions.
\begin{theorem}
    \label{thm:CohFamFromFun}
    Let $f_0\in \mathcal{D}(G)$ be a real function with $\int_M f_0(x) dx=0$ and let $\theta \in \Theta$. Define the sets 
    \begin{equation}\label{eq:coheren_set_in_thm}
        A^t_\theta := \{x \in M\: | \: [\Per_\theta^t f_0] (x) \geq 0\}.
    \end{equation}
    Then, for any $t\geq 0$, we find
    \begin{equation}\label{eq:FK17_estimate}
        \P_{x_0 \sim A_\theta^0}(x_s \in A_\theta^s, \: \forall s\in [0, t]) \geq \frac{1}{2}\norm{f_0}_\infty^{-1} \abs{A_\theta^0}^{-1} \norm{\Per_\theta^t f_0}_1.
\end{equation}
\end{theorem}
\begin{proof}
    See Appendix \ref{ap:proofs}.
\end{proof}

Recall that $\mathcal{D}(G)\subset H^2(M)$, cf.~\eqref{eq:generator_domain}, such that $\norm{f_0}_\infty$ is finite. Note that the theorem above bounds the decay of the survival probability by the $L^1$-norm of $\Per_\theta^t f$. Since we work on $\H_0 \subset L^2(M, \C)$, we need an additional estimation to express the decay in terms of the $L^2$-norm of $\Per_\theta^t f$. Bounding the $L^1$-norm $\norm{f}_1$ from below by $\norm{f}_\infty^{-1} \norm{f}_2^2$ via H\"older, we obtain
\begin{equation}\label{eq:suvival_prob_L2}
     \P_{x_0 \sim A_\theta^0}(x_s \in A_\theta^s, \: \forall s\in [0, t]) \geq \frac{1}{2}\norm{f}_\infty^{-2} \abs{A_\theta^0}^{-1}  \norm{\Per_\theta^t f}_2^2.
\end{equation}
The bound on the the $L^1$-norm in terms of the $\infty$-norm and the square of the $L^2$-norm is in general rather crude. Hence, \eqref{eq:suvival_prob_L2} is expected to underestimate the survival probability of $A_\theta^\bullet$. In Section \ref{sec:examples}, we compute and evaluate coherent sets for some examples. Discretizing the system introduces errors that decrease the coherence of the computed coherent sets. The results of Section \ref{sec:examples} indicate that the discretization error and the crude bound on the $L^1$-norm seem to cancel, in the sense that the simulated survival probabilities decay similarly to the right-hand side of \eqref{eq:suvival_prob_L2}.
\begin{remark}
    Theorem \ref{thm:CohFamFromFun} establishes a connection between the decay rate of a function $f\in C(M,\R)$ under $\Per$ and the escape rate of the corresponding coherent set $A_\theta^\bullet$ defined in \eqref{eq:coheren_set_in_thm}. We find
    \begin{align*}
        E(A_\theta^\bullet) &\geq \liminf_{t\to \infty} \frac{1}{t} \log\left(\frac{1}{2}\norm{f}_\infty^{-1} \abs{A_\theta^0}^{-1} \norm{\Per_\theta^t f}_1\right) \\
        &= \liminf_{t\to \infty} \frac{1}{t} \log\left(\norm{\Per_\theta^t f}_1\right).
    \end{align*}
    This quantity is called the \emph{Lyapunov exponent} (in the $L^1$-norm) of $f$ under $\Per$. The study of Lyapunov exponents to characterize coherent sets has been pioneered in~\cite{froyland2010coherent}.
\end{remark}

\subsection{Coherent sets from eigenfunctions}
In the previous works \cite{froyland2017estimating, froyland2020computation}, coherent sets were extracted via Theorem~\ref{thm:CohFamFromFun} from eigenfunctions of an augmented generator, as these provided functions that decayed slowly. In our setting, the augmented generator $\G$ and the Mather operators $\mathbf{M}^t$ do in general not possess eigenfunctions. However, in Section \ref{sec:numerics} we construct finite-dimensional approximations of the linear operator $\G$ and the Mather operators~$\mathbf{M}^t$. Naturally, these matrices admit eigenfunctions. In the following paragraph we assume the existence of an eigenfunction $\f$ of the Mather semigroup and its generator over $C(\Theta, \H_0)$ and show how to construct coherent sets from~$\f$. Hence, this subsection should be understood as the derivation of a heuristic method that can be applied to eigenfunctions of the discrete generator.

Assume there was an eigenfunction $\f\in C(\Theta, \H_0)$ of the augmented generator $\G$ with eigenvalue $z = \lambda + \eta i \in \C$. Consequently, $\f$ is an eigenfunction of $\mathbf{M}^t$ with eigenvalue $e^{zt}$ for every $t\geq 0$. In particular, we find
\begin{equation}\label{eq:eigenvalue_equation}
    \Per_\theta^t \f(\theta) = e^{z t} \f(\phi^t \theta),
\end{equation}
for all $\theta \in \Theta$ and $t\geq 0$. By Theorem \ref{thm:classical_solution}, we find $\f^R(\theta)\in \mathcal{D}(G)$ for all $\theta \in \Theta$. Let $\f^R:=\text{Re}(\f)$ be the real part of $\f$. 
Since $\Per$ is a real cocycle, \eqref{eq:eigenvalue_equation} implies
\begin{equation*}
   \Per_\theta^t \f^R(\theta) = e^{\lambda t}\text{Re}(e^{\eta t i}  \f(\phi^t\theta)).
\end{equation*}
We define a coherent family $A_\bullet^\bullet = (A_\theta^t)_{\theta \in \Theta, t\geq 0}$ by
\begin{equation}\label{eq:positive_level}
    A^t_\theta := \{x\in M \:|\: \text{Re}(e^{\eta t i}  \f(\phi^t\theta, x)) \geq 0\}.
\end{equation}
The estimate~\eqref{eq:suvival_prob_L2} yields
\begin{equation}\label{eq:eigenfunction_survival_prob}
    \P_{x_0 \sim A_\theta^0}(x_s \in A_\theta^s, \: \forall s\in [0, t]) \geq \frac{1}{2}e^{2\lambda t} \norm{\f^R(\theta)}_\infty^{-2} \abs{A_\theta^0}^{-1} \norm{\text{Re}(e^{\eta t i}  \f(\phi^t\theta))}_2^2.
\end{equation}
Assuming that the real part of the fibres $e^{\eta t i}  \f(\phi^t\theta)$ is uniformly bounded from below, we conclude that the survival probability of $A_\theta^\bullet$ decays at most at rate $e^{2\lambda t}$, up to constants. Hence, we expect
\begin{equation}\label{eq:cum_prob_estimate}
    E(A_\theta^\bullet) \geq 2\lambda, \quad C(A_\theta^\bullet) \gtrsim \int_0^\infty e^{2\lambda t} dt = -\frac{1}{2\lambda}.
\end{equation}
The results of Section \ref{sec:examples} show that $-\frac{1}{2\lambda}$ seems to be a good estimate of the cumulative survival probability that we computed in numerical examples.

Using \eqref{eq:positive_level} to define a coherent set enables us to derive theoretical bounds on the survival probability, but the resulting coherent sets might not provide much insight into the dynamical properties of the underlying vector field. The definition simply divides the space $M$ into two regions of roughly the same size. Considering the function $-\mathbf{f}$ shows that the complement $B_\theta^t := M \setminus A_\theta^t$ satisfies the same bounds on the survival probability. Hence, the method \eqref{eq:positive_level} is unable to identify multiple smaller coherent regions.
It was, however, extended in \cite[Proposition 4.6]{froyland2020computation} to use linear combinations of multiple eigenfunctions to obtain multiple coherent sets from them.
Here, in Section \ref{sec:numerics_coherent_sets} we propose alternative methods to identify coherent sets from the fibres of an eigenfunction~$\f$.

\subsection{Coherent sets from approximate eigenfunctions}\label{sec:coherent_from_approx}
In general, the augmented generator $\G$ does not possess eigenfunctions. However, the existence of \emph{approximate eigenfunctions} is guaranteed. While less applicable in practice, we show how approximate eigenfunctions can be utilized to construct coherent sets. We briefly introduce the approximate spectrum of an operator.

Let $\mathds B$ be a Banach space and let $T:\mathcal{D}(T) \to \mathds B$ be a closed operator, where $\mathcal{D}(T) \subset \mathds B$ is the domain of~$T$. The \emph{approximate spectrum} of $T$ is defined by
\begin{equation*}
    \sigma_\text{ap}(T) = \{z \in \C \mid \exists x_n \in \mathcal{D}(T) \text{ with } \norm{x_n}=1 : \: \norm{(T-z ) x_n} \to 0, \text{ as }n\to \infty\}.
\end{equation*}
We call the sequence $(x_n)_{n\in \N}$, and any element of it, an \emph{approximate eigenfunction} to the approximate eigenvalue~$z$. We find $\partial \sigma(T) \subset  \sigma_\text{ap}(T) \subset \sigma(T)$, where $\partial \sigma(T)$ is the boundary of the spectrum of~$T$. For a proof of this fact, see \cite[Chapter IV, Proposition 1.10]{engel2000one}.

Recall Theorem \ref{thm:spectral_mapping}, which characterized the spectrum of the Mather operators $\mathbf{M}^t$ and augmented generator $\G$. In particular, the spectrum of $\G$ is given by
\begin{equation}\label{eq:spectrum_G_reloaded}
    \sigma(\G) = \{\lambda + \eta i \mid \lambda \in \Sigma(\Per) , \: \eta \in \R \}.
\end{equation}
As mentioned at the end of Section \ref{sec:transfer_operator}, the Sacker--Sell spectrum $\Sigma(\Per)$ is contained in $(-\infty, \varrho]$ for the constant $\varrho < 0$ from Proposition \ref{prop:spectral_gap}. Assume that $\Sigma(\Per) \neq \emptyset$, and let $\lambda<0$ be its maximum value. By \eqref{eq:spectrum_G_reloaded}, we find $\lambda \in \partial \sigma(\G)$, and, therefore, $\lambda \in \sigma_{\text{ap}}(\G)$. This holds independently of whether we consider the Mather semigroup and its generator on $L^2(\Theta, \H_0)$ or $C(\Theta, \H_0)$. 

Since $\G$ is real and $\lambda \in \R$, we find
\begin{align*}
    \norm{(\G-\lambda)\f_n}_\infty &\geq \max\{ \norm{\text{Re}((\G-\lambda)\f_n)}_\infty, \:\norm{\text{Im}((\G-\lambda)\f_n)}_\infty \}\\
    &= \max \{ \norm{(\G-\lambda)\f_n^R}_\infty, \: \norm{(\G-\lambda)\f_n^I}_\infty\},
\end{align*}
where $\f_n^R$ and $\f_n^I$ are the real and imaginary parts of~$\f_n$. Note that both $\f_n^R$ and $\f_n^I$ are real functions with $\norm{\f_n}_\infty \leq \norm{\f_n^R}_\infty + \norm{\f_n^I}_\infty$.
Hence, w.l.o.g.~we can choose an approximate eigenfunction $(\f_n)_{n\in \N}$ such that $\f_n \in C(\Theta, L^2(M, \R))$ are real functions. In particular, since $\f_n \in \mathcal{D}(\G)$, we have that $\f_n(\theta) \in \mathcal{D}(G)$ for every $\theta \in \Theta$, cf.~\eqref{eq:Mather_generator_domain}. 
\begin{proposition}
\label{prop:CohFamfromApprox}
   For any $T>0$ and $\delta >0$, there is an $n \in \N$ such that the family of sets defined by
    \begin{equation}\label{eq:coherent_set_from_approx}
        A^t_\theta := \{x \in M\: | \: [\Per_\theta^t \f_n (\theta)] (x) \geq 0\}
    \end{equation}
    satisfies the following estimate for all $\theta \in \Theta$ and $t\in[0, T]$:
    \begin{equation}
    \label{eq:survival_from_ap_efun}
    \P_{x_0 \sim A_\theta^0}(x_s \in A_\theta^s, \: \forall s\in [0, t]) \geq \frac{1}{2}e^{2\lambda t}\norm{\f_n(\theta)}_\infty^{-2} \abs{A_\theta^0}^{-1} \big(\norm{\f_n(\phi^t \theta)}_2 - \delta\big)^2.
\end{equation}
\end{proposition}
\begin{proof}
    The unbounded, closed operator $(\G-\lambda)$ generates the strongly continuous semigroup~$e^{-\lambda}\mathbf{M}^t$. Since $\f_n \in \mathcal{D}(\G-\lambda)$, we apply \cite[Theorem 2.4]{pazy2012semigroups} to compute
    \begin{equation*}
        e^{-\lambda t} \mathbf{M}^t \f_n - \f_n = \int_0^t \mathbf{M}^s (\G - \lambda) \f_n \:ds.
    \end{equation*}
    We use the fact $\norm{\mathbf{M}^s}\le e^{\varrho s} \le 1$, by Proposition \ref{prop:spectral_gap}, to obtain
    \begin{equation*}
        \norm{ e^{-\lambda t} \mathbf{M}^t \f_n - \f_n}_\infty \leq t \norm{(\G - \lambda) \f_n}_\infty.
    \end{equation*}
    As $n\to \infty$, the right-hand side approaches $0$. For fixed $T>0$ and $\delta>0$ there is an $n\in \N$ such that
    \begin{equation}\label{eq:approximate_eFun_discrepancy}
        \norm{ e^{-\lambda t} \Per_\theta^t \f_n(\theta) - \f_n(\phi^t \theta)}_2 \leq \delta
    \end{equation}
    for all $\theta \in \Theta$ and $t\in [0, T]$.
    Rearranging the inequality yields
    \begin{equation}\label{eq:decay_approx_efun}
        \norm{\Per_\theta^t \f_n(\theta)}_2 \geq e^{\lambda t} \big(\norm{\f_n(\phi^t \theta)}_2 - \delta\big).
    \end{equation}
    Now, the statement follows from \eqref{eq:suvival_prob_L2}.
\end{proof}
Note that the bound \eqref{eq:survival_from_ap_efun} bound is almost identical to \eqref{eq:eigenfunction_survival_prob}, the bound on the survival probability of a coherent set extracted from an eigenfunction of $\G$, up the error term $\delta$. 

In conclusion, a coherent sets $A_\theta^\bullet$ extracted from an approximate eigenfunction is expected to have a similar survival probability as a coherent set extracted from an eigenfunction of~$\G$. However, the bounds on the survival probability of $A_\theta^\bullet$ only hold up to a finite time-horizon $T$. Hence, we cannot make any statement about the asymptotic decay rate $E(A_\theta^\bullet)$, but we expect the cumulative survival probability to be of order~$-\frac{1}{2 \lambda}$. Another drawback of extracting coherent sets from an approximate eigenfunction versus from an eigenfunction is that the definition of $A_\theta^\bullet$ in \eqref{eq:coherent_set_from_approx} requires knowledge of the cocycle evaluations $\Per_\theta^t \f_n (\theta)$, unlike the definition \eqref{eq:positive_level} of a coherent set extracted from an eigenfunction $\f$, which only required knowledge of the fibres of~$\f$.
 
\subsection{Coherent sets from spectral projections}
In the previous section we studied coherent sets extracted from approximate eigenfunctions of the augmented generator $\G$. However, we were only able to derive bounds on the decay rate up to a finite time horizon $[0, T]$. In this section we show how coherent sets can be extracted from spectral subspaces, given the existence of exponential dichotomies. The advantage of this approach is that we are able to obtain a uniform decay rate for $t\to \infty$. In particular, we are able to derive a bound on the asymptotic escape rate $E(A_\theta^\bullet)$.

By Proposition \ref{prop:Sacker_Sell_segments} the Sacker--Sell spectrum $\Sigma(\Per)$ consists of the union of closed segments $[r_k^-, r_k^+]\subset \R$ for $1\leq k \leq N$ and $N\in \N_0 \cup \{\infty\}$. 
Assume that $\Sigma(\Per) \neq \emptyset$ and that $r_1^- > -\infty$, i.e.~the interval $[r_1^-, r_1^+]$ is an isolated segment of~$\Sigma(\Per)$. Fix $\lambda < r_1^-$ close to $r_1^-$. There is an exponential dichotomy at $\lambda$ that splits the space $\H_0$ into a stable bundle $S(\theta)$ and an unstable bundle~$U(\theta)$. 
There is a constants $\beta>0$ and $C>0$ such that
\begin{equation*}
    \norm{\Per_\theta^t f}_2 \geq C^{-1} e^{(\lambda + \beta) t} \norm{f}_2,
\end{equation*}
for all $\theta \in \Theta$, $t\geq 0$ and $f\in U(\theta)$.
The next lemma shows that w.l.o.g.~we can choose $f\in U(\theta)$ to be a real function.
\begin{lemma}\label{lem:projection_real}
    If $f\in U(\theta)$, then both $f^R:=\textup{Re}(f)$ and $f^I:=\textup{Im}(f)$ lie in $U(\theta)$ as well.
\end{lemma}
\begin{proof}
    See Appendix \ref{ap:proofs}.
\end{proof}

Let $f\in U(\theta)$ be a real function with $\norm{f}_2=1$. By the definition of an exponential dichotomy, we find $f\in \text{ran}(\Per_{\phi^{-t} \theta}^t)$ for some $t>0$, and thereby, $f\in \mathcal{D}(G)$. Define a coherent set by
\begin{equation}
    A^t_\theta := \{x \in M\: | \: [\Per_\theta^t f] (x) > 0\}.
\end{equation}
By \eqref{eq:suvival_prob_L2}, the survival probability of this coherent set satisfies
\begin{equation}
     \P_{x_0 \sim A_\theta^0}(x_s \in A_\theta^s, \: \forall s\in [0, t]) \geq \frac{1}{2}\norm{f}_\infty^{-2}\abs{A_\theta^0}^{-1} C^{-1} e^{(\lambda + \beta)t}.
\end{equation}
This bound provides an estimate on the escape rate defined in \eqref{eq:def_escape_rate}:
\begin{equation}
    E(A_\theta^\bullet) \geq \lambda + \beta.
\end{equation}
Since $\lambda < r_1^-$ was arbitrary, we find coherent sets with escape rate arbitrarily close to $r_1^-$.

In Theorem \ref{thm:spectral_mapping}, it was shown that for any $t\geq 0$ the spectrum of the Mather operators $\mathbf{M}^t$ contains the isolated annulus $\Lambda^t:=\{e^{t\lambda + \eta i} \:|\: \lambda \in [r^-_1, r^+_1], \: \eta \in [0, 2\pi)\}$. By the contents of Section \ref{sec:spectral_projections}, in particular \eqref{eq:Riesz_range}, the range of the Riesz-projection $\mathfrak{P}_{\Lambda^t}$ consists of functions with fibres in $U(\theta)$.
Hence, given a function $\f \in \ran(\mathfrak{P}_{\Lambda^t})$, the nonzero fibres $\f(\theta)\in U(\theta)$ induce coherent sets with escape rate at least~$\lambda + \beta$ (note that due to the negative sign a larger escape rate corresponds to higher coherence).

\section{Numerical implementation}
\label{sec:numerics}

\subsection{The augmented generator in Fourier space}
\label{ssec:genFourier}
The theoretical results of the previous two section were valid for domains $M$ that were either open, bounded subsets of $\R^d$ with piecewise $C^4$ boundary, or 
a torus $\tor^d$. In the following, to enable the use of Fourier series, we consider the case $\Theta = \tor^{d_d}$ and~$M = \tor^{d_p}$. The dimensions $d_d$ and $d_p$ stand for the driving dimension and the physical dimension, respectively. For the ergodic driving dynamic we choose a quasi-periodic rotation on the torus, i.e.
\begin{equation*}
    \partial_t \theta_t = \Psi(\theta_t) = \alpha,
\end{equation*}
where $\alpha \in \R^{d_d}$ is a constant vector whose entries are rationally independent. We derive a finite-dimensional approximation of the Mather semigroup $\mathbf{M}^t$ and its generator $\G$. 

Note that the space $L^2(\Theta, \H)$ is isomorphic to $L^2(\Theta \times M, \C) = L^2(\tor^{d_d} \times \tor^{d_p}, \C)$ which is a Hilbert space with inner product
\begin{equation*}
    \langle \f, \mathbf{g} \rangle = \int_{\tor^{d_d}} \int_{\tor^{d_p}} \overline{\f(\theta, x)} \mathbf{g}(\theta, x) \: dx \: d\theta. 
\end{equation*}
An orthonormal basis, also called Hilbert basis, is given by the Fourier modes
\begin{equation*}
    F_{m,n}(\theta, x):= e^{2\pi i (\theta \cdot m + x \cdot n)},
\end{equation*}
for $m\in \Z^{d_d}$ and $n \in \Z^{d_p}$. We represent functions $\f\in L^2(\tor^{d_d} \times \tor^{d_p}, \C)$ with respect to this basis, i.e.~
\begin{equation*}
    \f = \sum_{m,n} \widehat{\f}(m,n) F_{m,n},
\end{equation*}
where the summation is over all $(m,n)\in \Z^{d_d}\times \Z^{d_p}$.
The Fourier coefficients $\widehat{\f}(m,n) \in \C$ can be computed by
\begin{equation*}
    \widehat{\f}(m,n) = \langle F_{m,n}, \f \rangle = \int_{\tor^{d_d}} \int_{\tor^{d_p}} e^{-2\pi i (\theta \cdot m + x \cdot n)} \f(\theta, x) \:dx\: d\theta.
\end{equation*}
Analogously, we can represent the vector field $v:\tor^{d_d} \times \tor^{d_p} \to \R^{d_p}$ in Fourier coordinates. Each of the $d_p$ components $v_j$ has its individual Fourier coefficients $\widehat{v_j} (m,n)$. We write $\widehat{v}(m,n) \in \C^{d_p}$ for the vector of Fourier coefficients.
\begin{lemma}\label{lem:div_free_fourier}
    Let $v:\tor^{d_d} \times \tor^{d_p} \to \R^{d_p}$ be a divergence-free vector field. For any $(m,n) \in \Z^{d_d} \times \Z^{d_p}$, we find
    \begin{align*}
        \widehat{v}(m,n) &= \overline{\widehat{v}(-m,-n)}, \\
        n \cdot \widehat{v}(m,n)&= 0.
    \end{align*}
\end{lemma}
\begin{proof}
    See Appendix \ref{ap:proofs}.
\end{proof}

Similar to how we write functions $\f\in L^2(\tor^{d_d} \times \tor^{d_p}, \C)$ in terms of their Fourier coefficients, we can be represent the augmented generator $\G$ as an infinite-dimensional matrix with rows and columns indexed by $\Z^{d_d}\times \Z^{d_p}$, respectively. Recall that we defined the Mather semigroup $(\mathbf{M}^t)_{t\geq 0}$ and its generator $\G$ over functions $\f$ with fibres $\f(\theta)$ in the Hilbert space $\H_0$ of $L^2$-functions that integrate to~$0$. This corresponds to only considering those Fourier modes $F_{m,n}$ for which $n\in \Z^{d_p}$ is nonzero. To avoid additional case distinctions, in this section we include functions that do not integrate to~$0$. We discuss this choice further in Remark \ref{rem:zero_modes}.
\begin{proposition}\label{prop:matrix_Gamma}
    Let $\f \in \mathcal{D}(\G)$. We find
    \begin{equation*}
        \widehat{\G \f}(m,n) = \sum_{m', n'} \Gamma(m, n \: ; \: m', n') \widehat{\f}(m',n').
    \end{equation*}
    The infinite-dimensional matrix $\Gamma \in \C^{(\Z^{d_d}\times \Z^{d_p}) \times (\Z^{d_d}\times \Z^{d_p})}$ is given by
    \begin{equation}\label{eq:def_Gamma}
        \Gamma(m, n \: ; \: m', n') := \begin{cases}
            \frac{1}{2} \varepsilon^2 (2\pi \norm{n}_e)^2 - 2\pi i \big(m \cdot \alpha + n\cdot\widehat{v}(0,0) \big), \quad & (m,n) = (m',n')\\
            -2\pi i (n \cdot \widehat{v}(m-m', n-n')), & (m,n) \neq (m',n').
        \end{cases}
    \end{equation}
\end{proposition}
\begin{proof}
    Let $\f \in \mathcal{D}(\G)$. Since the operator $\G$ is closed, we find
    \begin{equation*}
        \G \f = \sum_{m',n'} \widehat{\f}(m',n') \G F_{m',n'}\,.
    \end{equation*}
    We conclude
    \begin{equation}\label{eq:GF_matrix_product}
        \widehat{\G \f}(m,n) = \sum_{m', n'} \widehat{\G F_{m',n'}}(m,n) \widehat{\f}(m',n').
    \end{equation}
    In the following we show $\widehat{\G F_{m',n'}}(m,n) = \Gamma(m, n \: ; \: m', n')$, which will complete the proof. Using \eqref{eq:Mather_generator_explicit}, we explicitly compute~$\G F_{m',n'}$.
    \begin{align*}
        \G F_{m',n'}(\theta, x) &= \frac{1}{2} \varepsilon^2 \Delta_x F_{m',n'}(\theta, x) - \nabla_x F_{m',n'}(\theta, x) \cdot v(\theta, x) - \partial_\theta F_{m',n'}(\theta, x) \cdot \Psi(\theta) \\
        &= \left(\frac{1}{2} \varepsilon^2 (2\pi \norm{n'}_e)^2 - 2\pi i (n'\cdot v(\theta, x)) - 2\pi i (m' \cdot \alpha)\right) F_{m',n'}(\theta, x).
    \end{align*}
    From this formula, we can extract the Fourier coefficients, which are made up of three summands
    \begin{nalign} \label{eq:GF_threee_summands}
        \widehat{\G F_{m',n'}}(m,n) &= \langle F_{m,n} , \G F_{m',n'} \rangle \\
        &=\frac{1}{2} \varepsilon^2 (2\pi \norm{n'}_e)^2 \langle F_{m,n}, F_{m',n'} \rangle \\
        &\hphantom{=}- 2\pi i (m' \cdot \alpha) \langle F_{m,n}, F_{m',n'} \rangle\\
        &\hphantom{=}- 2 \pi i \langle F_{m,n} , (n' \cdot v) F_{m', n'} \rangle.
    \end{nalign}
    Since the Fourier functions form an orthonormal basis, $\langle F_{m,n}, F_{m',n'} \rangle$ is 1 if $(m,n) = (m', n')$ and $0$ otherwise. Therefore, the first two summands are $0$ unless $(m,n) = (m',n')$. In that case, we can replace $m'$ and $n'$ by $m$ and $n$. We compute the third summand
    \begin{nalign}\label{eq:third_summand}
        - 2 \pi i \langle F_{m,n} , (n' \cdot v) F_{m', n'} \rangle &= - 2\pi i \langle F_{m - m', n - n'}, n' \cdot v \rangle \\
        &= - 2\pi i (n' \cdot \langle F_{m - m', n - n'}, v \rangle) \\
        &= -2\pi i( n' \cdot \widehat{v} (m-m', n-n')) \\
        &= -2\pi i( n \cdot \widehat{v} (m-m', n-n')),
    \end{nalign}
    where in the last line we used Lemma \ref{lem:div_free_fourier}. Inserting \eqref{eq:third_summand} into \eqref{eq:GF_threee_summands} shows that $\widehat{\G F_{m',n'}}(m,n) = \Gamma (m,n\: ; \: m', n')$ as defined in \eqref{eq:def_Gamma}. 
\end{proof}

\begin{remark}\label{rem:zero_modes}
    An off-diagonal entry of $\Gamma(m,n\: ; \: m', n')$ is nonzero if and only if $n \cdot \widehat{v}(m-m', n-n') \neq 0$. By Lemma \ref{lem:div_free_fourier}, this is equivalent to $n' \cdot \widehat{v}(m-m', n-n') \neq 0$. Hence, if either $n=0$ or $n'=0$, then $\Gamma(m,n\: ; \: m', n') = 0$. This shows that Fourier modes $F_{m,n}$ with $n=0$ do not ``communicate'' with other Fourier modes under~$\Gamma$. This justifies including these Fourier modes in $\Gamma$, even though in~Section \ref{sec:transfer_operator} the augmented generator $\G$ was constructed over functions $\f$ whose fibres $\f(\theta)$ lie in $\H_0$. In particular, any Fourier mode $F_{m,0}$ is an eigenfunction of $\G$ with purely imaginary eigenvalue (essentially owing to the constant function in physical space being invariant under $\Per^t_{\theta}$ and Fourier modes on the base being eigenfunctions of the ergodic torus rotation). We include these Fourier modes to visualize the spectral gap between $\varrho$ and $0$ in the examples of Section~\ref{sec:examples}.
\end{remark}

\subsection{The discrete augmented generator}
\label{ssec:discr_gen}

We discretize the generator $\G$ by a Galerkin projection onto a chosen set of Fourier modes, e.g.~those with $m,n$ close to $0$. Let $S_d \subset \Z^{d_d}$ and $S_p \subset \Z^{d_p}$ be a selection of Fourier modes in the driving space and in physical space, respectively. To have the discrete operator mapping real-valued functions to real-valued ones again, Lemma~\ref{lem:div_free_fourier} suggests to choose the sets $S_d$ and $S_p$ such that $m \in S_d$ if and only if $-m \in S_d$ and $n\in S_p$ if and only if~$-n \in S_p$.
We project onto the Fourier modes given by
\begin{equation*}
    S := S_d \times S_p \subset \Z^{d_d}\times \Z^{d_p}.
\end{equation*}
A natural choice for $S_d$, respective $S_p$, are all modes in $\Z^{d_d}$, respective $\Z^{d_p}$, that are within a given radius of the origin, with respect to either the Euclidean or $\ell^\infty$-norm. In two of the examples of Section \ref{sec:examples}, we fix constants $K \in \N$ and $r>0$ and define
\begin{equation}\label{eq:S_d_S_p}
    S_d := \{m \in \Z^{d_d} \mid \norm{m}_\infty \leq K\}, \quad S_p := \{n \in \Z^{d_p} \mid \norm{n}_e \leq r\}.
\end{equation}

Given a set $S\subset\Z^{d_d}\times \Z^{d_p}$, we obtain the discrete generator $\Gamma_S \in \C^{S \times S}$ as the restriction of the matrix $\Gamma$ from Proposition \ref{prop:matrix_Gamma} to the Fourier modes in~$S$. We can write $\Gamma_S$ as the sum of two diagonal and a skew-Hermitian matrix
\begin{equation}\label{eq:gamma_DRA}
    \Gamma_S = D + R + A.
\end{equation}
The matrices $D,R\in \C^{S\times S}$ are diagonal and given by
\begin{nalign}\label{eq:def_DR}
    D&=\text{diag}\left(\frac{1}{2} \varepsilon^2 (2\pi \norm{n}_e)^2 \mid (m,n) \in S\right),    \\
    R&=\text{diag}\left(-2 \pi i(m\cdot \alpha) \mid (m,n) \in S\right).
\end{nalign}
Note that $D$ is real while $R$ is purely imaginary. The matrix $A\in \C^{S\times S}$ is given by
\begin{equation}\label{eq:def_A}
    A(m,n\: ; \: m', n') = -2\pi i (n \cdot \widehat{v}(m-m', n-n')).
\end{equation}
Using Lemma \ref{lem:div_free_fourier}, one can verify that $A$ is skew-hermitian, i.e., $A(m',n'\: ; \: m, n) = -\overline{A(m,n\: ; \: m', n')}$.

The matrix $\Gamma_S\in \C^{S\times S}$ is the Galerkin projection onto the subspace
\begin{equation}
    V_S := \text{span}\{F_{m,n} \: | \: (m,n) \in S\}.
\end{equation}
The orthogonal projection $\Pi_S:L^2(\tor^{d_d} \times \tor^{d_p}, \C) \to L^2(\tor^{d_d} \times \tor^{d_p}, \C)$ onto the subspace $S$ can be expressed as a concatenation $\Pi_S=P_S^* P_S$. Define the linear operator $P_S: L^2(\tor^{d_d} \times \tor^{d_p}, \C) \to \C^S$ by
\begin{equation}
    (P_S \f)_{m,n} = \widehat{\f}(m,n),
\end{equation}
where $(m,n)\in S$. Let $\widehat{\f}_S \in \C^S$. The adjoint $P_S^*: \C^S \to L^2(\tor^{d_d} \times \tor^{d_p}, \C)$ is given by
\begin{equation}
    P_S^* \widehat{\f}_S = \sum_{(m,n)\in S} \widehat{\f}_S(m,n) F_{m,n}.
\end{equation}
Note that $P_S^*P_S$ is the identity on $\C^S$. The matrix $\Gamma_S$ is a Galerkin projection of augmented generator $\G$ in the sense that
\begin{equation*}
    \Gamma_S = P_S \G P_S^*.
\end{equation*}
For any $t\geq 0$ we obtain a finite-dimensional approximation of the Mather operator $\mathbf{M}^t$ given by~$\exp (t \Gamma_S)$.
In general, projection and exponentiation will not commute, hence $P_S \mathbf{M}^t P_S^* \neq \exp (t \Gamma_S)$.

\subsection{Coherent sets}\label{sec:numerics_coherent_sets}

Let $\Gamma_S\in \C^{S\times S}$ be a finite-dimensional approximation of the augmented generator $\G$, like introduced in the previous section. Let $\widehat{\f}_S \in \C^S$ be an eigenfunction of $\Gamma_S$ with eigenvalue $z = \lambda + \eta i \in \C$. Consequently, $\widehat{\f}_S$ is an eigenfunction of the discrete Mather operator $\exp(t \Gamma_S)$ with eigenvalue $e^{tz}$. Recall the decomposition $\Gamma_S = D+R+A$ from \eqref{eq:gamma_DRA}. Since $A+R$ is skew-Hermitian and $D$ is diagonal with real, non-positive entries, all eigenvalues of $\Gamma_S$ have non-positive real part. The unit vector $e_{(0,0)} \in \C^S$ corresponding to the Fourier mode $F_{0,0}$ is in the kernel of $\Gamma_S$. In fact, all unit vector $e_{(m,0)}\in \C^S$ are eigenfunctions of $\Gamma_S$ with purely imaginary eigenvalues. A finite-dimensional vector $\smash{\widehat{\f}_S \in \C^S}$ represents a function $\smash{ \f := P_S^* \widehat{\f}_S \in L^2(\tor^{d_d} \times \tor^{d_p}, \C) }$. For the computation of coherent sets, we are only interested in those functions $\f$ whose fibres $\f(\theta)$ integrate to $0$. Hence, we only consider eigenfunctions $\widehat{\f}_S$ of $\Gamma_S$ to eigenvalues with negative real part $\lambda<0$. We extract coherent sets from $\f$ as if it was an eigenfunction of~$\G$.

In Section \ref{sec:coherent_sets}, we presented methods to compute coherent sets from from functions $\f$ that are either (approximate) eigenfunction of the generator $\G$ or that lie in the spectral set corresponding to an annulus in $\sigma(\mathbf{M}^t)$ with large radius. In the finite-dimensional case, any approximate eigenvalue is an eigenvalue and every spectral set contains at least one eigenfunction. Hence, heuristically, all three methods of Section \ref{sec:coherent_sets} apply to the function $\f = P_S^* \widehat{\f}_S$. In all three methods, we defined a coherent set $A_\theta^\bullet$ as the positive support of $\text{Re}(\Per_\theta^t \f(\theta))$. We treat $\f$ as if it was an eigenfunction of $\G$ and assume that $\Per_\theta^t \f(\theta) = e^{zti} \f(\phi^t \theta)$. Hence, we define coherent sets as the positive level sets $\text{Re}(e^{\eta t i}\f(\phi^t\theta))$. In practice, other methods to extract almost-invariant
sets from an eigenfunction have also proven successful~\cite{froyland2005statistically}. In the numerical example of Section~\ref{sec:examples} we implement three different methods to extract coherent sets from an eigenfunction $\f$ of the generator~$\G$.
\begin{enumerate}[({CS}1)]
    \item $A^t_\theta := \{x\in M \:|\: \text{Re}(e^{\eta t i} \f(\phi^t\theta, x)) > 0\}.$
    \item Fix $q>0$. Define $A^t_\theta := \{x\in M \:|\: \norm{\text{Re}(e^{\eta t i} \f(\phi^t\theta))}_1^{-1}\abs{\text{Re}(e^{\eta t i} \f(\phi^t\theta, x))} > q\}.$
    \item Fix $q>0$. Define $A^t_\theta := \{x\in M \:|\: \norm{\f(\phi^t\theta)}_1^{-1} \abs{\f(\phi^t\theta, x)} > q\}.$
\end{enumerate}
Note that due to the normalization applied, the choice $q=1$ is natural in (CS2) and~(CS3).
The first method is the one described in \eqref{eq:positive_level}. The second method is almost identical, but instead of identifying the positive support of $\text{Re}(e^{\eta t i} \f(\phi^t\theta, x))$ as part of the coherent set, we choose the regions in which $e^{\eta t i} \f(\phi^t\theta, x)$ has a large absolute real part. Thus, the resulting sets $A_\theta^t$ consist of multiple smaller regions, rather than one large region. Additionally, unlike for the first method, the complement $B_\theta^t = M \setminus A_\theta^t$ is in general not coherent. The third method is similar to the second, but instead of the absolute real part, the absolute value in $\C$ is considered. Note that the factor $e^{\eta t i}$ is omitted, as it does not influence the absolute value. The third method is special in the sense that the resulting coherent sets satisfy $A_\theta^t = A_{\phi^t \theta}^0$. Hence, the resulting coherent set is of the form $A_\bullet = (A_\theta)_{\theta \in \Theta}$ with
\begin{equation*}
    A_\theta := \{x\in M \:|\: \norm{\f(\theta)}_1^{-1} \abs{\f(\theta, x)} > q\}.
\end{equation*}
The third method is to be used with caution. It may happen that $\abs{\f(\theta, x)}=1$ for all $\theta \in \Theta$ and $x\in M$ such that $A_\theta$ is either empty or the entire set $M$, thus rendering the method futile.
In Section~\ref{sec:examples}, we present two examples for which the third method has proven successful.

Given an eigenfunction $\widehat{\f}_S$ of $\Gamma_S$ and having chosen one of the three methods described above, we can determine whether a given point $x\in M$ lies inside or outside the set $A_\theta^t$ by evaluating $\f(\phi^t\theta)$ in the point $x$. For fixed starting parameter $\theta \in \Theta$, this allows us to estimate the survival probability of the coherent set $A_{\theta}^\bullet$ experimentally. At time $t=0$, we initialize a large number of particles that are uniformly distributed in~$M$. As time progresses, the parameter $\theta$ evolves according to $\Psi$, while the particles evolve according to the nonautonomous vector field~$v$. In each time frame $t>0$ of the simulation, we evaluate $\f(\phi^t\theta)$ in the positions of each of the particles and determine whether they lie inside or outside of $A_{\theta}^t$. We obtain an estimate of the survival probability of $A_{\theta}^\bullet$ up to time $t$ by computing the fraction of particles that have remained inside the coherent set up to time $t$ divided by the number of particles that started inside~$A_{\theta}^0$.

\subsection{Computational aspects}\label{sec:computational_asp}
In this section, we discuss aspects such as how to use sparsity of $\Gamma_S$ to shorten the computational time it takes to construct $\Gamma_S$ and to compute its eigenvalues with largest real part.

Recall the representation $\Gamma_S=D+R+A$ as the sum of two diagonal, and one skew-Hermitian matrix, cf.~\eqref{eq:gamma_DRA}--\eqref{eq:def_A}. The matrix $D$ is diagonal with real, nonpositive entries; $R$ is diagonal with purely imaginary entries; and $A$ is skew-Hermitian. Therefore, the off-diagonal entries of $\Gamma_S$ are skew-Hermitian, i.e.~they satisfy $\Gamma_S(m,n\:;\: m', n') = \overline{\Gamma_S(m',n'\:;\: m, n)}$ for all $(m,n)\neq (m',n')\in S$. Hence, it suffices to compute the upper-right half of $\Gamma_S$. From \eqref{eq:def_A}, it follows that an off-diagonal entry $\Gamma_S(m,n\:;\: m', n')$ is nonzero if and only if
\begin{equation*}
    n' \cdot \widehat{v}(m-m', n-n') = n \cdot \widehat{v}(m-m', n-n') \neq 0.
\end{equation*}
The fact that $n$ can be exchanged for $n'$ follows from Lemma \ref{lem:div_free_fourier}. If $v$ only has a small number of nonzero Fourier modes, then $\Gamma_S$ only has a small number of entries in each column, which leads to a high degree of sparsity. The vector field $v$ is assumed to be smooth, and smooth functions are well-approximated by a small number of Fourier modes. Hence, we set $\widehat{v}(m,n)=0$ whenever $\abs{\widehat{v}_i(m,n)}$ is below a given threshold $err>0$ for all $1\leq i \leq d_p$. In the examples of the next section, we used a threshold of $err=10^{-4}$. Using this threshold, no vector field considered in Section~\ref{sec:examples} had more than 150 nonzero Fourier modes. 

We compute the eigenvalues of $\Gamma_S$ using the MATLAB function \texttt{eigs}, which approximates the $k$ smallest magnitude eigenvalues of a given matrix. This function is based on matrix inverse iterations which requires a shift 
if $\Gamma_S$ is singular. Hence, we perform an eigenvalue shift by~$1$. We use \texttt{eigs} to compute the $k$ smallest magnitude eigenvalues of $\Gamma_S - \Id$ and afterwards add $1$ to all computed eigenvalues. This effectively computes the $k$ eigenvalues of $\Gamma_S$ that are closest to~$1\in \C$. The shift 1 is chosen if the spectral gap of $\Gamma$ is unknown. To obtain eigenvalues with more diverse real parts, the heuristic of taking the upper bound $|\varrho|$ on the spectral gap from Proposition~\ref{prop:spectral_gap} suggests itself. Note that $|\varrho| \propto \frac{\varepsilon^2}{2}$.

Using these methods, we are able to construct $\Gamma_S$ and compute $100$ eigenvalues around the origin on a conventional computer within minutes for sets $S$ with more than $30.000$ entries.

\section{Examples}\label{sec:examples}
We apply the numerical methods described in the previous section to three different types of nonautonomous dynamics. In all of the following examples, we consider $\Theta = \tor^2$ and $M=\tor^2$. Fix $\alpha \in \R^2$ that is rationally independent. Throughout, the dynamic on $\Theta$ is given by the quasi-periodic torus rotation by $\alpha$. 

The first two examples consist of a 2$\times$2 grid of counter-rotating gyres that is translated in space by an offset that depends on $\theta$. The third example consists of two shears of oscillating strength. The first example is special in the sense that it is dynamically equivalent to an autonomous dynamics. Therefore, we are able to derive the existence of exponential dichotomies and of eigenfunctions of the augmented generator~$\G$. For the second example, the existence of eigenfunctions of $\G$ is no longer given. However, our methods are still able to determine coherent sets that have slowly decaying survival probabilities. In the third example, there are no clear coherent structures apparent. However, using our methods we are able to determine coherent sets whose survival probabilities are close to the theoretical bound.

\subsection{Translated gyres}
\label{ssec:trans_gyres}
\begin{figure}
    \centering
    \begin{subfigure}{0.35\textwidth}
        \includegraphics[width=\textwidth]{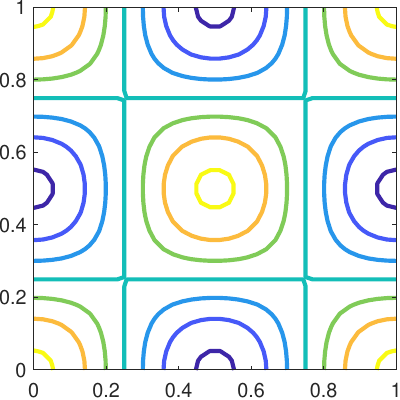}
    \end{subfigure}
    \hspace{0.1\textwidth}
    \begin{subfigure}{0.35\textwidth}
        \includegraphics[width=\textwidth]{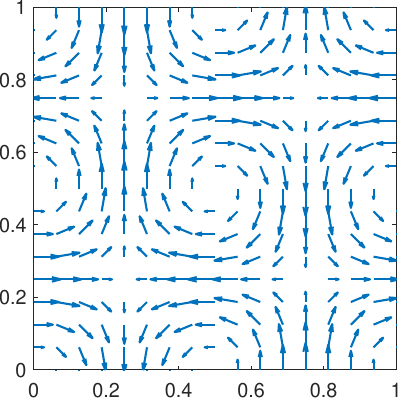}
    \end{subfigure}
    \caption{Left: Stream function  $\psi(x) = \frac{1}{2\pi}\cos(2\pi x_1) \cos(2\pi x_2)$. Right: The corresponding vector field $v_\text{aut}$.}
    \label{fig:Gyres}
\end{figure}

Consider the autonomous vector field on $M=\tor^2$,
\begin{equation}\label{eq:autonomous_gyres}
    v_\text{aut}(x) = \begin{pmatrix}
        \cos(2\pi x_1)\sin(2\pi x_2) \\
        -\sin(2\pi x_1) \cos(2\pi x_2)
    \end{pmatrix}.
\end{equation}
This vector field corresponds to the stream function $\psi(x) = \frac{1}{2\pi}\cos(2\pi x_1) \cos(2\pi x_2)$, and yields to a 2$\times$2 grid of rotating gyres, cf.~Figure \ref{fig:Gyres}. We construct a nonautonomous dynamics by translating $v_\text{aut}$ in space by $\theta \in \tor^2$. Define
\begin{equation}\label{eq:def_v_1}
    v(\theta, x) := v_\text{aut}(x + \theta).
\end{equation}
The dynamics on the parameter space $\Theta = \tor^2$ is given by a quasi-periodic torus rotation by~$\alpha \in \R^2$. 

\paragraph{Conjugacy to an autonomous flow.}
We claim that the vector field $v$ is dynamically equivalent to the vector field $v_\text{aut}$ with a constant drift in direction $\alpha$, i.e.
\begin{equation*}
    v_\text{eff}(x) = v_\text{aut}(x) + \alpha.
\end{equation*}
This dynamic equivalence can be understood as follows: As time progresses, the parameter $\theta\in \tor^2$ moves constantly in direction $\alpha$. Hence, the vector field $v(\theta, \cdot)$ is translated at a constant rate in the direction $-\alpha$. By changing the inertial frame of reference, this dynamics is equivalent to the vector field remaining constant and
its associated particle motion is superimposed with a translation of constant rate in the opposite direction~$\alpha$. 

This heuristic is rigorously described by Proposition~\ref{prop:commutative_diagram} below. Let $\Per_\theta^t$ be the solution operators to the Fokker--Planck equation with the nonautonomous vector field $v$, and let $\Per_\text{eff}^t$ be the solution operators to the Fokker--Planck equation with the autonomous vector field $v_\text{eff}$. For $\beta \in \tor^2$, let $T_\beta:\H \to \H$ be the translation by $\beta$, i.e.~
\begin{equation*}
    [T_\beta f](x) = f(x - \beta).
\end{equation*}
Note that $T_\beta$ is an isometric isomorphism with inverse $T_{-\beta}$.
\begin{proposition}\label{prop:commutative_diagram}
    The following diagram commutes for every $\theta \in \tor^2$ and $t\geq 0$.
    \begin{equation}\label{eq:commutative_diagram}
    \begin{tikzcd} \H \arrow[r, "\Per_\theta^t"]\arrow[d, "T_\theta"] &\H \arrow[d,"T_{\phi^t \theta}"]\\ \H \arrow[r,"\Per_\text{eff}^t"] & \H 
    \end{tikzcd}
    \end{equation}
\end{proposition}
\begin{proof}
    See Appendix \ref{ap:proofs}.
\end{proof}

\paragraph{Mather and Sacker--Sell spectra.}
The fact that the cocycle $\Per$ is cohomologous to the autonomous cocycle $\Per_\text{eff}$ implies that the Mather operators $\mathbf{M}^t$ admit eigenfunctions and that the Sacker--Sell spectrum $\Sigma(\Per)$ consists of discrete points. 

Recall Section \ref{sec:autonomous_case}, in which we studied autonomous cocycles, i.e.~one-parameter semigroups. The Sacker--Sell spectrum of $\Per_\text{eff}$ consists of discrete values $\Sigma(\Per_\text{eff})=\{\lambda_k \mid 1\leq k \leq N \}$, for $N\in \N_0 \cup \{\infty\}$. For each $\lambda_k$ there is $z_k\in \C$ with $\abs{z_k} = \lambda_k$ and a function $f\in \H_0$ such that $\Per_\text{eff}^t f = e^{z_kt}f$ for all $t\geq 0$. Let $f$ be such an eigenfunction. Define a function $\f\in L^2(\tor^2, \H)$ by
\begin{equation*}
    \f(\theta) = T_{-\theta} f.
\end{equation*}
The commutative diagram \eqref{eq:commutative_diagram} yields 
\begin{align*}
   \lbrack\mathbf{M}^t \f\rbrack (\theta) &= \Per_{\phi^{-t}\theta}^t \f(\phi^{-t}\theta) = (T_{-\theta} \circ \Per_\text{eff}^t \circ T_{\phi^{-t} \theta}) T_{-\phi^{-t}\theta} f \\
   &= e^{z_kt} T_{-\theta} f = e^{z_kt} \f(\theta).
\end{align*}
This shows that $\f$ is an eigenfunction of the Mather semigroup with eigenvalue $z_k$. 

\begin{figure}[htb]
    \centering
    \begin{subfigure}{0.45\textwidth}
         \includegraphics[width=\textwidth]{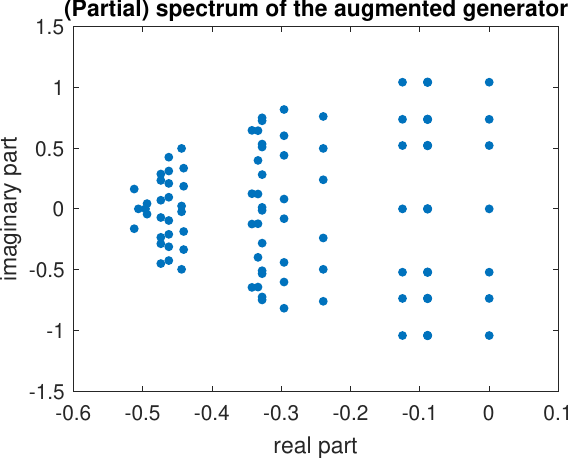}
    \end{subfigure}
    \hfill
    \begin{subfigure}{0.45\textwidth}
        \includegraphics[width=\textwidth]{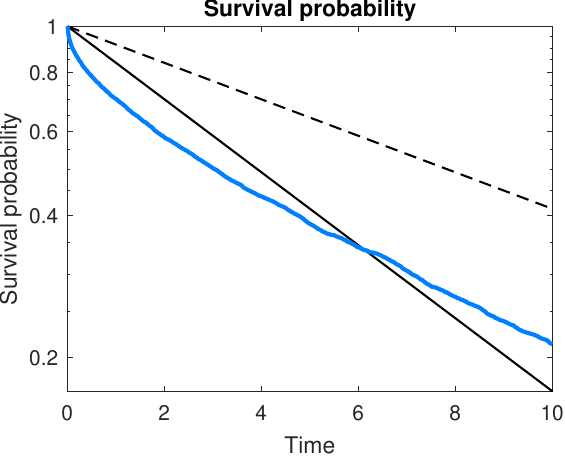}
    \end{subfigure}
   
    \caption{Translated gyre example. Left: 100 eigenvalues of the discrete generator $\Gamma_S$ around the origin. Right: In blue, the simulated survival probability of the coherent sets extracted from an eigenfunction $\f$ to an eigenvalue $z$ using method (CS3) with threshold $q=1$. The black lines show the curves $e^{\lambda t}$ (dashed) and $e^{2\lambda t}$ (solid) for comparison, where $\lambda$ is the real-part of $z$. The computed cumulative survival probability is $C(A_{\theta}^\bullet)= 6.070$.}
    \label{fig:Translated_Gyres_spec_sp}
\end{figure}

\begin{proposition}\label{prop:semi_aut_SS}
    We find $\Sigma(\Per) = \Sigma(\Per_\text{eff}) = \{\lambda_k \mid 1\leq k \leq N\}$.
\end{proposition}
\begin{proof}
    See Appendix \ref{ap:proofs}.
\end{proof}
By Theorem \ref{thm:spectral_mapping}, the spectrum of the augmented generator $\G$ is given by lines parallel to the imaginary axis with real parts~$\lambda_k$.
The spectrum of the discretized generator $\Gamma_S$, which we will construct next, is compatible with this structure
seen in Figure~\ref{fig:Translated_Gyres_spec_sp}.

\paragraph{Discrete generator.}
We construct a numerical approximation $\Gamma_S\in \C^{S\times S}$ of the augmented generator $\G$, as described in Section~\ref{sec:numerics}. This requires us to select a set of Fourier modes $S\subset \Z^4$. We use the sparsity structure of $\Gamma$ to choose a suited set $S$. Recall that an off-diagonal entry of $\Gamma$ is only nonzero if $n \cdot \widehat{v}(m-m', n-n')\neq 0$, cf.~\eqref{eq:def_Gamma}. The vector field $v$, defined by \eqref{eq:def_v_1}, only has four nonzero Fourier modes $\widehat{v}(m,n)$. Since $m,n\in \Z^2$, we write them as $\widehat{v}(m_1, m_2, n_1, n_2)$. The four nonzero Fourier modes of $v$ are
\begin{align*}
    \widehat{v}(+1,+1,+1,+1) &= \frac{1}{4} i \begin{pmatrix}-1 \\ +1 \end{pmatrix}, \qquad &
    \widehat{v}(+1,-1,+1,-1) &= \frac{1}{4} i \begin{pmatrix}+1 \\ +1 \end{pmatrix}, \\
    \widehat{v}(-1,+1,-1,+1) &= \frac{1}{4} i \begin{pmatrix}-1 \\ -1 \end{pmatrix}, \qquad &
    \widehat{v}(-1,-1,-1,-1) &= \frac{1}{4} i \begin{pmatrix}+1 \\ -1 \end{pmatrix}.
\end{align*}
The fact that these are the only nonzero Fourier modes implies that $\Gamma$ has a block diagonal structure after suitably reordering the indices, i.e.~there are sets $\mathcal{I}_k\subset \Z^4$ such that $\Gamma(\mathcal{I}_k; \mathcal{I}_\ell) = \Gamma(\mathcal{I}_\ell; \mathcal{I}_k)= 0$ for any $k\neq \ell$. More precisely, an off-diagonal entry $\Gamma(m,n;m',n')$ that is nonzero requires $m-m' = n-n' = (\pm1,\pm 1)$, and in particular, $m-n=m'-n'$.
The requirement $m_1-n_1 = m_1'-n_1'$ and $m_2 - n_2 = m_2'-n_2'$ suggests to define for $(k,\ell) \in \Z^2$ the class
\begin{equation}\label{eq:equivalence_class}
    \mathcal I_{k,\ell} := \{(n_1 + k, n_2 + \ell, n_1, n_2) \mid (n_1,n_2) \in \Z^2\} \subset \Z^4.
\end{equation}
Now for two different classes $\mathcal{I}, \mathcal{J} \subset\Z^4$ we obtain $\Gamma(\mathcal{I}; \mathcal{J}) = \Gamma(\mathcal{J}; \mathcal{I})= 0$.
We note that $\Gamma_S$ can be decomposed along % into
an even finer set of equivalence classes. 

Note that the blocks of $\Gamma$ are of infinite size. However, after choosing $S\subset\Z^4$, the resulting finite-dimensional matrix $\Gamma_S \in \C^{S\times S}$ inherits a block diagonal structure from~$\Gamma$. 
Hence, the spectrum of $\Gamma_S$ and the corresponding eigenfunctions are given by the spectrum and the eigenfunctions of each of the individual blocks. It is desirable to choose $S\subset \Z^4$ such that $\Gamma_S$ contains few diagonal block, each of high resolution, rather than many diagonal blocks of low resolution. We define $S$ in the following way: Fix $K\in \N$ and~$r>0$. Define
\begin{equation}\label{eq:translated_gyres_S}
    S:= \{(n_1 + k, n_2 + \ell, n_1, n_2) \in \Z^4 \mid k,\ell \in \{-K,\hdots, K\},\: \norm{(n_1, n  _2)}_e \leq r  \}.
\end{equation}
The set $S$ intersects $(2K+1)^2$ equivalence classes of the type~\eqref{eq:equivalence_class}. The resolution of each of the classes is determined by the radius~$r$. 

\begin{figure}
    \centering
    \begin{subfigure}{0.3\textwidth}
        \includegraphics[width=\textwidth]{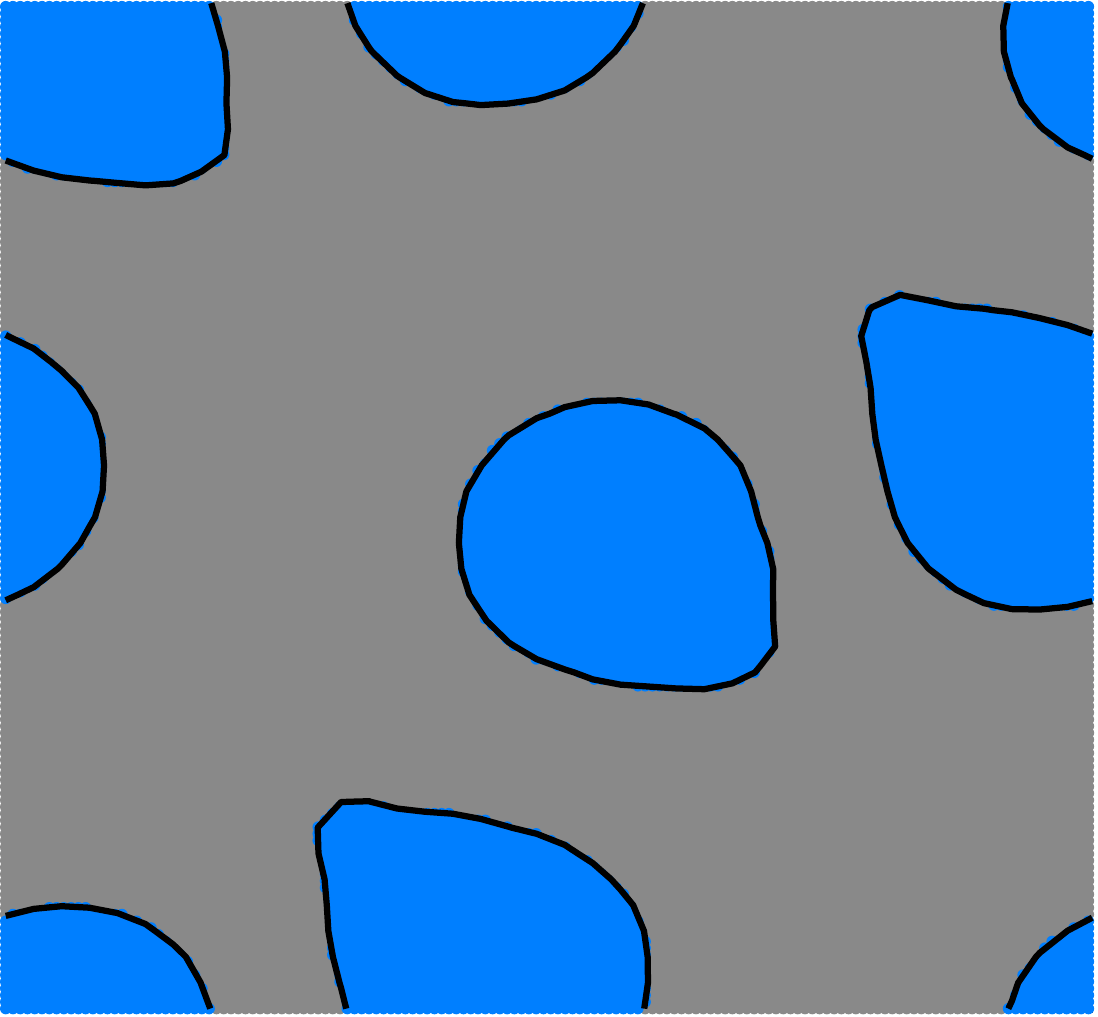}
    \end{subfigure}
    \hfill
    \begin{subfigure}{0.3\textwidth}
        \includegraphics[width=\textwidth]{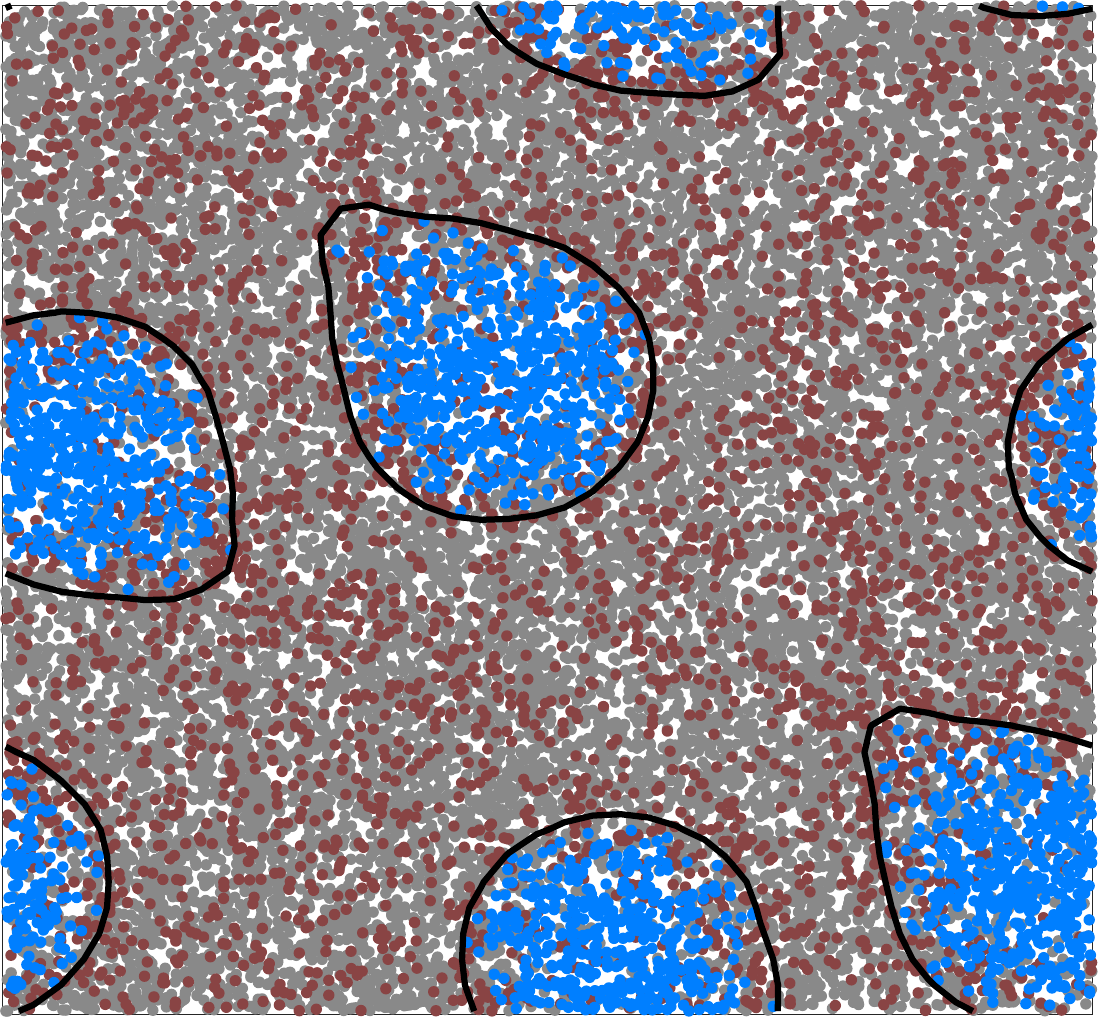}
    \end{subfigure}
    \hfill
    \begin{subfigure}{0.3\textwidth}
        \includegraphics[width=\textwidth]{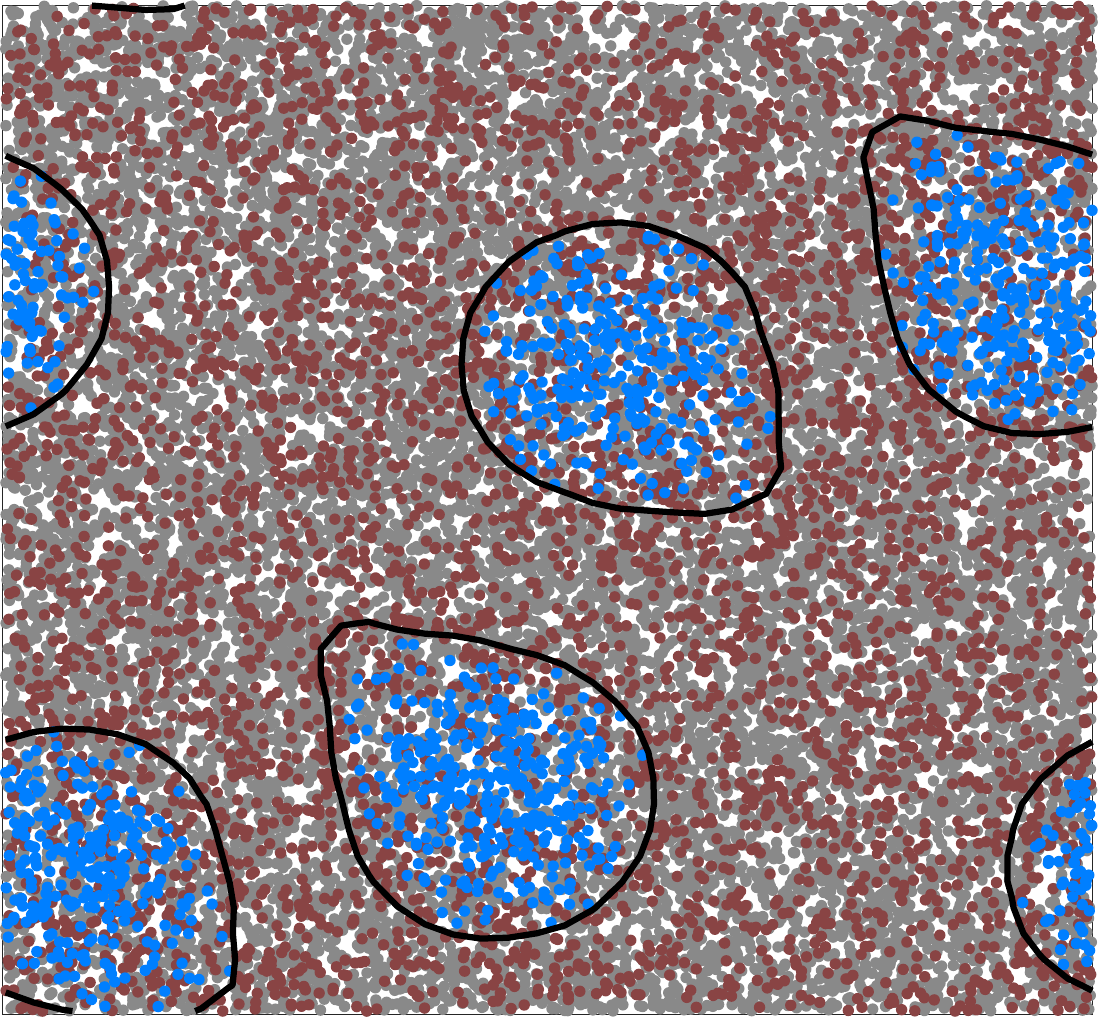}
    \end{subfigure}\\
    \vspace{3mm}
    \begin{subfigure}{0.3\textwidth}
        \includegraphics[width=\textwidth]{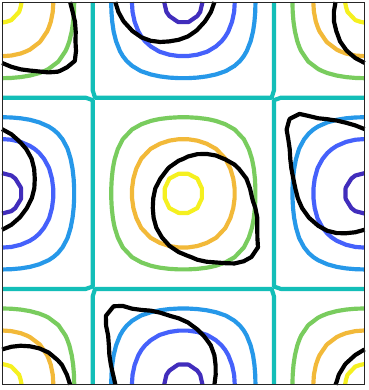}
    \end{subfigure}
    \hfill
    \begin{subfigure}{0.3\textwidth}
        \includegraphics[width=\textwidth]{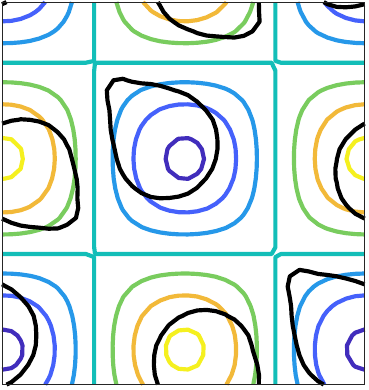}
    \end{subfigure}
    \hfill
    \begin{subfigure}{0.3\textwidth}
        \includegraphics[width=\textwidth]{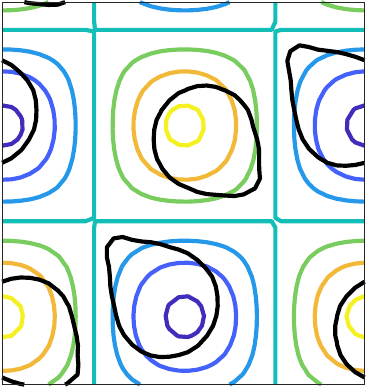}
    \end{subfigure}
    \caption{Translated gyre example: Evolution of the coherent sets, marked in black, determined using method (CS3) with threshold $q = 1$. The snapshots are taken at time $t=0,5,10$. The top row depicts a particle simulation. Particles that started inside $A_{\theta}^0$ and have remained there until time $t$ are colored in blue. Particles that have started inside $A_{\theta}^0$ but were outside the coherent set at some time $0<s\leq t$ are colored in red. Other particles are colored in grey. The bottom row locates the coherent set with respect to the position of the gyres, as indicated by the stream function at the respective time. A video of the particle simulation for the time interval $[0,10]$ can be found at \url{https://github.com/RobinChemnitz/MatherCoherent}.}
    \label{fig:Translated_Gyres_snapshots}
\end{figure}

\paragraph{Results.}
In the following we present experimental results for the described system. We choose $\alpha = 0.2 (1, \sqrt{2})$ as the direction of the quasiperiodic driving and $\varepsilon = 0.03$ as the diffusion constant. The set $S$ is constructed like described in \eqref{eq:translated_gyres_S} with constants $K=2$ and $r = 11$. This results in a set $S\subset \Z^4$ with $9425$ elements. We compute the discrete generator~$\Gamma_S$, which is a sparse matrix 
with 0.04$\%$ nonzero entries. In Figure~\ref{fig:Translated_Gyres_spec_sp} we show the spectrum of $\Gamma_S$ around the origin, computed as described in Section \ref{sec:computational_asp}. Note that the spectrum of $\Gamma_S$ is distributed on lines parallel to the imaginary axis. This fits to the description of the spectrum of $\G$, which we showed to consist of discrete lines parallel to the imaginary axis. We pick an eigenvalue $z = -0.089 - 1.041i$ of $\Gamma_S$ of largest nonzero real part. Hence, the real part is given by~$\lambda = -0.089$. Let $\widehat{\f}_S$ be the corresponding eigenfunction. As described in Section~\ref{sec:numerics_coherent_sets}, $\widehat{\f}_S$ represents a function $\f:= P_S^* \widehat{\f}_S\in L^2(\tor^2\times \tor^2, \C)$. For a given $\theta \in \tor^2$, we compute a coherent set $A_\theta^\bullet$ based on the real part of the fibres of~$\f$, cf.\ method~(CS2). Fix the threshold $q=1$ and define
\begin{equation*}
    A^t_\theta := \{x\in \tor^2 \:|\: \norm{\text{Re}(e^{\eta t i} \f(\phi^t\theta))}_1^{-1}\abs{\text{Re}(e^{\eta t i} \f(\phi^t\theta, x))} > q\}.
\end{equation*}
We estimate the survival probability for the starting parameter $\theta = (0, 0)$. At time $t=0$, we fill the domain $\tor^2$ with a $150\times150$ grid of particles and let them evolve in time according to the SDE \eqref{eq:SDE} with the nonautonomous vectorfield $v$. For each time~$t>0$, we obtain an experimental estimate on the survival probability of $A_{\theta}^\bullet$ up to time~$t$, as described in Section~\ref{sec:numerics_coherent_sets}. In Figure \ref{fig:Translated_Gyres_snapshots}, we show snapshots of the simulation for three selected times. In Figure \ref{fig:Translated_Gyres_spec_sp}, there is a plot of the experimental survival probability up to time~$t=10$. From \eqref{eq:eigenfunction_survival_prob}, we expect that the survival probability of $A_{\theta}^\bullet$ decays at most at rate $2\lambda$, up to constants. The results in Figure~\ref{fig:Translated_Gyres_spec_sp} agree with this estimate. For $t$ close to~$0$, the survival probability decays quickly but for later times~$t$, the decay of the survival probability is significantly slower than~$e^{2\lambda t}$. The computed cumulative survival probability is $C(A_{\theta}^\bullet)= 6.070$. This is more than the estimate $-\frac{1}{2 \lambda}= 5.682$ from \eqref{eq:cum_prob_estimate}. Determining coherent sets based on the absolute value of the fibres $\f(\theta)$ rather than their absolute real part, i.e.~based on method (CS3), yields almost identical results. Applying the method (CS1) also yields a coherent set with slowly decaying survival probability (not depicted here). 

\subsection{Oscillating gyres}
\label{ssec:osc_gyre}

\begin{figure}
    \centering
    \includegraphics[width=0.4\textwidth]{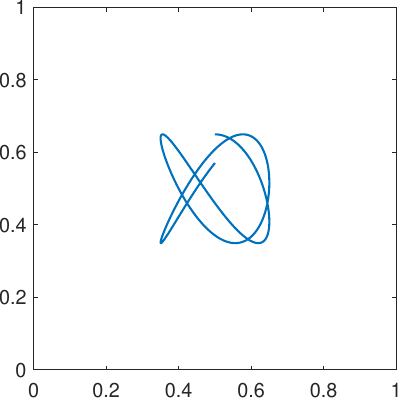}
    \caption{Oscillating gyre example: The trajectory up to time $t=10$ of the midpoint of the central gyre of the vector field $v$, defined in \eqref{eq:oscillating_gyres} for parameters $\alpha=0.2(1, \sqrt{2})$, $\delta = 0.15$, and strating parameter $\theta=(0,0)$.}
    \label{fig:oscillation}
\end{figure}

This example is defined similarly to the previous one. However, unlike in the previous case, the nonautonomous vector field we construct is not dynamically equivalent to an autonomous one. Like before, consider the autonomous vector field $v_\text{aut}$, defined in \eqref{eq:autonomous_gyres}, that consists of a 2$\times$2 grid of counter-rotating gyres, cf.\ Figure~\ref{fig:Gyres}.
We construct a nonautonomous vector field by letting the grid of gyres oscillate in a quasi-periodic manner. Let $\delta> 0$ be the amplitude of the oscillation. Define
\begin{equation}\label{eq:oscillating_gyres}
    v(\theta, x) = v_\text{aut}\left(x + \delta \begin{pmatrix}
        \sin(2\pi \theta_1) \\ \cos(2\pi \theta_2)
    \end{pmatrix}\right).
\end{equation}
Since $\theta \in \tor^2$ evolves constantly in direction $\alpha \in \R^2$, the grid of gyres oscillates in direction of $x_1$ with a period of $\alpha_1^{-1}$ and in direction $x_2$ with a period of~$\alpha_2^{-1}$. A trajectory of the oscillation is depicted in Figure~\ref{fig:oscillation}.

\begin{figure}
    \centering
    \begin{subfigure}{0.45\textwidth}
         \includegraphics[width=\textwidth]{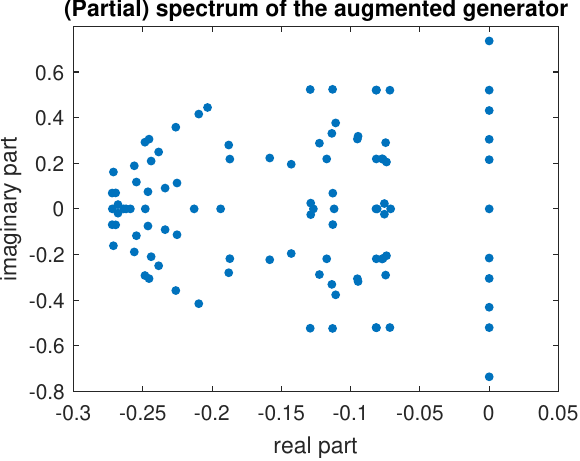}
    \end{subfigure}
    \hfill
    \begin{subfigure}{0.45\textwidth}
        \includegraphics[width=\textwidth]{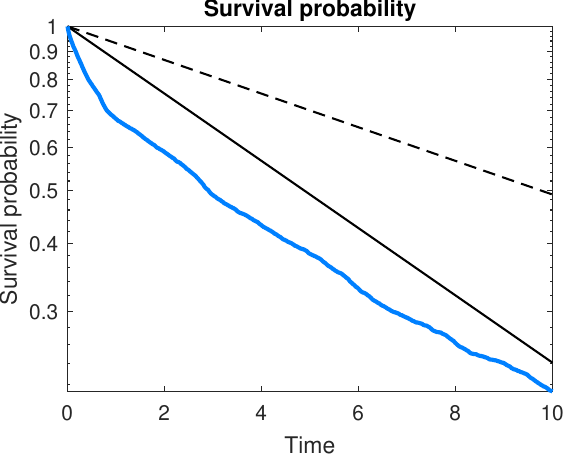}
    \end{subfigure}
   
    \caption{Oscillating gyre example. Left: 100 eigenvalues of the discrete generator $\Gamma_S$ around the origin. Right: In blue, the simulated survival probability of the coherent sets determined by the absolute value of an eigenfunction $\f$ to an eigenvalue $z$ with threshold $q=1$. The black lines show the curves $e^{\lambda t}$ (dashed) and $e^{2\lambda t}$ (solid) for comparison, where $\lambda$ is the real-part of $z$. The computed cumulative survival probability is $C(A_\theta^\bullet) =  6.158$.}
    \label{fig:Oscillating_Gyres_spec_sp}
\end{figure}

\paragraph{Mather and Sacker--Sell spectra.}
For $\delta$ close to $0$, the vector field is almost autonomous. Unlike in the previous example, this does not guarantee the existence of eigenfunctions of the Mather operator $\mathbf{M}^t$, but for $\delta$ small enough, we can prove the existence of exponential dichotomies. Let $\Per_\text{aut}$ be the autonomous cocycle of transfer operators corresponding to the autonomous vector field~$v_\text{aut}$. From Section~\ref{sec:autonomous_case}, we know that the Sacker--Sell spectrum $\Sigma(\Per_\text{aut})=\{\lambda_k \mid 1\leq k \leq N\}$, for $N\in \N_0 \cup \{\infty\}$, consists of discrete points. Let~$\lambda \notin \Sigma(\Per_\text{aut})$. Hence, the spectrum of the Mather operators $\mathbf{M}^t_\text{aut}$ of $\Per_\text{aut}$ satisfies~$\sigma(\mathbf{M}^t_\text{aut}) \cap e^{\lambda t} S^1 = \emptyset$. In other words, the ring $e^{\lambda t} S^1$ lies in the resolvent set of $\mathbf{M}^t_\text{aut}$ for each~$t>0$. By \cite[Corollary 6.44]{chicone1999evolution}, the nonautonomous cocycle $\Per$ has an exponential dichotomy at $\lambda$ if for some~$t>0$
\begin{equation}\label{eq:perturbation_bound}
    \sup_{\theta \in \Theta} \norm{\Per_\text{aut}^t - \Per_\theta^t} < \inf_{z\in e^{\lambda t} S^1} \norm{(\mathbf{M}^t_\text{aut} - z\Id )^{-1}}^{-1}.
\end{equation}
For fixed $t>0$, the right-hand side defines a constant that only depends on $t$ and~$\lambda$. As $\delta$ approaches $0$, the difference between the vector fields $v$ and $v_\text{aut}$ tends to~$0$ uniformly on~$\tor^2$. By Lemma~\ref{lem:uniform_regularity}, the left-hand side of \eqref{eq:perturbation_bound} vanishes as~$\delta \to 0$. We conclude that for sufficiently small~$\delta$, the nonautonomous cocycle $\Per$ has an exponential dichotomy at $\lambda$ and the spectrum of the augmented generator $\G$ does not contain the line $\{\lambda + \eta i \mid \eta \in \R\}$, cf.\ Theorem~\ref{thm:spectral_mapping}. Hence, for small values of~$\delta$, we expect the spectrum of $\G$ to consist of thin bands along the imaginary axis, at least for regions close to the origin. In Figure~\ref{fig:Oscillating_Gyres_spec_sp}, we show the spectrum of the discrete generator $\Gamma_S$ for $\delta=0.15$. This value of $\delta$ seems to be not small enough such that the spectrum consists of clear lines, nevertheless, the band structure of the spectrum is recognizable in the region close to~$0$.

\begin{figure}
    \centering
    \begin{subfigure}{0.3\textwidth}
        \includegraphics[width=\textwidth]{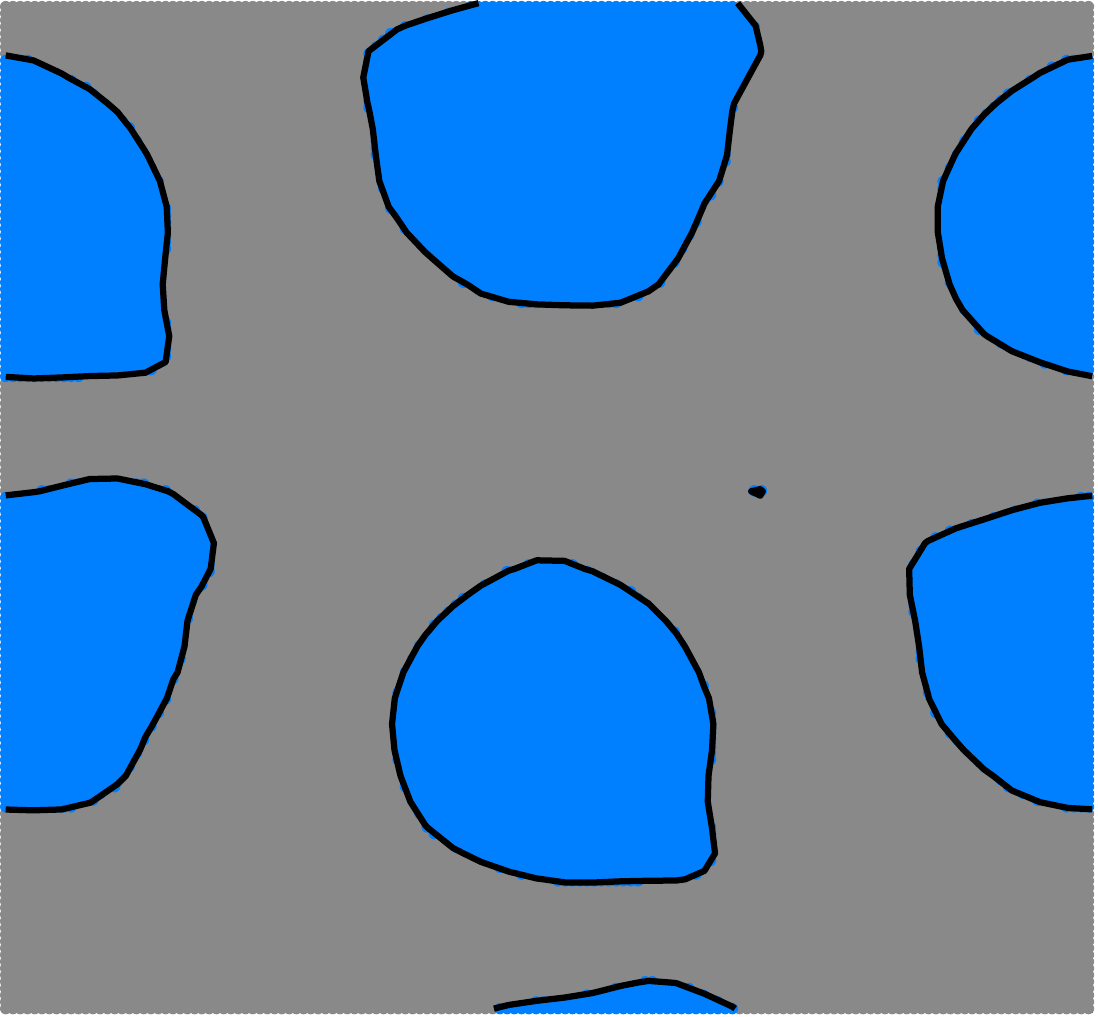}
    \end{subfigure}
    \hfill
    \begin{subfigure}{0.3\textwidth}
        \includegraphics[width=\textwidth]{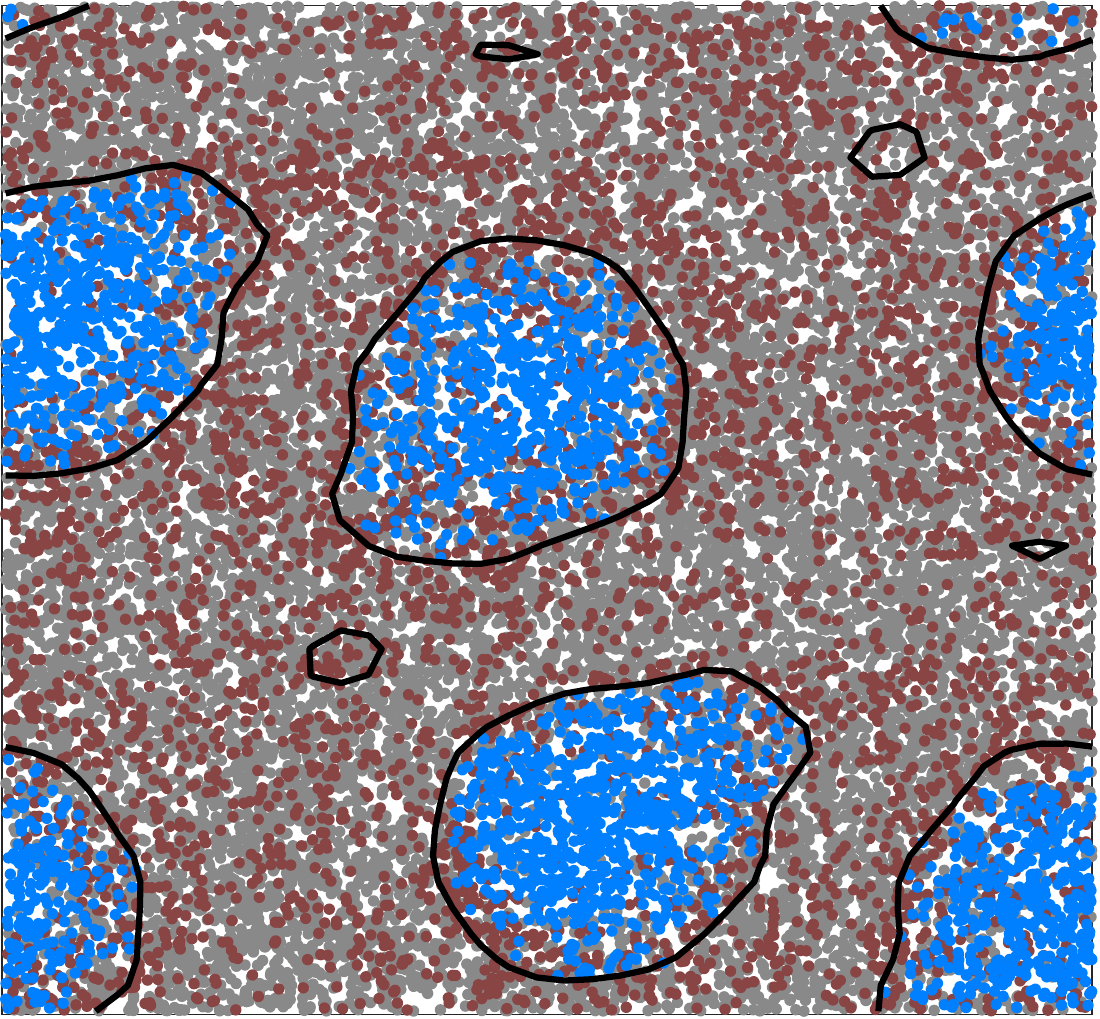}
    \end{subfigure}
    \hfill
    \begin{subfigure}{0.3\textwidth}
        \includegraphics[width=\textwidth]{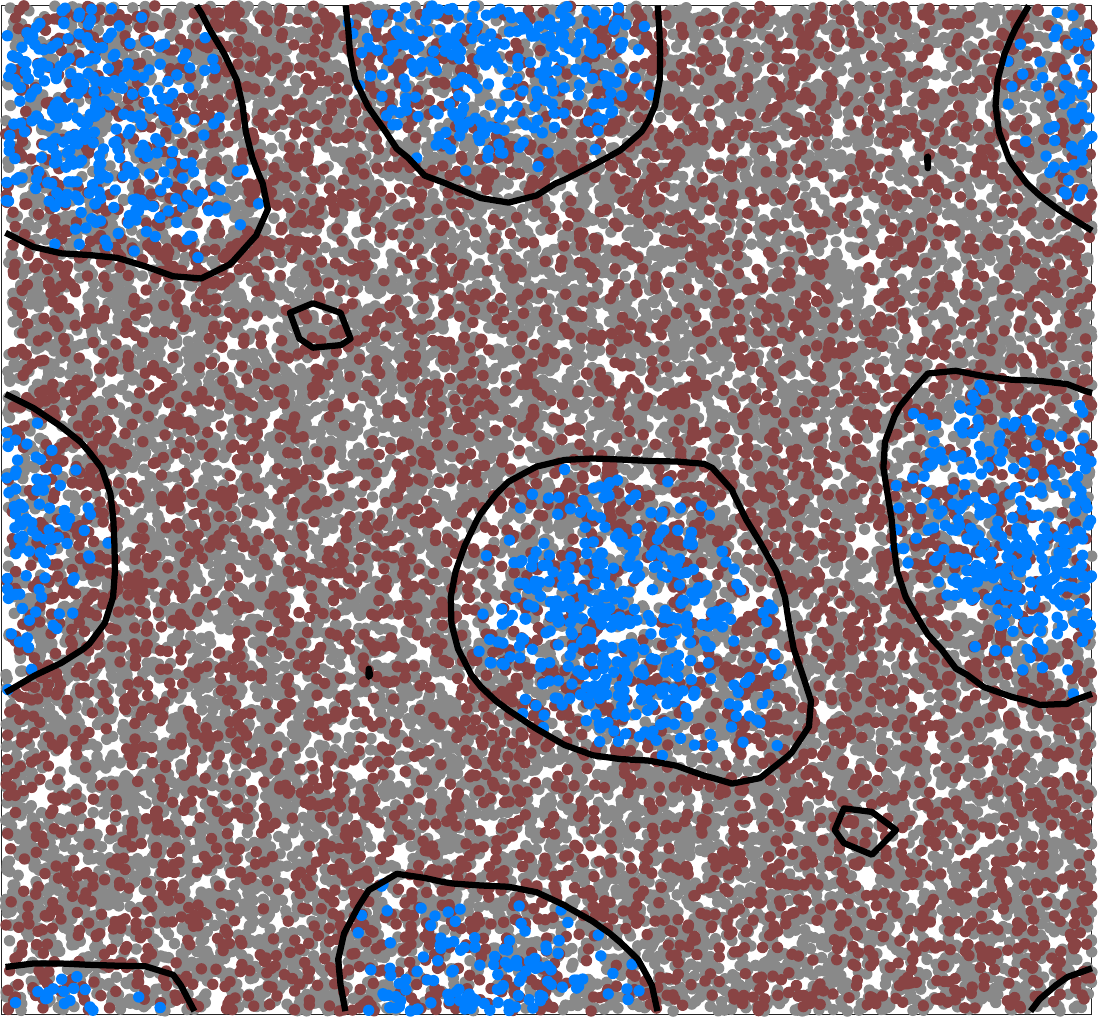}
    \end{subfigure} \\
    \vspace{3mm}
    \begin{subfigure}{0.3\textwidth}
        \includegraphics[width=\textwidth]{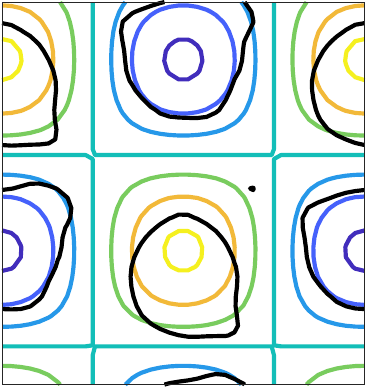}
    \end{subfigure}
    \hfill
    \begin{subfigure}{0.3\textwidth}
        \includegraphics[width=\textwidth]{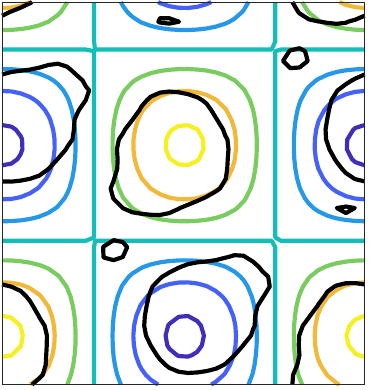}
    \end{subfigure}
    \hfill
    \begin{subfigure}{0.3\textwidth}
        \includegraphics[width=\textwidth]{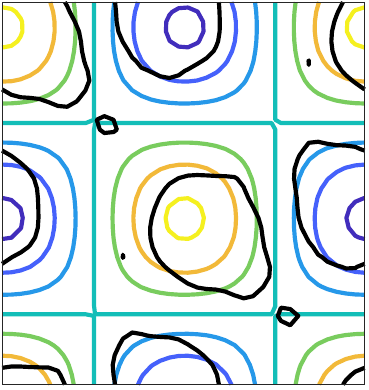}
    \end{subfigure}
    \caption{Oscillating gyre example: Evolution of the coherent sets, marked in black, determined using method (CS2) with threshold $q = 1$. The snapshots are taken at time $t=0,5,10$. The top row depicts a particle simulation. Particles that started inside $A_{\theta}^0$ and have remained there until time $t$ are colored in blue. Particles that have started inside $A_{\theta}^0$ but were outside the coherent set at some time $0<s\leq t$ are colored in red. Other particles are colored in grey. The bottom row locates the coherent set with respect to the position of the gyres, as indicated by the stream function at the respective time.  A video of the particle simulation for the time interval $[0,10]$ can be found at \url{https://github.com/RobinChemnitz/MatherCoherent}.}
    \label{fig:Oscillating_Gyres_snapshots}
\end{figure}

\paragraph{Results.}
In the following the present experimental results for the described system. We choose $\alpha = 0.2 (1, \sqrt{2})$ as the direction of the quasiperiodic driving, $\varepsilon = 0.03$ as the diffusion constant, and $\delta = 0.15$ as the strength of the oscillation. We define the set $S\in \Z^4$ of Fourier modes as described in~\eqref{eq:S_d_S_p}.
\begin{equation*}
S := \{(m_1, m_2, n_1, n_2)\in \Z^4 \mid \norm{(m_1, m_2)}_\infty \leq K, \: \norm{(n_1, n_2)}_e \leq r\},
\end{equation*}
where $K=6$ and $r=8$. The constant $K$ determines the resolution in the driving and $r$ determines the resolution in physical space. The chosen values of $K$ and $r$ result in a set $S\subset \Z^4$ with $33293$ elements. We compute the discrete generator~$\Gamma_S$, which is a sparse matrix 
with $0.31\%$ nonzero entries. In Figure \ref{fig:Oscillating_Gyres_spec_sp} we show the spectrum of $\Gamma_S$ around the origin, computed as described in Section~\ref{sec:computational_asp}. We pick a real eigenvalue $z = -0.071$ of $\Gamma_S$ of largest nonzero real part. Since $z$ is real, we have~$\lambda = z$. Let $\widehat{\f}_S$ be the corresponding eigenfunction. As described in Section~\ref{sec:numerics_coherent_sets}, $\widehat{\f}_S$ represents a function $\f:= P_S^* \widehat{\f}_S\in L^2(\tor^2\times \tor^2, \C)$. For a given $\theta \in \tor^2$, we compute a coherent set $A_\theta^\bullet$ based on the absolute value of the fibres of~$\f$, cf.~method (CS2). Fix the threshold $q=1$ and define
\begin{equation*}
    A^t_\theta := \{x\in \tor^2 \:|\: \norm{\f(\phi^t\theta)}_1^{-1} \abs{\f(\phi^t\theta, x)} > q\}.
\end{equation*}
Recall that these sets satisfy $A_{\theta}^t = A_{\phi^t \theta}^0$. Hence, the coherent family can be written as $A_\bullet$. We estimate the survival probability for the starting parameter $\theta = (0, 0)$. At time $t=0$, we fill the domain $\tor^2$ with a $150\times150$ grid of particles and let them evolve in time according to the SDE \eqref{eq:SDE} with the nonautonomous vectorfield~$v$. For each time $t>0$, we obtain an experimental estimate on the survival probability of $A_{\theta}^\bullet$ up to time $t$, as described in Section~\ref{sec:numerics_coherent_sets}. In Figure~\ref{fig:Oscillating_Gyres_snapshots}, we show snapshots of the simulation for three selected time points. In Figure~\ref{fig:Oscillating_Gyres_spec_sp}, there is a plot of the experimental survival probability up to time~$t=10$. From~\eqref{eq:eigenfunction_survival_prob}, we expect that the survival probability of $A_{\theta}^\bullet$ decays at most at rate $2\lambda$, up to constants. The results in Figure~\ref{fig:Oscillating_Gyres_spec_sp} agree with this estimate as the estimated survival probability decays very similarly to~$e^{2\lambda t}$. The computed cumulative survival probability is $C(A_\theta^\bullet) =  6.158$. This is less than the estimate $-\frac{1}{2 \lambda} = 7.042$ from~\eqref{eq:cum_prob_estimate}. Determining coherent sets based on the absolute real part of the fibres $\f(\theta)$ rather than their absolute value, i.e.~based on method (CS2), yields almost identical results. Applying the method (CS1) also yields a coherent set with slowly decaying survival probability (not depicted here). 

\subsection{Shear}\label{sec:shear}
In both of the previous examples, the coherent sets we identified were regions around the center of the four gyres, close to sets that could have been determined from considering the individual vector fields $v(\theta, \cdot)$ for each fixed~$\theta \in \tor^2$.
The example we present next is of different 
nature.

\begin{figure}
    \centering
    \begin{subfigure}{0.3\textwidth}
        \includegraphics[width=\textwidth]{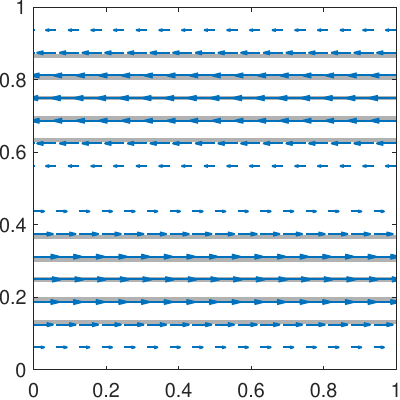}
    \end{subfigure}
    \hfill
        \begin{subfigure}{0.3\textwidth}
        \includegraphics[width=\textwidth]{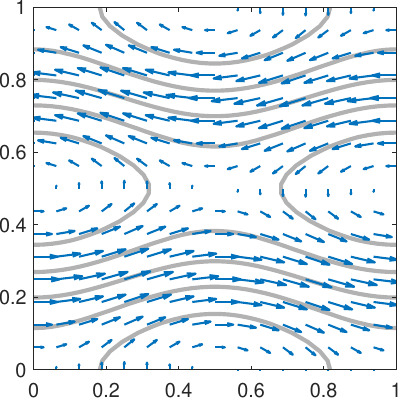}
    \end{subfigure}
    \hfill
        \begin{subfigure}{0.3\textwidth}
        \includegraphics[width=\textwidth]{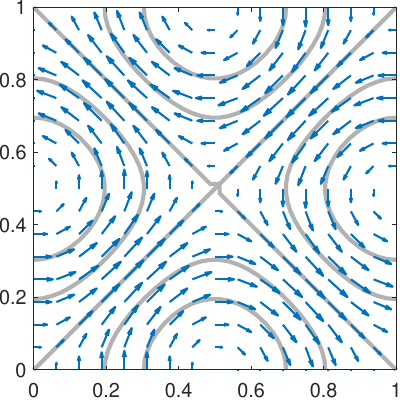}
    \end{subfigure}
    \caption{Shear example: The vector field defined in \eqref{eq:shear_vec} in blue and the corresponding stream function in grey for parameters $\theta=(0.25, 0)$, $(0.25, 0.05)$, $(0.25, 0.25)$. From left to right, the strength of the shear in $x$-direction stays constant at 1 while the shear in $y$-direction increases from 0 to~1.}
    \label{fig:shear_quiver}
\end{figure}

A well-studied dynamical feature that is known to result in chaotic behavior is shear. In this example we overlay two shears of oscillating strength, by considering
the nonautonomous vector field
\begin{equation}\label{eq:shear_vec}
    v(\theta, x) = \sin(2\pi \theta_1) \begin{pmatrix}
         \sin(2\pi x_2) \\ 0 \end{pmatrix} + \sin(2\pi \theta_2) \begin{pmatrix} 0 \\ \sin(2\pi x_1)
    \end{pmatrix}.
\end{equation}
Figure~\ref{fig:shear_quiver} shows a representation of this vector field for selected values of~$\theta$. Determining coherent sets of $v$ is non-trivial, since vector fields with shears are known for their mixing properties~\cite{ChSo89,blumenthal2023exponential, FrKo23}, while coherent sets are characterized by experiencing little mixing with their surrounding. Despite this feature, which opposes the existence of coherent sets, our methods are able to identify a coherent set whose survival probability decays at a similar rate as the survival probability of the previous example of the oscillating gyres. Note that this comparison is fair in the sense that the velocity field strength is of the same order of magnitude and below we take the same diffusion coefficient $\varepsilon$ as in the previous examples.

\begin{figure}
    \centering
    \begin{subfigure}{0.45\textwidth}
         \includegraphics[width=\textwidth]{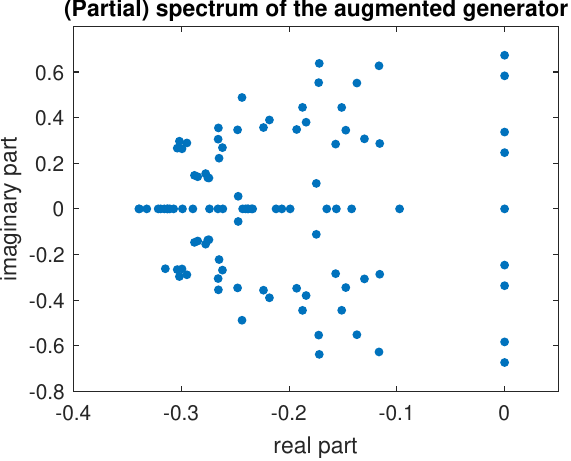}
    \end{subfigure}
    \hfill
    \begin{subfigure}{0.45\textwidth}
        \includegraphics[width=\textwidth]{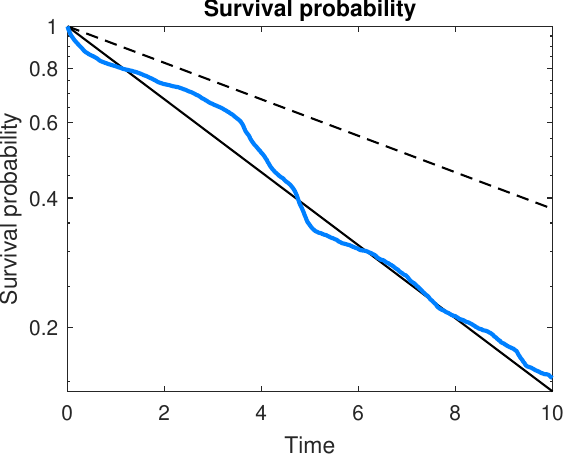}
    \end{subfigure}
   
    \caption{Shear example. Left: 100 eigenvalues of the discrete generator $\Gamma_S$ around the origin. Right: In blue, the simulated survival probability of the coherent sets determined by positive support of an eigenfunction $\f$ to an eigenvalue $z$. The black lines show the curves $e^{\lambda t}$ (dashed) and $e^{2\lambda t}$ (solid) for comparison, where $\lambda$ is the real-part of $z$. The computed cumulative survival probability is $C(A_{\theta}^\bullet)= 5.680$.}
    \label{fig:Shear_spec_sp}
\end{figure}

\begin{figure}
    \centering
    \begin{subfigure}{0.3\textwidth}
        \includegraphics[width=\textwidth]{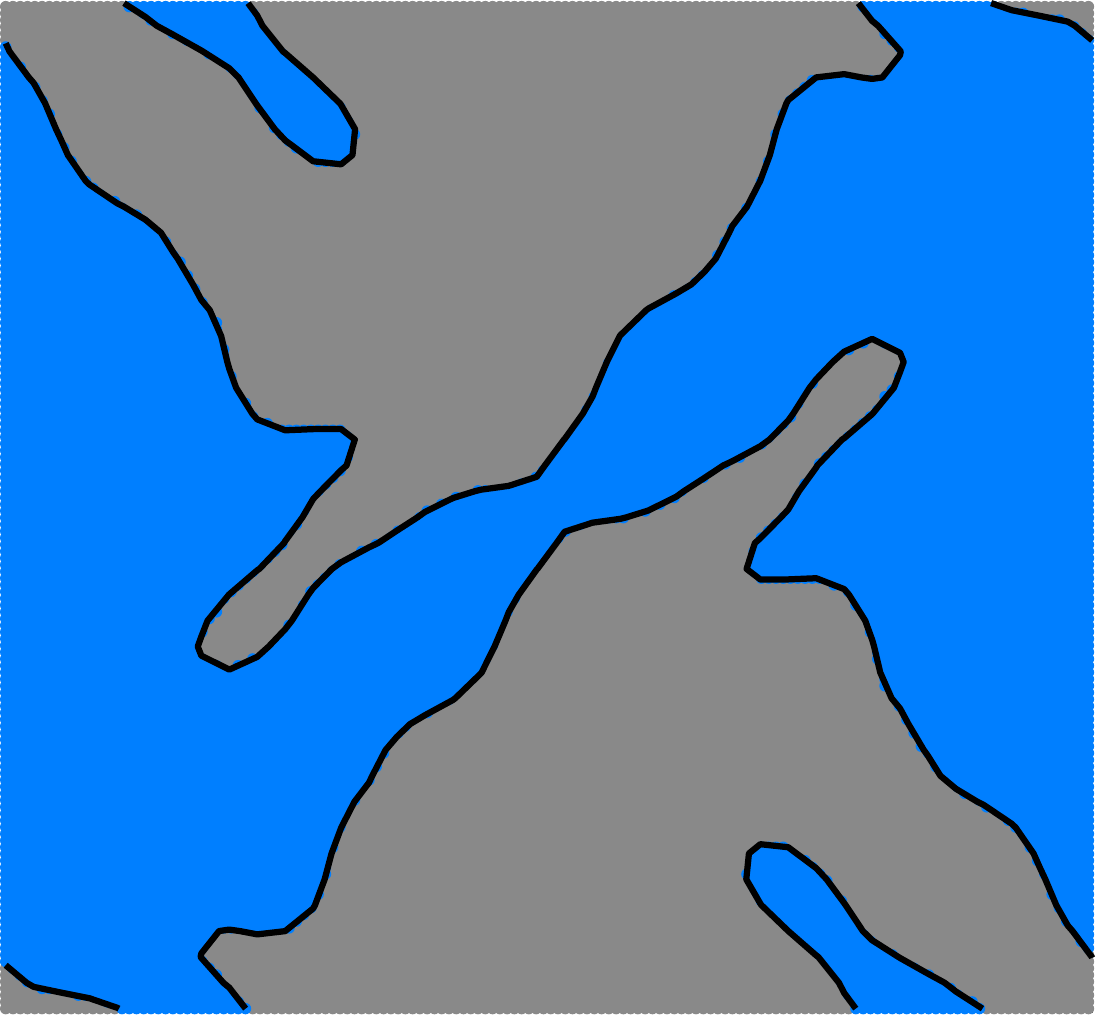}
    \end{subfigure}
    \hfill
    \begin{subfigure}{0.3\textwidth}
        \includegraphics[width=\textwidth]{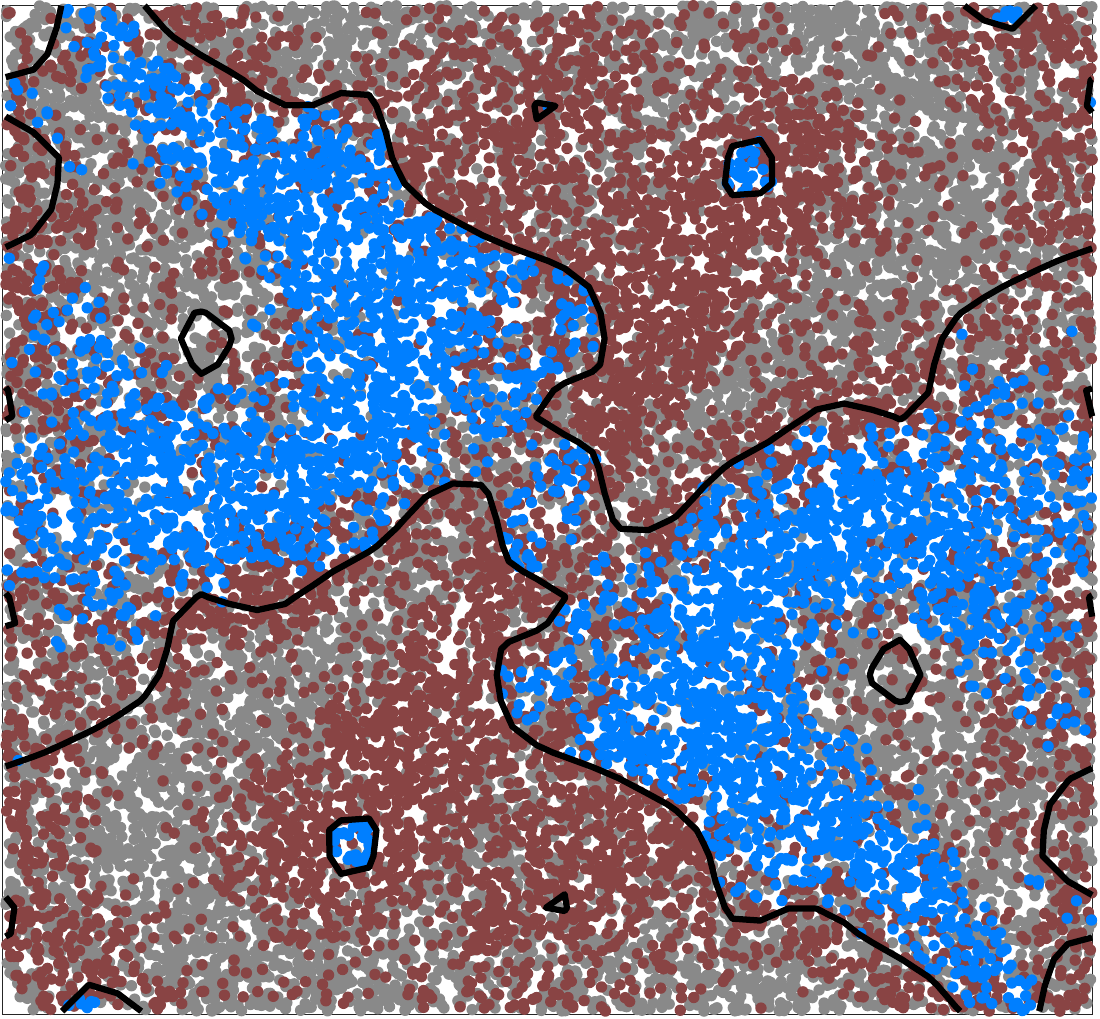}
    \end{subfigure}
    \hfill
    \begin{subfigure}{0.3\textwidth}
        \includegraphics[width=\textwidth]{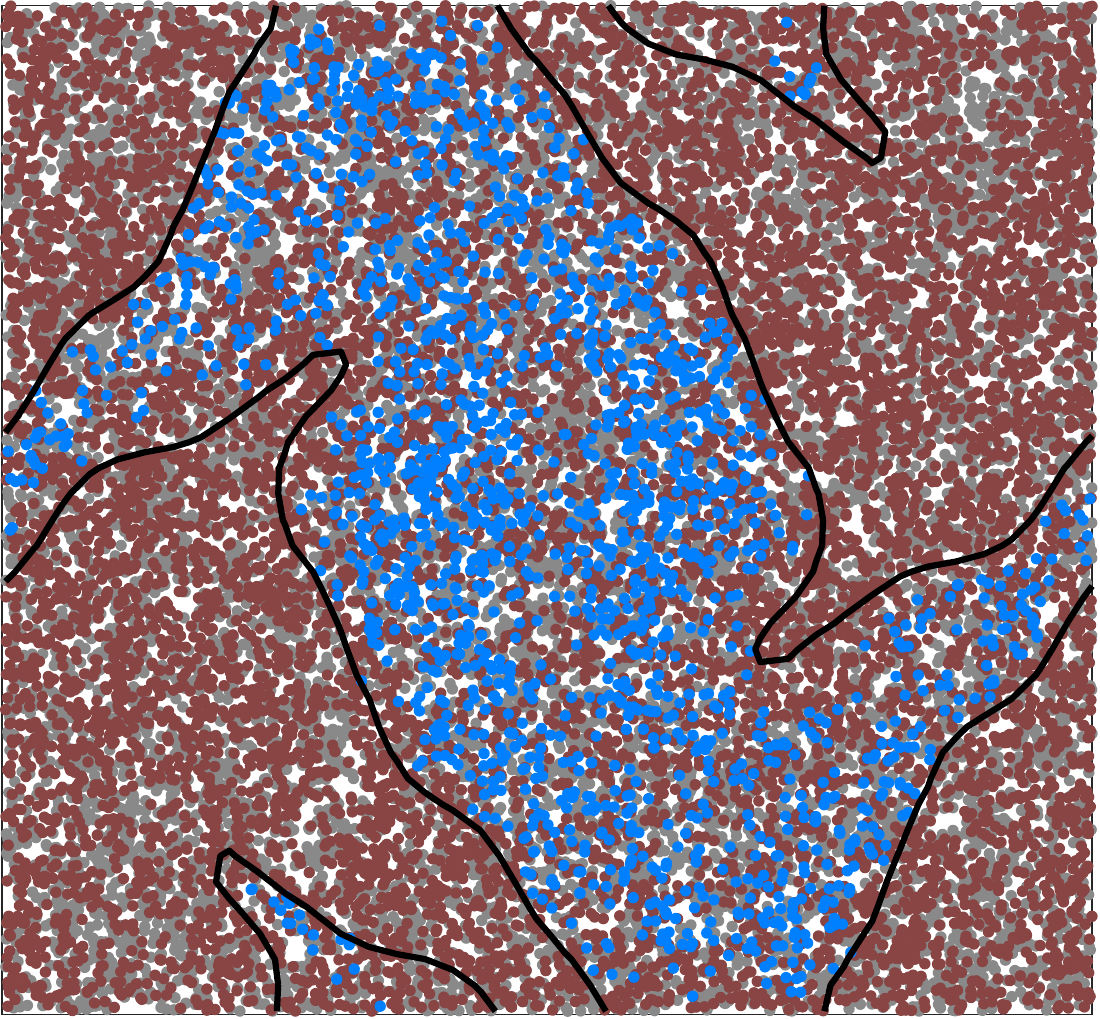}
    \end{subfigure}
    \caption{Shear example: Evolution of the coherent sets, marked in black, determined using method (CS1). The snapshots are taken at time $t=0,5,10$. Particles that started inside $A_{\theta}^0$ and have remained there until time $t$ are colored in blue. Particles that have started inside $A_{\theta}^0$ but were outside the coherent set at some time $0<s\leq t$ are colored in red. Other particles are colored in grey.  A video of the particle simulation for the time interval $[0,10]$ can be found at \url{https://github.com/RobinChemnitz/MatherCoherent}.}
    \label{fig:Shear_snapshots}
\end{figure}

\paragraph{Results.}
In the following we present experimental results for the described system. We choose $\alpha = 0.2 (1, \sqrt{2})$ as the direction of the quasiperiodic driving and $\varepsilon = 0.03$ as the diffusion constant. We define the set $S\in \Z^4$ of Fourier modes like described in \eqref{eq:S_d_S_p}:
\begin{equation*}
    S := \{(m_1, m_2, n_1, n_2)\in \Z^4 \mid \norm{(m_1, m_2)}_\infty \leq K, \: \norm{(n_1, n_2)}_e \leq r\},
\end{equation*}
where $K=6$ and~$r=8$. The constant $K$ determines the resolution in the driving and $r$ determines the resolution in physical space. The chosen values of $K$ and $r$ result in a set $S\subset \Z^4$ with $33293$ elements. We compute the discrete generator~$\Gamma_S$, which is a sparse matrix 
with $0.023\%$ nonzero entries. In Figure~\ref{fig:Shear_spec_sp} we show the spectrum of $\Gamma_S$ around the origin, computed as described in Section~\ref{sec:computational_asp}. Unlike in the previous two examples, there are no recognizable bands of eigenvalues parallel to the imaginary axis. This could indicate the absence of exponential dichotomies for this example. We pick a real eigenvalue $z = -0.097$ of $\Gamma_S$ with the largest nonzero real part. Since $z$ is real, we have $\lambda = -0.097$ and~$\eta = 0$. Let $\widehat{\f}_S$ be the corresponding eigenfunction. As described in Section \ref{sec:numerics_coherent_sets}, $\widehat{\f}_S$ represents a function $\f:= P_S^* \widehat{\f}_S\in L^2(\tor^2\times \tor^2, \C)$. For a given $\theta \in \tor^2$, we compute a coherent set $A_\theta^\bullet$ based on the positive support of the fibres of~$\f$, cf.~method (CS1). Define
\begin{equation*}
    A^t_\theta := \{x\in M \:|\: \text{Re}(e^{\eta t i} \f(\phi^t\theta, x)) > 0\}.
\end{equation*}
Since $\eta = 0$, these sets satisfy $A_{\theta}^t = A_{\phi^t \theta}^0$. Hence, the coherent family can be written as $A_\bullet$. We estimate the survival probability for the starting parameter~$\theta = (0, 0)$. At time $t=0$, we fill the domain $\tor^2$ with a $150\times150$ grid of particles and let them evolve in time according to the SDE \eqref{eq:SDE} with the nonautonomous vectorfield $v$. For each time $t>0$, we obtain an experimental estimate on the survival probability of $A_{\theta}^\bullet$ up to time~$t$, as described in Section~\ref{sec:numerics_coherent_sets}. In Figure~\ref{fig:Shear_snapshots}, we show snapshots of the simulation for three selected time instances. In Figure~\ref{fig:Shear_spec_sp}, there is a plot of the experimental survival probability up to time~$t=10$. From \eqref{eq:eigenfunction_survival_prob}, we expect that the survival probability of $A_{\theta}^\bullet$ decays at most at rate $2\lambda$, up to constants. The results in Figure~\ref{fig:Shear_spec_sp} agree with this estimate as the estimated survival probability decays very similarly to~$e^{2\lambda t}$. The computed cumulative survival probability is~$C(A_{\theta}^\bullet)= 5.680$. This is more than the estimate $-\frac{1}{2 \lambda}= 5.133$ from~\eqref{eq:cum_prob_estimate}. Using method (CS2) or (CS3) did not result in coherent sets with a slowly decaying survival probability.

\section{Conclusion}
Our numerical implementation to compute coherent sets is mainly based on heuristic derivations and experimental evidence. There are two aspects we would like to discuss further. Additionally, we provide an outlook on the application of our methods to more general flows.

\paragraph{Discussion.} The selection of the set of Fourier modes $S$, which is a crucial choice for the construction of the discrete generator $\Gamma_S$, has so far been made on the basis of experimental evidence. In both examples of Section \ref{ssec:osc_gyre} and Section \ref{sec:shear}, we chose $S=S_d \times S_p \subset \Z^{d_d} \times \Z^{d_p}$ with $S_d$ a ball of radius $K$ in the $\ell^\infty$-norm and $S_p$ a ball of radius $r$ in the Euclidean norm. The radii that performed best in these examples, while keeping the computational effort feasible, were $K=6$ and $r=6$. This results in $\abs{S_d} = 169$ and $\abs{S_p}=197$, i.e.~the resolution in the driving and the physical space is similar. The reason we believe that the Euclidean norm works best in the physical space is that the decay of a Fourier mode $F_{m,n}$ due to diffusion scales proportionally to the Euclidean norm of~$n$, cf.~\eqref{eq:def_Gamma}. Hence, Fourier modes $F_{m,n}$, where $n$ has large Euclidean norm decay fast and are thus expected to contribute less to the eigenfunctions of~$\Gamma_S$. In the driving, there is no diffusion, which is a possible explanation why the $\ell^\infty$-norm worked best in our cases. However, these observations have been made only for a small number of examples at a medium resolution $\abs{S} \sim 10^5$. To apply or methods to more complex dynamics at higher resolutions, it is worth to investigate theoretically what choices of $S_d$ and $S_p$ lead to an optimal approximation of the spectrum of $\G$ by the spectrum of~$\Gamma_S$. 

The other aspect we discuss is the extraction of coherent sets from an eigenfunction $\widehat{\f}_S$ of $\Gamma_S$ with eigenvalue~$z$. In the derivation of our method, we treated $\widehat{\f}_S$ as if it was an eigenfunction of $\Gamma$, i.e.~as if $P_S^* \widehat{\f}_S$ was an eigenfunction of $\G$. However, even for a high resolution of $S$, we cannot guarantee that $\smash{\|(\Gamma- z) \widehat{\f}_S\|_2}$ is small. Hence, $\smash{P_S^* \widehat{\f}_S}$ is a priori not even an approximate eigenfunction of~$\G$. The reason for this is that, unlike in the physical space, there is no diffusion in the driving. Due to diffusion in physical space, an eigenfunction $\widehat{\f}_S$ is expected to have low Fourier coefficients $\smash{\widehat{\f}_S(m,n)}$ for large $n$; however, it may have large Fourier coefficients $\widehat{\f}_S(m,n)$ for large values of $m$, but small values of~$n$. Hence, at the boundary of $S$, where $m$ is large, but $n$ is small, there may be a significant discrepancy between $\smash{\Gamma_S \widehat{\f}_S}$ and~$\smash{\Gamma \widehat{\f}_S}$. Numerically our method still performed well, which could indicate that in practice this problem plays a minor role. Finding theoretical bounds for the survival probability of coherent sets extracted from eigenfunctions of the discrete generator $\Gamma_S$ is a challenging task that can be the topic of future research.  

\paragraph{Outlook.} The numerical discretization of the Mather operators $\mathbf{M}^t$ and their generator $\G$ using a Galerkin projection onto Fourier modes required $\Theta = \tor^{d_d}$ and $M = \tor^{d_p}$. For simplicity, we additionally restricted ourselves to the case that the driving dynamics on $\Theta$ is given by a quasiperiodic rotation. In Section~\ref{sec:transfer_operator}, the transfer operator cocycle $\Per$ and the corresponding Mather semigroup have been constructed for more general domains $\Theta$ and $M$ than just tori, and for general ergodic driving. Our results can be viewed as a proof of concept that coherent sets of a nonautonomous dynamics with ergodic driving of the form \eqref{eq:SDE} can be extracted from eigenfunctions of a discretization of the augmented generator $\G$. For domains $\Theta$ and $M$ which are not tori, other discretization methods, like an Ulam discretization, are needed. Hybrid methods are possible as well~\cite{GiDa20}, such as using a Fourier discretization in $\Theta$ and an Ulam discretization in $M$~\cite{froyland2017estimating}. A significant advantage of the Fourier discretization in the case $\Theta=\tor^{d_d}$ and $M = \tor^{d_p}$ is that a finite-dimensional vector $\smash{\widehat{\f}_S \in \C^S}$ represents a function $\f := P_S^*\widehat{\f}_S \in L^2(\Theta \times M, \C)$ in an infinite-dimensional space. Hence, coherent sets $A_\bullet^\bullet$ extracted from $\f$, as described in Section \ref{sec:numerics_coherent_sets}, have a continuous resolution in both $\theta$ and physical space.

Recently, methods have been proposed to faithfully approximate (parts of) the spectrum of transfer-type operators, such as the Koopman operator, see~\cite{GiDa19,korda2020data,giannakis2021delay,colbrook2021rigorous,colbrook2023mpedmd,valva2023consistent}. Since some of these methods explicitly aim at approximating continuous spectra and associated spectral projections, their adaptation to the augmented generator $\G$ could be an interesting future avenue.

\section*{Acknowledgments}

This research has been partially supported by Deutsche Forschungsgemeinschaft (DFG) through grant CRC 1114 ``Scaling Cascades in Complex Systems'', Project Number 235221301, Project A08 ``Characterization and prediction of quasi-stationary atmospheric states''. The authors thank G.~Froyland for engaging discussions and for suggesting the oscillating gyre example in Section~\ref{ssec:osc_gyre}. The third mentioned author thanks I.~Mezi{\'c} for a suggestion regarding literature which proved valuable for this work.

\bibliography{bibliography.bib}
\bibliographystyle{myalpha}

\appendix
\section{Proofs}\label{ap:proofs}

\textbf{Proof of Proposition \ref{prop:spectral_gap}.}
    We use classical arguments to establish the decay of the evolution family. Since the operators $\Per_\theta^t$ are real, it suffices to verify the bound for real functions. We compute the decay of $\norm{f(t)}_2^2$, where $f(t)$ is a solution to the Fokker--Planck equation~\eqref{eq:FP},
    \begin{align*}
        \partial_t \norm{f(t)}_2^2 &= \partial_t \int_M f(t,x)^2 \: dx\\
        &= \int_M 2 f(t,x) \partial_t f(t,x) \: dx \\
        &= \int_M \varepsilon^2f(t,x) \Delta_x f(t, x)  - 2 f(t,x) \nabla_x f(t, x) \cdot v(t, x) \: dx \\
        &= -\varepsilon^2 \int_M \nabla_x f(t, x) \cdot \nabla_x f(t,x) \: dx - \int_M \nabla_x (f(t,x)^2) \cdot v(t,x) \:dx \\
        &= -\varepsilon^2 \norm{\nabla_x f(t)}_2^2 - \int_M \textup{div}_x(f(t,\cdot)^2 v(t,\cdot)) \: dx.
    \end{align*}
    By the divergence theorem, the second integral vanishes. Hence, the decay of $\norm{f(t)}_2^2$ is given by
    \begin{equation}\label{eq:decay_rate}
        \partial_t \norm{f(t)}_2^2 = - \varepsilon^2 \norm{\nabla_x f(t)}_2^2.
    \end{equation}
    The Poincar\'e--Wirtinger inequaility states that there is a constant $c>0$, depending only on $M$, such that for each $f\in \H_0$ we find
    \begin{equation*}
        \norm{\nabla_x f}_2 \geq c \norm{f}_2.
    \end{equation*}
    If $f(0) \in \H_0$, then $f(t)\in \H_0$ for all $t\geq 0$. Hence, we find
    \begin{equation*}
        \partial_t \norm{f(t)}_2^2 \leq - \varepsilon^2 c^2 \norm{f(t)}_2^2.
    \end{equation*}
    By the Grönwall's lemma, we obtain the bound
    \begin{equation*}
        \norm{f(t)}_2^2 \leq \norm{f(0)}_2^2 e^{-\varepsilon^2 c^2 t}
    \end{equation*}
    Lastly, let $f=f(0)\in \H_0$ have $\norm{f}_2 = 1$. By definition of $\Per_\theta^t$, we have $\Per_\theta^t f = f(t)$, and consequently
    \begin{equation*}
        \norm{\Per_\theta^t f}_2 \leq e^{-\frac{\varepsilon^2 c^2}{2} t}.
    \end{equation*}
    This shows that $\norm{\Per_\theta^t}  \leq e^{\varrho t}$ for $\varrho=-\frac{\varepsilon^2 c^2}{2}$, which completes the proof. \hfill $\square$\\

\textbf{Proof of Lemma \ref{lem:uniform_regularity}.}
    Let $f(0) = f \in \H_0$ with $\norm{f}_2 = 1$. The functions $f_i(t) := \Per_{\theta_i}^t f \in \H_0$, for $i=1,2$, are the unique classical solutions to the Fokker--Planck equation with starting parameter $\theta_i$. In particular, we have for $t>0$
    \begin{nalign}\label{eq:generator_f_i}
        \partial_t f_i(t, x) &= \big[G(\phi^t \theta_i) f_i(t)] (x)\\
        &= \frac{1}{2} \varepsilon^2\Delta_x f_i(t, x) - \nabla_x f_i(t, x) \cdot v(\phi^t \theta_i, x).
    \end{nalign}
    We can rewrite the equation \eqref{eq:generator_f_i} for $f_1(t)$, to be
    \begin{nalign}\label{eq:variational}
        \partial_t f_1(t, x) &=  \frac{1}{2} \varepsilon^2 \Delta_x f_1(t, x) - \nabla_x f_1(t, x) \cdot v(\phi^t \theta_2, x) + \nabla_x f_1(t, x) \cdot \big(v(\phi^t \theta_2, x) - v(\phi^t \theta_1, x) \big)\\
        &= \big[ G(\phi^t \theta_2) f_1(t) \big] (x) + g(t),
    \end{nalign}
    where the nonhomogeneous part $g(t)\in \H_0$ is given by
    \begin{equation}\label{eq:def_g}
        g(t, x) := \nabla_x f_1(t, x) \cdot \big(v(\phi^t \theta_2, x) - v(\phi^t \theta_1, x) \big).
    \end{equation}
    The nonautonomous generator $G(\phi^t \theta_2)$ generates the operators $\Per_{\theta_2}^t$. Since $f_1$ is a classical solution to \eqref{eq:variational}, the function $\Per_{\theta_1}^t f = f_1(t)$ is given by the mild solution
    \begin{equation*}
        \big[\Per_{\theta_1}^t f\big] (x) = \big[\Per_{\theta_2}^t f\big] (x) + \int_0^t \big[\Per_{\theta_2}^{t-s} g(s)\big] (x) ds.
    \end{equation*}
    Hence, we can bound the distance between the solutions $\Per_{\theta_1}^t f$ and $\Per_{\theta_2}^t f$.
    \begin{nalign}\label{eq:uniform_regularity_Per}
        \norm{\Per_{\theta_1}^t f - \Per_{\theta_2}^t f}_2^2 &= \int_M \left( \int_0^t [\Per_{\theta_2}^{t-s} g(s)](x) ds\right)^2 dx\\
        &\leq \int_M t \int_0^t [\Per_{\theta_2}^{t-s} g(s)](x)^2\:ds\: dx\\
        &= t \int_0^t \int_M [\Per_{\theta_2}^{t-s} g(s)](x)^2 \:dx\:ds\\
        &= t \int_0^t \norm{\Per_{\theta_2}^{t-s} g(s)}_2^2 \:ds \\
        &< t \int_0^t \norm{g(s)}_2^2\: ds,
    \end{nalign}
    where in the second line we used Jensen's inequality and in the third line the Fubini--Tonelli theorem. We claim that as the distance between $\theta_1$ and $\theta_2$ approaches 0, the integral in the last line tends to 0. By definition of $\delta$ and $g(s)$ in \eqref{eq:def_g}, we find
    \begin{equation*}
        \abs{g(s,x)} \leq \norm{\nabla_xf_1(s,x)}_e \delta,
    \end{equation*}
    and consequently $\norm{g(s)}_2^2 \leq \delta^2 \norm{\nabla_x f_1(s)}_2^2$ for $s\in [0,t]$. In the proof of Proposition \ref{prop:spectral_gap}, namely equation \eqref{eq:decay_rate}, we showed 
    \begin{equation*}
        \norm{\nabla_x f_1(s)}_2^2 = - \frac{1}{\varepsilon^2} \partial_s \norm{f_1(s)}_2^2.
    \end{equation*} 
    This yields
    \begin{equation*}
        \norm{g(s,x)}_2^2 \leq - \frac{\delta^2}{\varepsilon^2} \partial_s \norm{f_1(s)}_2^2.
    \end{equation*}
    Inserting this bound into \eqref{eq:uniform_regularity_Per} yields
    \begin{align*}
        \norm{\Per_{\theta_1}^t f - \Per_{\theta_2}^t f}_2^2 &\leq - \frac{\delta^2 t}{\varepsilon^2} \int_0^t \partial_s \norm{f_1(s)}_2^2 \:ds \\
        &= \frac{\delta^2 t}{\varepsilon^2} \big[\norm{f_1(0)}_2^2 - \norm{f_1(t)}_2^2\big] \\
        &\leq \frac{\delta^2 t}{\varepsilon^2},
    \end{align*}
    where we used $\|f_1(0)\|_2=\|f\|_2 = 1$ and that~$\|f_1(t)\|_2 \le \|f_1(0)\|_2$, because $\Per^t_{\theta}$ is a contraction. This estimate holds for any function $f\in \H_0$ with $\norm{f}_2 = 1$. We conclude
    \begin{equation*}
        \norm{\Per_{\theta_1}^t - \Per_{\theta_2}^t} \leq \frac{\delta \sqrt{t}}{\varepsilon}.
    \end{equation*}
    This completes the proof. \hfill $\square$ \\

\textbf{Proof of Theorem \ref{thm:CohFamFromFun}.}
    Our proof closely follows \cite[A.6]{froyland2017estimating}. Their proof requires that the family of sets $A_\theta^\bullet$ is \emph{sufficiently nice}. We use arguments from \cite{froyland2020computation} to avoid this assumption.

    Firstly, by Proposition \ref{prop:continuous_regularity}, the solution 
    $f(t) = \Per_\theta^t f_0$ is represented by a function $ f \in C(\R_{\geq 0} \times M, \R)$ such that $A_\theta^\bullet$ is well-defined. Fix a time $T\geq 0$ and let $(t_k)_{k\in \N}$ be a dense sequence in $[0, T]$ with $t_1 = T$. Let $x_t(\omega)$ be a random trajectory of the SDE \eqref{eq:SDE}, where $\Omega$ is the probability space on which the SDE is defined. Define the events
    \begin{align*}
        \mathcal{E}_n &:= \{\omega \in \Omega \mid x_{t_k}(\omega) \in A_\theta^{t_k}, \: \forall k\in\{1,\hdots, n\}\}, \\
        \mathcal{E} &:= \{\omega \in \Omega \mid x_t(\omega) \in A_\theta^t, \: \forall t\in [0, T]\}.
    \end{align*} 
    Since $f$ is continuous and $(t_k)_{k\in \N}$ is dense in $[0, T]$, we find
    \begin{equation}
        f(t, x_t(\omega)) \geq 0, \: \forall t\in [0,T] \quad \iff \quad f(t_k, x_{t_k}(\omega)) \geq 0, \: \forall k\in \N.
    \end{equation}
    We conclude $\mathcal{E}_n \downarrow \mathcal{E}$. 
    
    Given a probability measure $\pi$ on $M$, let $\P_{x_0 \sim \pi}$ denote the probability measure on $\Omega$, conditioned on the initial condition $x_0 \sim \pi$. If $\pi$ is a finite measure that is not a probability measure, we define $\P_{x_0 \sim \pi} = \pi(\Omega) \P_{x_0 \sim \pi/\pi(\Omega)}$. If $\pi$ is a finite signed measure, consider the Hahn decomposition $\pi = \pi^+ - \pi^-$, and define $\P_{x_0 \sim \pi} = \P_{x_0 \sim \pi^+} - \P_{x_0 \sim \pi^-}$. Observe that the statement of the theorem does not depend on the scaling of $f_0$. Hence, w.l.o.g.~we may assume $\norm{f_0  }_1 = 2$. In particular, the positive and negative part of $f_0$, which we denote by $f^+_0$, respectively $f^-_0$, integrate to 1. Hence, $f^+_0$ and $f_0^-$ are the densities of probability measures $\nu^+$ and $\nu^-$ on~$M$. In particular, $f_0$ is the density of the signed measure $\nu = \nu^+ - \nu^-$.
    
    Since $\Per_\theta^T f_0$ is the distribution of particles after time $T$ with initial distribution $\nu$, we find
    \begin{equation*}
        \P_{x_0 \sim \nu} (x_T \in A_\theta^T) = \int_{A_\theta^T} [\Per_\theta^T f_0](x) \: dx = \frac{1}{2} \norm{\Per_\theta^T f_0}_1.
    \end{equation*}
    The last equality follows from the definition of $A_\theta^T$ and the fact that $\Per_\theta^T f_0$ integrates to~$0$. For $n\in \N$, consider the decomposition
    \begin{align*}
        \P_{x_0 \sim \nu} (x_T \in A_\theta^T) = &\hspace{3pt} \P_{x_0 \sim \nu} (x_{t_k} \in A_\theta^{t_k}, \: \forall k \in \{1, \hdots, n\}) \\
        &+ \sum_{j = 2}^n \underbrace{\P_{x_0 \sim \nu}(x_{t_j} \notin A_\theta^{t_j}, \: x_{t_k} \in A_\theta^{t_k}, \: \forall t_k > t_j )}_{=:p_j}.
    \end{align*}
    One can show \cite{froyland2013estimating} that $p_j \leq 0$. We conclude
    \begin{equation*}
        \P_{x_0 \sim \nu}(\mathcal{E}_n) \geq \frac{1}{2} \norm{\Per_\theta^T f_0}_1.
    \end{equation*}
    The $\sigma$-additivity of $\P_{x_0 \sim \nu}$ implies $\lim_{n\to \infty}\P_{x_0 \sim \nu}(\mathcal{E}_n) = \P_{x_0 \sim \nu}(\mathcal{E})$. This shows that also $\P_{x_0 \sim \nu}(\mathcal{E})$ is bounded from below by $\frac{1}{2} \norm{\Per_\theta^T f_0}_1$.
    
    It remains to argue that the measure $\nu$ is dominated by $\norm{f}_\infty \abs{A_\theta^0} \1_{A_\theta^0}$, where $\1_{A_\theta^0}$ is the uniform distribution on $A_\theta^0$. This concludes the proof. \hfill $\square$\\

\textbf{Proof of Lemma \ref{lem:projection_real}.}
    Let $\Pi:\Theta \to \L(\H_0)$ be the strongly continuous projection-valued function corresponding to the exponential dichotomy at $\lambda$. We have $S(\theta) = \text{ran}(\Pi(\theta))$ and $U(\theta)=\text{ker}(\Pi(\theta))$.

    Since $f\in U(\theta)$, the inverse $g_t:=\smash{(\Per_{\phi^{-t} \theta}^t)^{-1}f} \in U(\phi^{-t}\theta)$ is well-defined for each $t\geq 0$. In particular, we find $\Per_{\phi^{-t} \theta}^tg_t^R=f^R$ and $\Per_{\phi^{-t} \theta}^tg_t^I=f^I$. Now, assume $f^R \notin U(\theta)$. Consider the decompositions
    \begin{align*}
        f^R &= s_f + u_f, \hspace{6pt}& s_f&:= \Pi(\theta) f^R\in S(\theta), \hspace{6pt}&u_f&:=(\Id - \Pi(\theta)) f^R \in U(\theta),\\
        g^R_t &= s_{g_t} + u_{g_t}, \hspace{6pt}& s_{g_t}&:= \Pi(\phi^{-t}\theta) g^R_t\in S(\phi^{-t}\theta), \hspace{6pt}&u_{g_t}&:=(\Id - \Pi(\phi^{-t}\theta)) g^R_t \in U(\phi^{-t}\theta),
    \end{align*}
    in particular $s_f\neq 0$ by our assumption.
    By definition of $\Pi$, we have $\Per_{\phi^{-t}\theta}^t s_f = s_{g_t}$, $\Per_{\phi^{-t}\theta}^t u_f = u_{g_t}$ and
    \begin{align*}
        \norm{s_f}_2 &= \norm{\Per_{\phi^{-t}\theta}^t s_{g_t}}_2 \leq C e^{(\lambda-\beta)t}\norm{s_{g_t}}_2, \\
        \norm{u_{g_t}}_2 &= \norm{(\Per_{\phi^{-t}\theta}^t)^{-1} u_f}_2 \leq C e^{-(\lambda+\beta)t}\norm{u_f}_2 .
    \end{align*}
    This shows that $\norm{s_{g_t}}_2$ grows at least at rate $-(\lambda-\beta)$ in $t$, while $u_{g_t}$ grows at most at rate $-(\lambda + \beta)$ . We conclude
    \begin{equation}
        \label{eq:PerBound}
    \begin{aligned}
        \norm{(\Per_{\phi^{-t} \theta}^t)^{-1} f}_2 &\geq \norm{(\Per_{\phi^{-t} \theta}^t)^{-1} f^R}_2 = \norm{g_t^R}_2 \geq \norm{s_{g_t}}_2 - \norm{u_{g_t}}_2 \\
        &\geq C^{-1} e^{-(\lambda -\beta) t} \norm{s_f}_2 - Ce^{-(\lambda + \beta)t}\norm{u_f}_2.
    \end{aligned}
    \end{equation}
    Definition~\ref{def:exp_dichotomy} implies that $\| (\Per^t_{\theta}\vert_U)^{-1}\| \le Ce^{-(\lambda+\beta) t}$ for $t\ge 0$ \emph{uniformly} in~$\theta\in\Theta$. Since $f\in U(\theta)$, it follows that $e^{(\lambda+\beta) t} \| (\Per^t_{\phi^{-t}\theta})^{-1} f\|_2$ is bounded by the constant $C\|f\|_2$ for every~$t\ge 0$. This contradicts \eqref{eq:PerBound}, which predicts that $\smash{ e^{(\lambda+\beta) t} \| (\Per^t_{\phi^{-t}\theta})^{-1} f\|_2 }$ grows to infinity at least with rate~$2\beta > 0$, neglecting additive constants.
    We conclude $s_f = 0$ and, therefore, $f^R\in U(\theta)$. The proof for $f^I \in U(\theta)$ is analogous. \hfill $\square$ \\

\textbf{Proof of Lemma \ref{lem:div_free_fourier}.}
    The first equality is standard and holds for any real-valued function. We express the divergence of $v$ is terms of its Fourier coefficients.
    \begin{align}
        \text{div}_x(v)(\theta, x) &= \sum_{j=1}^{d_p} \partial_{x_j} v_j(\theta, x) \\
        &= \sum_{j=1}^{d_p}\partial_{x_j} \sum_{m,n} F_{m,n}(\theta, x) \widehat{v_j}(m,n) \\
        &= -2\pi i \sum_{j=1}^{d_p} \sum_{m,n} n_j F_{m,n}(\theta, x) \widehat{v_j}(m,n) \\
        &= -2\pi i \sum_{m,n} (n \cdot \widehat{v}(m,n)) F_{m,n}(\theta, x).
    \end{align}
    This shows that the Fourier coefficients of the divergence of $v$ are given by $n \cdot \widehat{v}(m,n)$. Since the divergence of $v$ is zero as a function, all the Fourier coefficients of the divergence must vanish. We conclude $n \cdot \widehat{v}(m,n) = 0$ for all $m,n\in \Z^{d_d} \times \Z^{d_p}$. \hfill $\square$ \\

\textbf{Proof of Proposition \ref{prop:commutative_diagram}.}
    Fix a fixed starting parameter $\theta \in \tor^2$ and a starting distribution $f\in \H_0$. Define $f_1(t) := (T_{\phi^t\theta} \circ \Per_{\theta}^t) f$ and $f_2(t) := (\Per_\text{eff}^t \circ T_\theta) f$. Both functions satisfy $f_1(0) = f_2(0) = T_\theta f$. We compute the time-derivatives of both functions using the Fokker--Planck equation \eqref{eq:FP}
    \begin{align*}
        \partial_t f_1(t, x) &= \partial_t \big\lbrack [\Per_\theta^t f](x - \theta - \alpha t)\big \rbrack \nonumber\\
        &= \frac{1}{2} \varepsilon^2 \Delta_x [\Per_\theta^t f](t, x-\phi^t\theta) - \nabla_x [\Per_\theta^t f](t,x - \phi^t\theta) \cdot \big(v(\phi^t\theta, x - \phi^t \theta)  + \alpha\big)\\
        &= \frac{1}{2} \varepsilon^2 \Delta_x f_1(t, x) - \nabla_x f_1(t,x) \cdot \big(v_\text{aut}(x)  + \alpha\big).\nonumber\\     
        \partial_t f_2(t, x) &= \partial_t [\Per_\text{eff}^t (T_\theta f)](x) \nonumber\\
        &=\frac{1}{2} \varepsilon^2 \Delta_x [\Per_\text{eff}^t (T_\theta f)](x) - \nabla_x [\Per_\text{eff}^t (T_\theta f)](x) \cdot v_\text{eff}(x)\\
        &= \frac{1}{2} \varepsilon^2 \Delta_x f_2(t, x) - \nabla_x f_2(t,x) \cdot \big(v_\text{aut}(x) + \alpha \big).\nonumber
    \end{align*}
    Since $f_1$ and $f_2$ solve the same NACP, which is known to have a unique solution, we conclude $f_1(t) = f_2(t)$ for all $t\geq 0$. This concludes the proof. \hfill $\square$\\

\textbf{Proof of Proposition \ref{prop:semi_aut_SS}.} 
    Consider $\lambda \notin \Sigma(\Per_\text{eff})$, such that there is a strongly continuous projection-valued function $\Pi_\text{eff}:\tor^2 \to \L(\H)$ that defines an exponential dichotomy of $\Per_\text{eff}$ at $\lambda$. By Section \ref{sec:autonomous_case}, the function $\Pi_\text{eff}$ is constant, i.e.~does not depend on $\theta$. Define another projection-valued function $\Pi:\tor^2 \to \L(\H_0)$ by $\Pi(\theta) f = T_{-\theta} \circ \Pi_\text{eff} \circ T_{\theta}$. We show that $\Pi$ defines an exponential dichotomy of $\Per$. Let $\theta \in \Theta$ and $t\geq 0$. Using the commutativity of \eqref{eq:commutative_diagram}, we compute
    \begin{align*}
        \Per_\theta^t \Pi(\theta) &= T_{\phi^t \theta}^{-1} \circ \Per_\text{eff}^t \circ T_\theta \circ T_{-\theta}\circ  \Pi_\text{eff} \circ T_\theta \\
        &= T_{-\phi^t \theta} \circ \Pi_\text{eff} \circ \Per_{\text{eff}}^t \circ T_\theta \\
        &= T_{-\phi^t \theta} \circ \Pi_\text{eff}\circ T_{\phi^t \theta} \circ T_{-\phi^t \theta} \circ \Per_{\text{eff}}^t \circ T_\theta\\
        &= \Pi(\phi^t \theta) \Per_\theta^t.
    \end{align*}
    Let $S(\theta)$ and $U(\theta)$ be the range and kernel of $\Pi(\theta)$. Analogously, let $S_\text{eff}$ and $U_\text{eff}$ be the range and kernel of $\Pi_\text{eff}$. By definition of $\Pi$, we find $U(\theta) = T_{-\theta} U_\text{eff}$ and $S(\theta) = T_{-\theta} S_\text{eff}$. Hence the restriction $\Per_\theta^t |_U : U(\theta) \to U(\phi^t \theta)$ is given by
    \begin{equation}\label{eq:Per_eff_U}
        \Per_\theta^t |_U = (T_{\phi^t \theta}^{-1} \circ \Per_\text{eff}^t \circ T_\theta) |_U = T_{\phi^t \theta}^{-1} \circ \Per_\text{eff}^t|_{U_\text{eff}} \circ T_\theta.
    \end{equation}
    Since $\Per_\text{eff}^t|_{U_\text{eff}}$ is invertible, this shows that $\Per_\theta^t |_U$ is invertible. 
    
    Let $C>0$ and $\beta>0$ be the constants of the exponential dichotomy $\Pi_\text{eff}$. Consider $f\in S(\theta)$, such that $T_\theta f\in S_\text{eff}$. We find
    \begin{equation*}
        \norm{\Per_\theta^t f}_2 = \norm{T_{\phi^t \theta}^{-1} \circ \Per_\theta^t \circ T_\theta f}_2 \leq C e^{(\lambda - \beta) t} \norm{T_\theta f}_2.
    \end{equation*}
    Since $T_\theta$ is an isometry, this proves $\norm{\Per_\theta^t|_S}\leq C e^{(\lambda - \beta) t}$. Now, consider $f \in U(\phi^t \theta)$, such that $T_{\phi^t\theta} f \in U_\text{eff}$. From \eqref{eq:Per_eff_U}, we derive
    \begin{equation}
        \norm{(\Per_\theta^t)^{-1} f}_2 = \norm{T_\theta^{-1} \circ (\Per_\text{eff}^t)^{-1} \circ T_{\phi^t \theta } f}_2 \geq C e^{-(\lambda + \beta) t} \norm{ T_{\phi^t \theta} f}_2.
    \end{equation}
    Since $T_{\phi^t \theta}$ is an isometry, this proves $\norm{\Per_\theta^t|_U}^{-1}\geq C^{-1} e^{(\lambda + \beta) t}$. This shows that $\Pi$ defines an exponential dichotomy of $\Per$ at $\lambda$. \hfill $\square$

\end{document}